\documentclass{amsart}
\usepackage{amsmath,amssymb,amscd,amsthm}
\newtheorem{formula}{}[subsection]
\newtheorem{proposition}[formula]{Proposition}
\newtheorem{corollary}[formula]{Corollary}
\newtheorem{lemma}[formula]{Lemma}
\newtheorem{theorem}[formula]{Theorem}
\newtheorem{problem}[formula]{Problem}
\theoremstyle{definition}
\newtheorem{definition}[formula]{Definition}
\newtheorem{example}[formula]{Example}
\newtheorem{construction}[formula]{Construction}
\theoremstyle{remark}
\newtheorem*{remark}{Remark}

\renewcommand{\b}{\beta}
\newcommand{\D}{\Delta}
\newcommand{\e}{\varepsilon}
\renewcommand{\t}{\theta}
\renewcommand{\L}{\Lambda}
\renewcommand{\l}{\lambda}
\newcommand{\bl}{\lambda\kern-0.53em\lambda}
\newcommand{\bmu}{\mu\kern-0.55em\mu}
\newcommand{\bnu}{\nu\kern-0.51em\nu}
\newcommand{\f}{\varphi}
\renewcommand{\O}{\Omega}

\newcommand{\A}{\mathcal A}
\newcommand{\F}{\mathcal F}

\renewcommand{\P}{\mathcal P}

\newcommand{\M}{{\mathrm M}}

\newcommand{\C}{\mathbb C}
\newcommand{\N}{\mathbb N}
\newcommand{\Q}{\mathbb Q}
\newcommand{\R}{\mathbb R}
\newcommand{\Z}{\mathbb Z}

\renewcommand{\k}{\mathbf k}

\newcommand{\<}{\langle}
\renewcommand{\>}{\rangle}

\newcommand{\mb}[1]{\textbf {\textit#1}}
\newcommand{\sbr}[2]{{\textstyle\genfrac{[}{]}{}{}{#1}{#2}}}
\newcommand{\bin}[2]{{\textstyle\binom{#1}{#2}}}

\newcommand{\cs}{\mathbin{\#}}

\newcommand{\bs}{\mathop{\rm bs}}
\newcommand{\bc}{\mathop{\rm bc}}

\newcommand{\cone}{\mathop{\rm cone}}
\newcommand{\core}{\mathop{\rm core}}
\newcommand{\link}{\mathop{\rm link}}
\renewcommand{\star}{\mathop{\rm star}}

\newcommand{\cc}{\mathop{\rm cc}}
\newcommand{\cub}{\mathop{\rm cub}}

\newcommand{\ind}{\mathop{\rm ind}\nolimits}

\newcommand{\bideg}{\mathop{\rm bideg}}
\newcommand{\id}{\mathrm{id}}
\newcommand{\hd}{\mathop{\rm hd}}
\newcommand{\Kd}{\mathop{\rm Kd}}
\newcommand{\Ker}{\mathop{\rm Ker}}
\newcommand{\Tor}{\mathop{\rm Tor}\nolimits}

\newcommand{\zp}{\mathcal Z_P}
\newcommand{\ma}{\mathop{\rm ma}}
\newcommand{\zk}{\mathcal Z_K}
\newcommand{\wk}{\mathcal W_K}
\newcommand{\bk}{B_T\zk}
\newcommand{\bp}{B_T\zp}

\newcommand{\sign}{\mathop{\rm sign}}
\newcommand{\td}{\mathop{\rm td}}

\begin{document}

\title{Torus actions, combinatorial topology and homological algebra}
\author{Victor M. Buchstaber}
\author{Taras E. Panov}
\thanks{Partially supported by the Russian Foundation for Fundamental
Research, grant no. 99-01-00090}
\address{Department of Mathematics and Mechanics, Moscow
State University, 119899 Moscow RUSSIA}
\email{buchstab@mech.math.msu.su}
\email{tpanov@mech.math.msu.su}

\begin{abstract}
The paper surveys some new results and open problems connected with such
fundamental combinatorial concepts as polytopes, simplicial complexes,
cubical complexes, and subspace arrangements. Particular attention is paid to
the case of simplicial and cubical subdivisions of manifolds and, especially,
spheres. We describe important constructions which allow to study all these
combinatorial objects by means of methods of commutative and homological
algebra. The proposed approach to combinatorial problems relies on the theory
of moment-angle complexes, currently being developed by the authors. The
theory centres around the construction that assigns to each simplicial
complex $K$ with $m$ vertices a $T^m$-space $\zk$ with a special bigraded
cellular decomposition.  In the framework of this theory, the well-known
non-singular toric varieties arise as orbit spaces of maximally free actions
of subtori on moment-angle complexes corresponding to simplicial spheres. We
express different invariants of simplicial complexes and related
combinatorial-geometrical objects in terms of the bigraded cohomology rings
of the corresponding moment-angle complexes. Finally, we show that the new
relationships between combinatorics, geometry and topology result in
solutions to some well-known topological problems.
\end{abstract}

\maketitle

\tableofcontents

\section*{Introduction}
We survey the results and open problems in the vast field, which was being
shaped during the last two decades and incorporated various aspects of
combinatorics of polytopes, combinatorial and algebraic topology, homological
algebra, group actions on topological spaces, algebraic geometry of toric
varieties and symplectic geometry. The main aim of this review is to show
that the theory of moment-angle complexes proposed by the authors allows to
substantially extend the relationships between the above listed branches of
mathematical science and thereby obtain solutions to some well-known
problems. Each section of the survey refers to a separate subject and
contains the necessary introductory remarks. Below we schematically overview
the contents of the article; this can be considered as the guide to the whole
survey.

Chapter 1 contains the necessary combinatorial, geometrical and topological
facts about polytopes, simplicial and cubical complexes and manifolds. We
describe both classical and original constructions, which allow to study
the combinatorial objects by methods of commutative and homological
algebra.

Section 1 of chapter 2 is the overview of constructions and results on
algebraic toric varieties and related combinatorial objects, which are
necessary for the rest of the paper. In section~2.2 we describe the
topological analogues of toric varieties, the quasitoric manifolds introduced
by Davis and Januszkiewicz. The main construction and results about the
quasitoric manifolds can be also found there. In section~2.3 we present the
solution to the quasitoric analogue of the Hirzebruch problem about connected
representatives in the cobordism classes of stably complex manifolds, which
was recently obtained by Buchstaber and Ray. In section~2.4 we give the
obtained by Panov combinatorial formulae for Hirzebruch genera of quasitoric
manifolds. The last section of chapter~2 describes the known results on the
classification of toric and quasitoric manifolds over a given simple
polytope.

The theory of moment-angle complexes being developed by the authors is the
subject of chapters~3 and~4. To each simplicial complex $K$ with $m$
vertices there is assigned the moment-angle complex~$\zk$ (see
section~3.2). The complex $\zk$ caries a canonical action of the
torus~$T^m$ with quotient the cone over $K$ and possesses the canonical
bigraded cellular decomposition (see section~3.3). A series of important
results of the theory arises from the fact that $\zk$ is a manifold provided
that~$K$ is a simplicial sphere. At the same time in the more general case of
simplicial manifold $K$ the singular points of~$\zk$ form an orbit of the
torus action, and the complement of the neighbourhood of this orbit is a
manifold with boundary. Using the bigraded cellular structure, in chapter~4
we calculate the bigraded cohomology ring of the moment-angle
complex~$\zk$. This calculation reveals new relationships with some
well-known constructions from homological algebra and opens the way to
solution of many combinatorial problems.

In chapter 5 we apply the theory of moment-angle complexes to the
well-known problem of calculation the cohomology of subspace arrangement
complements. We concentrate on coordinate subspace arrangements and
diagonal subspace arrangements (sections~5.2 and~5.3 correspondingly) and
survey the results obtained by the authors in this direction. In
particular, we calculate the cohomology ring of the complement of a
coordinate subspace arrangement by reducing to the cohomology of the
moment-angle complex~$\zk$. We also reduce the problem of calculating the
cohomology of the diagonal subspace arrangement complement to calculating
the cohomology of the loop space~$\O\zk$.

Almost all new concepts in our survey are provided by the corresponding
examples. We also give a lot of examples of particular computations, which
illustrate the general theorems. All the results in our paper are of three
types. Firstly, interpretations of classical results, which sometimes are
given without references. Secondly, results of other authors, which are
always provided by the corresponding references. And thirdly, results that
either have been obtained recently by the present authors, or are
extensions of results from the authors papers \cite{BP1}--\cite{BP5},
\cite{Pan1},~\cite{Pan2} and papers \cite{BR1},~\cite{BR2} by N.~Ray and
the first author.

The authors wish to express special thanks to Levan~Alania, Yusuf~Civan,
Natalia~Dobrinskaya, Nikolai~Dolbilin, Mikhail~Farber, Ivan~Izmestiev,
Oleg~Musin, Sergey~Novikov, Nigel~Ray, Elmer~Rees, Mikhail~Shtan'ko,
Mikhail~Shtogrin, Vladimir~Smirnov, Neil~Strickland, Sergey~Tarasov,
Victor~Vassiliev, Volkmar~Wel\-ker, Sergey~Yuzvinsky, G\"unter~Ziegler for
the insight gained from numerous discussions of different questions
relevant to the survey. We also grateful to all participants of the
seminar ``Topology and computational geometry" being held by O.\,R.~Musin
and the authors at the Department of Mathematics and Mechanics, Moscow
State University.

\section{Algebraic, combinatorial and geometrical background}
\subsection{Polytopes}
\label{poly}
Both combinatorial and geometrical aspects of the theory of convex polytopes
are exposed in a vast number of textbooks, monographs and papers. We just
mention Ziegler's book~\cite{Zi}, where a host of further references can be
found. In this section we review some basic concepts and constructions
used in the rest of the paper.

There are two different ways to define a convex polytope in $n$-dimensional
affine space~$\R^n$.

\begin{definition}
\label{pol1}
A {\it convex polytope\/} is the convex hull of a finite set of points in
some~$\R^n$.
\end{definition}

\begin{definition}
\label{pol2}
A {\it convex polyhedron\/} $P$ is an intersection of finitely many
half-spaces in some~$\R^n$:
\begin{equation}
\label{ptope}
  P=\bigl\{\mb x\in\R^n\::\:\<\mb l_i,\mb x\>\ge-a_i,
  \; i=1,\ldots,m\bigr\},
\end{equation}
where $\mb l_i\in(\R^n)^*$, $i=1,\ldots,m$, are some linear functions and
$a_i\in\R$. A~{\it (convex) polytope\/} is a bounded convex polyhedron.
\end{definition}
Nevertheless, the above two definitions produce the same geometrical object,
i.e. the subset of $\R^n$ is a convex hull of a finite point set if and only
if it is a bounded intersection of finitely many half-spaces.
This fact is proved in many textbooks on polytopes and convex geometry, see
e.g.~\cite[Theorem~1.1]{Zi}.

The {\it dimension\/} of a polytope is the dimension of its affine hull;
without loss of generality we may consider only $n$-dimensional polytopes
$P^n$ in $n$-dimensional space~$\R^n$. A {\it supporting hyperplane\/} of
$P^n$ is an affine hyperplane $H$ which intersects $P^n$ and for which the
polytope is contained in one of the two closed half-spaces determined by the
hyperplane. The intersection $P^n\cap H$ is a {\it face\/} of the polytope. We
also regard the polytope $P^n$ itself as a face; other faces are called {\it
proper faces\/}. The {\it boundary\/} $\partial P^n$ is the union of all
proper faces of~$P^n$. Each face of $n$-dimensional polytope $P^n$ is itself
a polytope of dimension~$\le n$. 0-dimensional faces are called {\it
vertices\/}, 1-dimensional faces are called {\it edges\/}, and codimension
one faces are called {\it facets\/}. The faces of $P^n$ of all dimensions
form a partially ordered set ({\it poset\/}) with respect to the inclusion.
This poset is called the {\it face lattice\/} of~$P^n$.

Two polytopes $P_1\in\R^{n_1}$ and $P_2\in\R^{n_2}$ of the same dimension are
said to be {\it affinely equivalent\/} if there is an affine map
$\R^{n_1}\to\R^{n_2}$ taking one polytope to another. Two polytopes are {\it
combinatorially equivalent\/} if there is a bijection between their faces
that preserves the inclusion relation. In other words, two polytopes are
combinatorially equivalent if their face lattices are isomorphic as posets.

\begin{example}[simplex and cube]
\label{simcub}
  An $n$-dimensional {\it simplex\/} $\D^n$ is the convex hull of $(n+1)$
  points of $\R^n$ that do not lie on a common affine hyperplane. All faces
  of an $n$-simplex are simplices of dimension~$\le n$. All $n$-simplices are
  affinely equivalent. The {\it standard\/} $n$-simplex is the convex hull of
  points $(1,0,\ldots,0)$, $(0,1,\ldots,0),\ldots,(0,\ldots,0,1)$, and
  $(0,\ldots,0)$ in $\R^n$. Alternatively, the standard $n$-simplex is
  defined by $(n+1)$ inequalities
  \begin{equation}
  \label{stsim}
    x_i\ge0,\;i=1,\ldots,n,\quad\text{and}\quad-x_1-\ldots-x_n\ge-1.
  \end{equation}
  The {\it regular\/} $n$-simplex is the convex hull of $n+1$
  points $(1,0,\ldots,0)$, $(0,1,\ldots,0)$, $\ldots$, $(0,\ldots,0,1)$
  in~$\R^{n+1}$.

  The {\it standard
  $q$-cube\/} is the convex polytope $I^q\subset\R^q$ defined by
  \begin{equation}
  \label{cube}
    I^q=\{(y_1,\ldots,y_q)\in\R^q\,:\;0\le y_i\le1,\,i=1,\ldots,q\}.
  \end{equation}
  Alternatively, the standard $q$-cube is the convex hull of $2^q$ points in
  $\R^q$ having only zero and unit coordinates.
\end{example}

The following construction identifies a convex
$n$-polytope with $m$ facets as the intersection of the positive cone
\begin{equation}
\label{pcone}
  \R^m_+=\bigl\{(y_1,\ldots,y_m)\in\R^m\::\:y_i\ge0,\;
  i=1,\ldots,m\bigr\}\subset\R^m
\end{equation}
with a certain $n$-dimensional plane.
\begin{construction}
\label{dist}
  Let $P\in\R^n$ be a convex $n$-polytope given by~(\ref{ptope}) with some
  $\mb l_i\in(\R^n)^*$, $a_i\in\R$, $i=1,\ldots,m$. Introduce $n\times
  m$-matrix $L$ whose columns are vectors $\mb l_i$ written in the standard
  basis of $(\R^n)^*$, i.e. $(L)_{ji}=(\mb l_i)_j$. Note that $L$ is of
  rank~$n$. Likewise, let $\mb a=(a_1,\ldots,a_m)^t\in\R^m$ be the column
  vector with entries~$a_i$. Then~(\ref{ptope}) shows that
  \begin{equation}
  \label{pdiag}
    P=\bigl\{\mb x\in\R^n\::\:(L^t\mb x+\mb a)_i\ge0,
    \; i=1,\ldots,m\bigr\},
  \end{equation}
  where $L^t$ is the transposed matrix and $\mb x=(x_1,\ldots,x_n)^t$ is the
  column vector. Consider the affine map
  \begin{equation}
  \label{af}
    A_P:\R^n\to\R^m,\quad A_P(\mb x)=L^t\mb x+\mb a\in\R^m.
  \end{equation}
  Its image is an $n$-dimensional plane in $\R^m$, and~(\ref{pdiag}) shows
  that $A_P(P)$ is the intersection of this plane with the positive
  cone~$\R_+^m$.
  Now, let $W$ be an $m\times(m-n)$-matrix of rank $(m-n)$ such that
  $LW=0$. Then it is easy to see that
  $$
    A_P(P)=\bigl\{\mb y\in\R^m\::\:W^t\mb y=W^t\mb a,\;y_i\ge0,\quad
    i=1,\ldots,m\bigr\}.
  $$
Note that the polytopes $P$ and $A_P(P)$ are affinely equivalent.
\end{construction}

\begin{example}\label{simdist}
Consider the standard $n$-simplex $\D^n\subset\R^n$ defined by
inequalities~(\ref{stsim}).
It has $m=n+1$ facets, and $\mb
l_1=(1,0,\ldots,0)^t$, $\ldots$, $\mb l_n=(0,\ldots,0,1)^t$,
$\mb l_{n+1}=(-1,\ldots,-1)^t$, $a_1=\ldots=a_n=0$, $a_{n+1}=1$
(see (\ref{ptope})).
One can take $W=(1,\ldots,1)^t$ in Construction~\ref{dist}. Hence, $W^t\mb
y=y_1+\ldots+y_m$, $W^t\mb a=1$, and we have
$$
  A_{\D^n}(\D^n)=\bigl\{\mb y\in\R^{n+1}\::\:y_1+\ldots+y_{n+1}=1,\;
  y_i\ge0,\quad i=1,\ldots,n\bigr\}.
$$
This is the regular $n$-simplex in $\R^{n+1}$.
\end{example}

\begin{remark}
Convex polytopes introduced in definitions~\ref{pol1} and~\ref{pol2} are {\it
geometrical\/} objects. However, one could also consider {\it combinatorial\/}
polytopes which are classes of combinatorially equivalent polytopes. In fact,
a combinatorial polytope is the face lattice (regarded as a poset) of a
(geometrical) polytope. In this review we deal with both geometrical and
combinatorial polytopes.
\end{remark}

Two different definitions of a convex polytope lead to two
different notions of {\it generic\/} polytopes.

A set of $m>n$ points in $\R^n$ is in {\it general position\/} if no $(n+1)$
of them lie on a common affine hyperplane. From the viewpoint of
Definition~\ref{pol1}, the convex polytope is generic if it is the convex
hull of a set of points in general position. This implies that all proper
faces of the polytope are simplices, i.e. every facet has the minimal number
of vertices (namely,~$n$). Such polytopes are called {\it simplicial\/}.

A set of $m>n$ hyperplanes $\<\mb l_i,\mb x\>=-a_i$, $\mb l_i\in(\R^n)^{*}$,
$\mb x\in\R^n$, $a_i\in\R$, $i=1,\ldots,m$, is in {\it general position\/} if
no point lies in more than $n$ of them. From the viewpoint of
Definition~\ref{pol2}, the convex polytope $P^n$ is generic if its bounding
hyperplanes (see~(\ref{ptope})) are in general position, i.e. there exactly
$n$ facets meet at each vertex of~$P^n$. Such polytopes are called {\it
simple\/}.

For any convex polytope $P\subset\R^n$ define its {\it polar set\/}
$P^*\subset(\R^n)^*$ as
$$
  P^*=\{\mb x'\in(\R^n)^{*}\::\:\;\langle\mb x',\mb x\rangle\ge-1
  \:\text{ for all }\:\mb x\in P\}.
$$
\begin{remark}
Our definition of the polar set is that used in the algebraic geometry of
toric varieties, not the classical one from the convex geometry. The latter
is obtained by replacing the inequality ``$\ge-1$" above by ``$\le1$".
One polar set is taken into another by the symmetry with centre in~0.
\end{remark}

It is well known in convex geometry that the polar set $P^*$ is a convex
polytope in the dual space $(\R^n)^*$ and $(P^*)^*=P$ provided that $0\in P$.
The polytope $P^*$ is called the {\it polar\/} (or {\it dual\/}) of~$P$.
Moreover, there is a one-to-one order reversing correspondence between face
lattices of $P$ and~$P^*$. In particular, if $P$ is simple, then $P^*$ is
simplicial, and vice versa.

\begin{example}
  Any polygon (2-polytope) is simple and simplicial at the same time.
  In dimensions $\ge3$ only polytope that is simultaneously simple and
  simplicial is the simplex. The cube is a simple polytope.
  The polar of a simplex is again a simplex. The polar of the cube is called
  the {\it cross-polytope\/}.
  The 3-dimensional cross-polytope is the octahedron.
\end{example}

In the sequel simple polytopes are denoted by the Latin letters $P$,
$Q$ etc., while simplicial ones are denoted by the same letters with
asterisque: $P^{*}$, $Q^{*}$ etc.

Each face of a simple polytope is a simple polytope. The {\it
product\/} $P_1\times P_2$ of two simple polytopes $P_1$ and $P_2$ is a simple
polytope as well. Let $P_1^{*}$ and $P_2^*$ be the corresponding polar
simplicial polytopes. Then $P_1^*\circ P_2^*:=(P_1\times P_2)^*$ is again a
simplicial polytope. The operation $\circ$ on simplicial polytopes can be
directly described as follows. Realise $P_1^*$ in $\R^{n_1}$ and $P_2^*$ in
$\R^{n_2}$ in such way that $0\in P^*_1$ and $0\in P^*_2$.
Then $P_1^*\circ P_2^*\subset\R^{n_1}\times\R^{n_2}$ is the convex hull of
the union of $P_1^*\subset\R^{n_1}\times0$ and $P_2^*\subset0\times\R^{n_2}$.

\begin{construction}[Connected sum of two simple polytopes]
\label{consum}
Suppose we are given two simple polytopes $P^n$ and $Q^n$, both of dimension
$n$, with distinguished vertices $v$ and $w$ respectively. The informal way
to get the {\it connected sum\/} $P^n\cs_{v,w} Q^n$ of $P^n$ at $v$ and $Q^n$
at $w$ is as follows. We ``cut off" $v$ from $P^n$ and $w$ from $Q^n$; then,
after a projective transformation, we can ``glue" the rest of $P^n$ to the
rest of $Q^n$ along the new facets to obtain $P^n\cs_{v,w} Q^n$.  Below we
give the formal definition, following~\cite[\S6]{BR2}; this definition will
be used later.

First, we introduce an $n$-polyhedron $\varGamma^n$, which will be used as a
template for the construction; it arises by considering the standard
$(n-1)$-simplex $\D^{n-1}$ in the subspace $\{x:x_1=0\}$ of $\R^n$, and
taking its cartesian product with the first coordinate axis. The facets $G_r$
of $\varGamma^n$ therefore have the form $\R\times D_r$,
where $D_r$, $1\leq r\leq n$, are the faces of~$\D^{n-1}$.
Both $\varGamma^n$ and the $G_r$ are divided into positive and negative
halves, determined by the sign of the coordinate~$x_1$.

We order the facets of $P^n$ meeting in $v$ as $E_1,\ldots,E_n$, and the
facets of $Q^n$ meeting in $w$ as $F_1,\ldots,F_n$.  Denote the complementary
sets of facets by $\mathcal{C}_v$ and $\mathcal{C}_w$; those in
$\mathcal{C}_v$ avoid $v$, and those in $\mathcal{C}_w$ avoid~$w$.

We now choose projective transformations $\phi_P$ and $\phi_Q$ of~$\R^n$,
whose purpose is to map $v$ and $w$ to $x_1=\pm\infty$ respectively.  We
insist that $\phi_P$ embeds $P^n$ in $\varGamma^n$ so as to satisfy two
conditions; firstly, that the hyperplane defining $E_r$ is identified with
the hyperplane defining~$G_r$, for each $1\leq r\leq n$, and secondly, that
the images of the hyperplanes defining $\mathcal{C}_v$ meet $\varGamma^n$ in
its negative half. Similarly, $\phi_Q$ identifies the hyperplane defining
$F_r$ with that defining~$G_r$, for each $1\leq r\leq n$, but the images of
the hyperplanes defining $\mathcal{C}_w$ meet $\varGamma^n$ in its positive
half.  We define the {\it connected sum\/} $P^n\cs_{v,w} Q^n$ of $P^n$ at $v$
and $Q^n$ at $w$ to be the simple convex $n$-polytope determined by the
images of the hyperplanes defining $\mathcal{C}_v$ and $\mathcal{C}_w$ and
hyperplanes defining $G_r$, $r=1,\ldots,n$. It is defined only up to
combinatorial equivalence; moreover, different choices for either of $v$
and~$w$, or either of the orderings for $E_r$ and~$F_r$, are likely to affect
the combinatorial type. When the choices are clear, or their effect on the
result irrelevant, we use the abbreviation $P^n\cs Q^n$.

The related construction of connected sum $P\cs Q^*$ of a simple polytope $P$
and a simplicial polytope $Q^*$ is described in~\cite[Example~8.41]{Zi}.
\end{construction}

A simplicial polytope $P^*$ is called {\it $k$-neighbourly\/} if any $k$
vertices span a face of~$P^*$. Likewise, a simple polytope $P$ is called
$k$-neighbourly if any $k$ facets of $P$ have non-empty intersection (i.e.
share a common codimension-$k$ face). Obviously, every simplicial (or simple)
polytope is 1-neighbourly. It can be shown (\cite[Corollary~14.5]{Br}, see
also Example~\ref{hneib} below) that if $P^*$ is a $k$-neighbourly simplicial
$n$-polytope and $k>\sbr n2$, then $P^*$ is an $n$-simplex. In particular,
any 2-neighbourly simplicial 3-polytope is a simplex. However, there exist
simplicial $n$-polytopes with any number of vertices which are $\sbr
n2$-neighbourly. Such polytopes are called {\it neighbourly\/}. In particular,
there exist a simplicial 4-polytope (different from 4-simplex) any two
vertexes of which are connected by an edge.

\begin{example}[neighbourly 4-polytope]
Let $P=\D^2\times\D^2$ be the product of two triangles. Then $P$ is a simple
polytope, and it is easy to see that any two facets of $P$ share a common
2-face. Hence, $P$ is 2-neighbourly. The polar $P^*$ is a neighbourly
simplicial 4-polytope.
\end{example}

More generally, it is easy to see that if a simple polytope $P_1$ is
$k_1$-neighbourly and a simple polytope $P_2$ is $k_2$-neighbourly, then
the product $P_1\times P_2$ is a $\min(k_1,k_2)$-neighbourly simple polytope.
It follows that $(\D^n\times\D^n)^*$ and $(\D^n\times\D^{n+1})^*$ provide
examples of neighbourly simplicial $2n$- and $(2n+1)$-polytopes. Another
example of neighbourly polytopes, with arbitrary number of vertices, is as
follows.

\begin{example}[cyclic polytopes]
\label{cyclic}
Define the {\it moment curve\/} in $\R^n$ as
$$
  \mb x:\R\longrightarrow\R^n,\qquad t\mapsto\mb x(t)=(t,t^2,\ldots,t^n)
  \in\R^n.
$$
For any $m>n$ define
the {\it cyclic polytope\/} $C^n(t_1,\ldots,t_m)$ as the convex hull of $m$
distinct points $\mb x(t_i)$, $t_1<t_2<\ldots<t_m$, on the moment curve. It
follows from the Vandermonde determinant identity that no $n+1$ points
on the moment curve lie on a common affine hyperplane. Hence,
$C^n(t_1,\ldots,t_m)$ is a simplicial $n$-polytope. It can be shown
(see~\cite[Theorem~0.7]{Zi}) that $C^n(t_1,\ldots,t_m)$ has exactly $m$
vertices $\mb x(t_i)$, the combinatorial type of cyclic polytope does not
depend on the specific choice of the parameters $t_1,\ldots,t_m$, and
$C^n(t_1,\ldots,t_m)$ is a neighbourly simplicial $n$-polytope. We denote the
combinatorial cyclic $n$-polytope with $m$ vertices by~$C^n(m)$.
\end{example}

Let $P$ be a simple $n$-polytope. Denote by $f_i$ the number of faces of $P$
of codimension $(i+1)$ (i.e. of dimension $(n-i-1)$). In particular, $f_0$
is the number of facets of $P$, which we will denote $m(P)$ or just~$m$.
\begin{definition}
\label{fhvect}
The integer vector $\mb f(P)=(f_0,\ldots,f_{n-1})$ is called the {\it
$f$-vector\/} of simple polytope~$P$. The integer vector
$(h_0,h_1,\ldots,h_n)$ whose components $h_i$ are defined from the equation
\begin{equation}
\label{hvector}
  h_0t^n+\ldots+h_{n-1}t+h_n=(t-1)^n+f_0(t-1)^{n-1}+\ldots+f_{n-1}
\end{equation}
is called the {\it $h$-vector\/} of~$P$. Finally, the integer vector
$(g_0,g_1,\ldots,g_{\sbr n2})$, where $g_0=1$, $g_i=h_i-h_{i-1}$, $i>0$, is
called the {\it $g$-vector\/} of~$P$.
\end{definition}
We also put $f_{-1}=1$, which means that the polytope itself is a face of
codimension~0. The $f$-vector and the $h$-vector determine each other by
means of linear relations, namely
\begin{equation}
\label{hf}
  h_k=\sum_{i=0}^k(-1)^{k-i}\bin{n-i}{n-k}f_{i-1},\quad
  f_{n-1-k}=\sum_{q=k}^n\bin qk h_{n-q},\quad k=0,\ldots,n.
\end{equation}
In particular, $h_0=1$ and
$h_n=(-1)^n\bigl(1-f_0+f_1+\ldots+(-1)^nf_{n-1}\bigr)$. By Euler's theorem,
\begin{equation}\label{euler}
  f_0-f_1+\cdots+(-1)^{n-1}f_{n-1}=1+(-1)^{n-1},
\end{equation}
which is equivalent to $h_n=1=h_0$.  In the case of simple polytopes Euler's
theorem admits the following generalisation.

\begin{theorem}[Dehn--Sommerville relations]
\label{ds}
  The $h$-vector of any simple $n$-polytope is symmetric, i.e.
  $$
    h_i=h_{n-i},\quad i=0,1,\ldots,n.
  $$
\end{theorem}
There are a lot of different ways to prove the Dehn--Sommerville equation. We
present a proof which uses a Morse-theoretical argument, firstly appeared
in~\cite{Br}. We will return to this argument in chapter~2.
\begin{proof}[Proof of Theorem~{\rm \ref{ds}}.] Let $P^n\subset\R^n$ be a
simple
polytope. Choose a linear function $\varphi:\R^n\to\R$ which is generic,
that is, it takes different values at all vertices of~$P^n$. For this
$\varphi$ there is a vector $\bnu$ in $\R^n$ such that $\varphi(\mb
x)=\langle\bnu,\mb x\rangle$. Note that $\bnu$ is parallel to no edge
of~$P^n$.  Now we view $\varphi$ as a height function on~$P^n$. Using
$\varphi$, we make the 1-skeleton of $P^n$ a directed graph by orienting each
edge in such a way that $\varphi$ increases along it (this can be done since
$\varphi$ is generic). For each vertex $v$ of $P^n$ define its index,
$\ind(v)$, as the number of incident edges that point towards~$v$. Denote the
number of vertices of index~$i$ by $I_\nu(i)$.  We claim that
$I_\nu(i)=h_{n-i}$. Indeed, each face of $P^n$ has a unique top vertex (a
maximum of the height function $\varphi$ restricted to the face) and a unique
bottom vertex (the minimum of $\varphi$). Let $F^k$ be a $k$-face of~$P^n$,
and $v_F$ its top vertex. Since $P^n$ is simple, there exactly $k$ edges of
$F^k$ meet at $v_F$, whence $\ind(v_F)\ge k$. On the other hand, each vertex
of index $q\ge k$ is the top vertex for exactly $\bin qk$ faces of
dimension~$k$. It follows that $f_{n-1-k}$ (the number of $k$-faces) can be
calculated as
$$
  f_{n-1-k}=\sum_{q\ge k}\bin qk I_\nu(q).
$$
Now, the second
identity from (\ref{hf}) shows that $I_\nu(q)=h_{n-q}$, as claimed. In
particular, the number $I_\nu(q)$ does not depend on~$\bnu$. On the other
hand, since $\ind_\nu(v)=n-\ind_{-\nu}(v)$ for any vertex $v$, one has
$$
  h_{n-q}=I_\nu(q)=I_{-\nu}(n-q)=h_q.
$$
\end{proof}
Using~(\ref{hf}), we can rewrite the Dehn--Sommerville equations
in terms of the $f$-vector as follows
\begin{equation}
\label{DSf}
  f_{k-1}=\sum_{j=k}^n(-1)^{n-j}\bin jk f_{j-1},\quad k=0,1,\ldots,n.
\end{equation}
The Dehn--Sommerville equations were established by Dehn for
$n\le5$ in 1905, and by Sommerville in the general case in 1927
(see~\cite{So}) in the form similar to~(\ref{DSf}).

\begin{example}
Let $P_1^{n_1}$ and $P_2^{n_2}$ be simple polytopes. Any face of
$P_1\times P_2$ is the product of a face of $P_1$ and a face of $P_2$, whence
$$
  f_k(P_1\times P_2)=
  \sum_{i=-1}^{n_1-1}f_i(P_1)f_{k-i-1}(P_2), \quad k=-1,0,\ldots,n_1+n_2-1.
$$
Set $h(P;t)=h_0+h_1t+\dots+h_nt^n$. Then it follows from the above formula
and~(\ref{hvector}) that
\begin{equation}\label{hvprod}
  h(P_1\times P_2;t)=h(P_1;t)h(P_2;t).
\end{equation}
\end{example}

\begin{example}
Let us express the $f$-vector and the $h$-vector of the connected sum $P^n\cs
Q^n$ in terms of that of $P^n$ and~$Q^n$. It follows from
Construction~\ref{consum} that
\begin{align*}
  f_i(P^n\cs Q^n)&=f_i(P^n)+f_i(Q^n)-\bin n{i+1},\quad i=0,1,\ldots,n-2;\\
  f_{n-1}(P^n\cs Q^n)&=f_{n-1}(P^n)+f_{n-1}(Q^n)-2
\end{align*}
(note that $\bin n{i+1}=f_i(\D^{n-1})$). Then it follows from (\ref{hf})
that
\begin{align*}
  h_0(P^n\cs Q^n)&=h_n(P^n\cs Q^n)=1;\\
  h_i(P^n\cs Q^n)&=h_i(P^n)+h_i(Q^n),\quad i=1,2,\ldots,n-1.
\end{align*}
\end{example}

Thus, $h_i$, $i=1,\ldots,n-1$, define integer-valued functions on the set of
simple polytopes which are linear with respect to the connected sum operation.

\begin{problem}
Describe all integer-valued functions on the set of
simple polytopes which are linear with respect to the connected sum operation.
\end{problem}

The $f$-vector of a {\it simplicial\/} polytope $P^*$ is defined as
$\mb f(P^*)=(f_0,f_1,\ldots,f_{n-1})$, where $f_i$ is the number of $i$-faces
($i$-simplices) of~$P^*$. The $h$-vector $h(P^*)=(h_0,h_1,\ldots,h_n)$ is
determined by identity~(\ref{hvector}). Note that if $P^*$ is the simplicial
polytope polar to a simple polytope $P$, then $f_i(P^*)=f_i(P)$. In
particular, the Dehn--Sommerville equations hold for simplicial polytopes as
well.

\begin{example}
\label{hneib}
Suppose $P^*$ is a $q$-neighbourly simplicial $n$-polytope. Then
$f_{k-1}(P^*)=\bin mk$, $k\le q$. From~(\ref{hf}) we obtain
\begin{equation}\label{fhneib}
  h_k(P^*)=\sum_{i=0}^k(-1)^{k-i}\bin{n-i}{k-i}\bin mi=
  \bin{m-n+k-1}k,\quad k\le q.
\end{equation}
(The last identity is obtained by calculating the coefficient of $t^k$ from
two sides of $\frac1{(1+t)^{n-k+1}}(1+t)^m=(1+t)^{m-n+k-1}$.) Suppose that
$P^*$ is not a simplex. Then $m>n+1$, which together with~(\ref{fhneib})
gives $h_0<h_1<\dots<h_q$. It follows from the Dehn--Sommerville
equations that $q\le\sbr n2$.
\end{example}

A natural question arises: which integer vectors may appear as the
$f$-vectors of simple (or, equivalently, simplicial) polytopes? The
Dehn--Sommerville relations provide a necessary condition.

\begin{proposition}[\cite{Kl}]
The Dehn--Sommerville relations are the most general linear equations
satisfied by the $f$-vectors of all simple (or simplicial) polytopes.
\end{proposition}
\begin{proof}
In \cite{Kl} the statement was proved directly, using~$f$-vectors. Instead,
we use $h$-vectors, which somewhat simplifies the proof. It is sufficient to
prove that the affine hull of the $h$-vectors $(h_0,h_1,\ldots,h_n)$ of
simple $n$-polytopes is an $\sbr n2$-dimensional plane. Set
$Q_k:=\D^k\times\D^{n-k}$, $k=0,1\ldots,\sbr n2$. Since $\mb
h(\D^k)=1+t+\dots+t^k$, the formula~(\ref{hvprod}) gives
$$
  \mb h(Q_k)=\frac{1-t^{k+1}}{1-t}\cdot\frac{1-t^{n-k+1}}{1-t}.
$$
It follows that $\mb h(Q_{k+1})-\mb h(Q_k)=t^{k+1}+\dots$,
$k=0,1,\ldots,\sbr n2-1$. Therefore, the vectors $\mb h(Q_k)$,
$k=0,1\ldots,\sbr n2$, are affinely independent.
\end{proof}

We mention also that the identity~(\ref{euler}) is the only linear relation
satisfied by the vectors of face numbers of {\it general\/} convex polytopes.

The conditions characterising the $f$-vectors of simple (or simplicial)
polytopes, now know as the {\it $g$-theorem\/}, were conjectured by
McMullen~\cite{McM1} in 1970 and proved by Stanley~\cite{St2}
(necessity) and Billera, Lee~\cite{BL} (sufficiency) in 1980. Besides the
Dehn--Sommerville equations, the $g$-theorem contains two groups of
inequalities, one linear and one non-linear. To
state the $g$-theorem completely, we need the following construction. For
any two positive integers $a$, $i$ there exists a unique {\it binomial
$i$-expansion\/} of $a$ of the form
$$
  a=\bin{a_i}i+\bin{a_{i-1}}{i-1}+\dots+\bin{a_j}j,
$$
where $a_i>a_{i-1}>\dots>a_j\ge j\ge1$. Define
$$
  a^{\langle i\rangle}=
  \bin{a_i+1}{i+1}+\bin{a_{i-1}+1}i+\dots+\bin{a_j+1}{j+1},\quad
  0^{\langle i\rangle}=0.
$$

\begin{theorem}[$g$-theorem]
\label{gth}
An integer vector $(f_0,f_1,\ldots,f_{n-1})$ is the $f$-vector of a simple
$n$-polytope if and only if the corresponding sequence $(h_0,\ldots,h_n)$
determined by~{\rm(\ref{hvector})} satisfies the following three conditions:
\begin{enumerate}
\item[(a)] $h_i=h_{n-i}$, $i=0,\ldots,n$ (the Dehn--Sommerville equations);
\item[(b)] $h_0\le h_1\le\dots\le h_{\sbr n2}$, i.e. $g_i\ge0$,
  $i=0,1,\ldots,\sbr n2$.
\item[(c)] $h_0=1$, $h_{i+1}-h_i\le(h_i-h_{i-1})^{\langle i\rangle}$, i.e.
$g_{i+1}\le g_i^{\langle i\rangle}$, $i=1,\ldots,\sbr n2-1$.
\end{enumerate}
\end{theorem}
\begin{remark}
Obviously, the same conditions characterise the $f$-vectors of {\it
simplicial\/} polytopes.
\end{remark}

An integral sequence $(k_0,k_1,\ldots,k_r)$ satisfying $k_0=1$
and $0\le k_{i+1}\le k_i^{\langle i\rangle}$ for $i=1,\ldots,r-1$ is called
an {\it $M$-vector\/} (after M.~Macaulay). Conditions~(b) and~(c) from the
$g$-theorem imply that the $g$-vector $(g_0,g_1,\ldots,g_{\sbr n2})$ of a
simple $n$-polytope is an $M$-vector. On the other hand, the notion of
$M$-vector appears in the following classification result of commutative
algebra.
\begin{theorem}[Macaulay]
\label{mvect}
An integral sequence $(k_0,k_1,\ldots,k_r)$ is an $M$-vector if and only if
there exists a commutative graded algebra $A=A^0\oplus A^2\oplus\dots\oplus
A^{2r}$ over a field $\k=A^0$, generated (as an algebra) by degree-two
elements, such that the dimension of $2i$-th graded component of $A$ equals
$k_i$, i.e. $\dim_\k A^{2i}=k_i$, $i=1,\ldots,r$.
\end{theorem}
\noindent The proof can be found in~\cite{St1}.

To prove the sufficiency of the $g$-theorem, Billera and Lee presented a
remarkable combinatorial-geometrical construction of a simplicial polytope
with any prescribed $M$-sequence as its $g$-vector. Stanley's proof of the
necessity of $g$-theorem (i.e. that the $g$-vector of a simple polytope is
an $M$-vector) relies upon deep results from algebraic geometry: the
Hard Lefschetz theorem for the cohomology of toric varieties.
We give the ideas of Stanley's proof in section~\ref{tori}.
After 1995 several more elementary combinatorial proofs of the $g$-theorem
appeared. The first elementary proof by McMullen~\cite{McM2}, which uses the
{\it polytope algebra\/} instead of the cohomology algebra of toric variety,
was still very involved. Recently Timorin~\cite{Ti} found much more simple
elementary proof of the $g$-theorem, which relies on the interpretation of
McMullen's polytope algebra as the algebra of differential operators (with
constant coefficients) vanishing on the {\it volume polynomial\/} of the
polytope.

It follows from the results
of~\cite{BL} that the part~(b) of the $g$-theorem, i.e. the inequalities
\begin{equation}
\label{glbt}
  h_0\le h_1\le\dots\le h_{\sbr n2},
\end{equation}
give the most general linear
inequalities satisfied by the $f$-vectors of simple (or simplicial)
polytopes. These inequalities are now known as the {\it Generalised Lower
Bound Theorem\/} (GLBT) for simple (simplicial) polytopes. During the last
two decades a lot of work was done and progress achieved in extending the
Dehn--Sommerville equations, GLBT, and $g$-theorem to more general objects
than simplicial polytopes. However, there are still a lot of intriguing open
problems here. Some of them are presented in the survey article by
Stanley~\cite{St5}. In the present paper we also review some related
questions (see the comments in the next section).

The $g$-theorem has the following important corollary.
\begin{theorem}[{Upper Bound Theorem (UBT) for simplicial polytopes,
\cite[Theorem~8.23, Corollary~8.38]{Zi}}]
\label{ubt}
  From all simplicial $n$-polytopes $P^*$ with $m$ vertices the cyclic
  polytope $C^n(m)$ (Example~{\rm \ref{cyclic}}) has the maximal number of
  $i$-faces, $2\le i\le n-1$. That is, if $f_0(P^*)=m$, then
  $$
    f_i(P^*)\le f_i\bigl(C^n(m)\bigr), \qquad 2\le i\le n-1.
  $$
\end{theorem}
\noindent Note that $f_i(C^n(m))=\bin m{i+1}$ for $0\le i<\sbr n2$.
The above theorem was conjectured by Motzkin in 1957. It was proved by
McMullen in 1970 and motivated him to conjecture the $g$-theorem. McMullen
showed also that the Upper Bound theorem is equivalent to the following
inequalities for the $h$-vector $\mb h(P^*)=(h_0,h_1,\ldots,h_n)$:
$$
  h_i(P^*)\le\bin{m-n+i-1}i,\qquad 0\le i\le\sbr n2
$$
(compare this with Example~\ref{hneib}).

To conclude these remarks about polytopes, we introduce an important
algebraic invariant of a (combinatorial) simple polytope, which appears
many times throughout this review. Let $P$ be a simple $n$-polytope with $m$
facets $F_1,\ldots,F_m$. Fix a commutative ring $\k$ with unit.
Let $\k[v_1,\ldots,v_m]$ be the polynomial algebra over $\k$ on $m$
generators. We make it a graded algebra by setting $\deg(v_i)=2$.
\begin{definition}
\label{frpol}
The {\it face ring\/} (or the {\it Stanley--Reisner ring\/}) of a
simple polytope $P$ is the quotient ring
$$
  \k(P)=\k[v_1,\ldots,v_m]/\mathcal I_P,
$$
where $\mathcal I_P$ is the ideal generated by all square-free
monomials $v_{i_1}v_{i_2}\cdots v_{i_s}$, $i_1<\dots<i_s$, such that
$F_{i_1}\cap\cdots\cap F_{i_s}=\emptyset$ in~$P$.
\end{definition}
\noindent Obviously, $\k(P)$ is a graded $\k$-algebra.

We mention that the Stanley--Reisner ring, the $f$-vector, and the
$h$-vector are invariants of a {\it combinatorial\/} simple polytope: they
depend only on the face lattice and do not depend on a particular
geometrical realisation.

\subsection{Simplicial complexes: topology and combinatorics}
\label{sim1}
Let $[m]$ denote the index set $\{1,\ldots,m\}$. For any subset $I\subset[m]$
denote by $\#I$ its cardinality.
\begin{definition}
\label{absimcom}
An (abstract) {\it simplicial complex\/} on the set $[m]$ is a
collection $K=\{I\}$ of subsets of $[m]$ such that for each
$I=\{i_1,\ldots,i_k\}\in K$ all subsets of $I$ (including $\emptyset$) also
belong to~$K$. Subsets $I\in K$ are called (abstract) {\it simplices\/} of~$K$.
\end{definition}
\noindent Similarly one defines a simplicial complex on any set~$\mathcal S$.
One-element subsets from $K$ are called {\it vertices\/} of~$K$.
If $K$ contains all one-element subsets of~$[m]$, then we say that $K$ is a
simplicial complex {\it on the vertex set~$[m]$\/}.
The {\it dimension\/} of an abstract simplex $I=\{i_1,\ldots,i_k\}\in K$ is
its cardinality minus one, i.e. $\dim I=\#I-1$. The dimension of an abstract
simplicial complex is the maximal dimension of its simplices.

\begin{definition}
\label{polyhed}
A {\it geometrical simplicial complex\/} (or {\it polyhedron\/}) is
a subset $\mathcal P\subset\R^n$ represented as a union of simplices of
any dimensions in such a way that the intersection of any two simplices
is a face of each. (By simplices here we mean
convex polytopes defined in Example~\ref{simcub}.)
\end{definition}

In the sequel we denote by $\D^{m-1}$ both the abstract simplex
(i.e. the simplicial complex consisting of {\it all\/} subsets of~$[m]$)
and the corresponding polyhedron.

\begin{remark}
The notion of polyhedron from Definition~\ref{pol2} is not the
same as that from the above definition. The first meaning of the term
``polyhedron" (i.e. the ``unbounded polytope") is adopted in the convex
geometry, while the second one (i.e. the ``geometrical simplicial
complex") is used in the combinatorial topology. Since both meanings became
classical in the appropriate science, we do not change their names. In this
review, by ``polyhedron" we will usually mean a geometrical simplicial
complex. The ``unbounded polytope" will be referred to as ``convex
polyhedron". Anyway, it will be always clear from the context which
``polyhedron" is under consideration.
\end{remark}

Only finite geometrical simplicial complexes (polyhedrons) are considered
in this review. The dimension of a polyhedron is the maximal dimension of its
simplices. It is a classical fact~\cite{Pon} that any $n$-dimensional abstract
simplicial complex $K$ admits a {\it geometrical realisation\/} $|K|$ as an
$n$-dimensional polyhedron in $\R^{2n+1}$
(abstract simplices of $K$ correspond to (polytopal) simplices of~$|K|$). A
geometrical realisation $|K|$ of an abstract simplicial complex $K$ is unique
up to a piecewise-linear homeomorphism.

\begin{construction}
One can construct a geometric realisation of a simplicial complex $K$ on
the vertex set $[m]$ in $m$-dimensional space as follows. Let $\mb e_i$ denote
the $i$-th unit coordinate vector in~$\R^m$. For each subset $I\subset[m]$
denote by $\D_I$ the convex hull of vectors $\mb e_i$ with $i\in I$. Obviously,
$\D_I$ is a (geometric) simplex. Then one has
$$
  |K|=\bigcup_{I\in K}\D_I\subset\R^m.
$$
\end{construction}

A {\it simplicial map\/} $f:|K_1|\to|K_2|$ of two polyhedrons is any mapping
of the set of vertices of $|K_1|$ to the vertices of~$|K_2|$, extended
linearly on simplices of $|K_1|$ to the whole of~$|K_1|$. A polyhedron $|K'|$
is called a {\it subdivision\/} of polyhedron $|K|$ if each simplex of $|K|$
is a union of finitely many simplices of~$|K'|$. A {\it piecewise linear
($PL$) map\/} $f:|K_1|\to|K_2|$ is a map that is simplicial between some
subdivisions of $|K_1|$ and~$|K_2|$. The standard reference for the $PL$
topology is~\cite{RoSa}.

\begin{example}[associated simplicial complex]
\label{dual}
Let $K$ be a simplicial complex on the set~$[m]$. Suppose that $K$ is not the
$(m-1)$-simplex. Define
$$
  \widehat{K}:=\bigl\{I\subset[m]\::\:[m]\setminus I\notin K\bigr\}.
$$
Obviously, $\widehat{K}$ is a simplicial complex on~$[m]$. It is called the
{\it simplicial complex associated with\/}~$K$.
\end{example}

\begin{construction}[join of simplicial complexes]
\label{join}
Let $K_1$, $K_2$ be simplicial complexes on $[m_1]$ and $[m_2]$ respectively.
Identify $[m_1]\cup[m_2]$ with $[m_1+m_2]$. The
{\it join\/} of $K_1$ and $K_2$ is the simplicial complex
$$
  K_1*K_2:=\bigl\{ I\subset[m_1+m_2]\::\:I=I_1\cup I_2,\;
  I_1\in K_1, I_2\in K_2\bigr\}
$$
on the set $[m_1+m_2]$.
\end{construction}

\begin{example}
1. If $K_1=\D^{m_1-1}$, $K_2=\D^{m_2-1}$, then $K_1*K_2=\D^{m_1+m_2-1}$.

2. The simplicial complex $\D^0*K$ (the join of $K$ and a point) is called
the {\it cone\/} over $K$ and denoted $\cone(K)$.

3. Let $S^0$ be disjoint union of two vertices. Then
$S^0*K$ is called the {\it suspension\/} over $K$ and denoted $\Sigma K$.
\end{example}
The geometric realisation of $\cone(K)$ (of $\Sigma K$) is the topological
cone (suspension) over~$|K|$.

The {\it barycentric subdivision\/} of an abstract simplicial complex $K$ is
the simplicial complex $\bs(K)$ on the set $\{I, I\in K\}$
of simplices of $K$ such that $\{I_1,\ldots,I_r\}\in\bs(K)$ if and only
if $I_1\subset I_2\subset\dots\subset I_r$ (after possible re-ordering).
The {\it barycentre\/} of a (polytopal) simplex $\D^n\in\R^n$ with vertices
$v_1,\ldots,v_{n+1}$ is the point $\bc(\D^n)=\frac
1{n+1}(v_1+\dots+v_{n+1})\in\D^n$. The {\it barycentric subdivision\/} of a
polyhedron $\mathcal P$ is the polyhedron $\bs(\mathcal P)$ defined as
follows. The vertices of $\bs(\mathcal P)$ are barycentres of simplices of
$\mathcal P$ of all dimensions. The set of vertices
$\{\bc(\D_1^{i_1}),\ldots,\bc(\D_r^{i_r})\}$ spans a simplex of
$\bs(\mathcal P)$ if and only if $\D_1^{i_1}\subset\ldots\subset\D_r^{i_r}$
in~$\mathcal P$. Obviously, one has $|\bs(K)|=\bs(|K|)$
for any abstract simplicial complex~$K$.

\begin{example}[order complex of a poset]
\label{oc}
Let $(\mathcal S,\prec)$ be any poset. Define $K(\mathcal S)$ to be the set
of all chains $x_1\prec x_2\prec\cdots\prec x_k$, $x_i\in\mathcal S$.
It is easy to see that $K(\mathcal S)$ is a simplicial complex. This complex
is called the {\it order complex\/} of the poset $(\mathcal S,\prec)$. In
particular, if $(\mathcal S,\varsubsetneq)$ is the face poset (face lattice)
of a simplicial complex~$K$, then $K(\mathcal S)$ is the barycentric
subdivision of~$K$.
\end{example}

The {\it missing face\/} of simplicial complex $K$ on the set $[m]$ is a
subset $I\subset[m]$ such that $I\notin K$, but every proper subset of $I$ is
a simplex of $K$. A {\it flag complex\/} is a simplicial complex for which
every missing face has two elements. Order complexes of posets
(in particular, barycentric subdivisions) are flag complexes due to the
transitivity relation.

For any subset $I\subset[m]$ denote by $K_I$ the subcomplex of $K$ consisting
of all simplices $J\in K$ such that $J\subset I$. The {\it link\/} of a
simplex $I$ of simplicial complex $K$ is the subcomplex $\link_K I\subset K$
consisting of all simplices $J\in K$ such that $I\cup J\in K$ and $I\cap
J=\emptyset$. For any vertex $\{i\}\in K$ the cone over $\link_K\{i\}$ (with
vertex~$\{i\}$) is naturally a subcomplex of~$K$, which is called the {\it
star\/} of $\{i\}$ and denoted $\star_K\{i\}$. The polyhedron
$|\star_K\{i\}|$ consists of all (polytopal) simplices of $|K|$ that
contain~$\{i\}$. When it is clear from the context which $K$ is under
consideration we write $\link I$ and $\star\{i\}$ instead of $\link_K I$ and
$\star_K\{i\}$ respectively.  Define $\core[m]=\{i\in[m]:\star\{i\}\ne K\}$.
The {\it core\/} of $K$ is the subcomplex $\core K=K_{\core[m]}$. Thus, the
core is the maximal subcomplex containing all vertices whose stars do not
coincide with~$K$.

\begin{example}
1. $\link(\emptyset)=K$.

2. Let $K$ be the simplex on four vertices $1,2,3,4$, and $I=\{1,2\}$. Then
$\link(I)$ is a subcomplex consisting of two vertices 3 and~4.

3. Let $K$ be the cone over $K'$ with vertex~$p$. Then $\link\{p\}=K'$,
$\star\{p\}=K$, and $\core K\subset K'$.
\end{example}

A {\it simplicial $n$-sphere\/} is a simplicial complex homeomorphic to
$n$-sphere~$S^n$. (Here and below by saying ``simplicial complex $K$ is
homeomorphic to~$X$" we mean that the geometric realisation $|K|$ is
homeomorphic to~$X$.) A {\it $PL$ (piecewise-linear) sphere\/} is a simplicial
sphere which is piecewise-linear homeomorphic to the boundary of a simplex.
The boundary of a simplicial $n$-polytope is an $(n-1)$-dimensional $PL$
sphere. $PL$ spheres that can be obtained in such way are called {\it
polytopal spheres\/}. Hence, we have the following inclusions of classes of
combinatorial objects:
$$
  \text{polytopal spheres }\subset\text{\ $PL$ spheres }\subset \text{\
  simplicial spheres.}
$$
In dimension 2 any simplicial sphere is polytopal (see
e.g.~\cite[Theorem~5.8]{Zi}). However, in higher dimensions both above
inclusions are strict. First examples of {\it non-polytopal\/} $PL$ 3-spheres
were found by Gr\"unbaum, and the smallest such sphere has 8 vertices. Good
description of these examples can be found in~\cite{Ba1}. A {\it non-$PL$\/}
simplicial sphere is presented in Example~\ref{non-PL} below.

The $f$-vector and the $h$-vector of an $(n-1)$-dimensional
simplicial complex $K^{n-1}$ are defined in the same way as for simplicial
polytopes: $\mb f(K^{n-1})=(f_0,f_1,\ldots,f_{n-1})$, where $f_i$ is the
number of $i$-dimensional simplices of~$K^{n-1}$, and
$\mb h(K^{n-1})=(h_0,h_1,\ldots,h_n)$, where $h_i$ are determined
by~(\ref{hvector}). Here we also assume $f_{-1}=1$. If $K^{n-1}=\partial P^*$
is the boundary of a simplicial $n$-polytope~$P^*$, then one obviously has
$\mb f(K^{n-1})=\mb f(P^*)$.

Since the $f$-vector of a polytopal sphere coincides with the $f$-vector of
the corresponding (simplicial) polytope, the $g$-theorem (Theorem~\ref{gth})
holds for polytopal spheres. So, it is natural to ask whether the $g$-theorem
extends to simplicial spheres. This question was posed by
McMullen~\cite{McM1} as an extension of his conjecture for simplicial
polytopes. After 1980, when the proof of McMullen's conjecture for simplicial
polytopes was found by Billera, Lee, and Stanley, the following became
perhaps the main open combinatorial-geometrical problem concerning the
$f$-vectors of simplicial complexes.

\begin{problem}[$g$-conjecture for simplicial spheres]
\label{gconj}
Does the $g$-theorem (Theorem~{\rm\ref{gth}}) hold for simplicial spheres?
\end{problem}

The $g$-conjecture is open even for $PL$ spheres. We note that only the
necessity of $g$-theorem (i.e. that the $g$-vector is an $M$-vector) should
be verified for simplicial spheres. If correct, the $g$-conjecture would
imply a characterisation of $f$-vectors of simplicial spheres.

The first part of Theorem~\ref{gth} (the Dehn--Sommerville equations) is
known to be true for simplicial spheres (see Corollary~\ref{dsgor} below).
The first inequality $h_0\le h_1$ from the second part of $g$-theorem is
equivalent to $1\le f_0-n$, which is obvious. The inequality $h_1\le h_2$
($n\ge4$) is equivalent to the lower bound $f_1\ge nf_0-\bin{n+1}2$ for the
number of edges, which is also known for simplicial spheres (see~\cite{Ba2},
in fact, the proof of the lower bound for the number of edges of
simplicial polytopes can be adapted for simplicial spheres). All these facts
together imply that the $g$-conjecture is true for simplicial spheres of
dimension~$\le4$. The inequality $h_2\le h_3$ ($n\ge6$) from the GLBT (the
second part of Theorem~\ref{gth}) is open. A lot of attempts to prove the
$g$-conjecture were made during the last two decades. Though unsuccessful,
these attempts resulted in some very interesting reformulations of the
$g$-conjecture. We just mention the results of Pachner~\cite{Pac1},
\cite{Pac2} reducing the $g$-conjecture (for $PL$-spheres) to some properties
of {\it bistellar moves\/}, and the results of~\cite{TWW} showing that the
conjecture follows from the {\it skeletal $r$-rigidity\/} of simplicial
$(n-1)$-sphere for $r\le\sbr n2$. It was shown independently by Kalai and
Stanley \cite[Corollary~2.4]{St3} that the GLBT holds for the boundary of an
$n$-dimensional ball that is a subcomplex of the boundary complex of a
simplicial $(n+1)$-polytope. However, it is not clear now which simplicial
complexes occur in this way. The lack of progress in proving the
$g$-conjecture motivated Bj\"orner and Lutz to launch the computer-aided seek
for counter examples~\cite{BjLu}. Though their computer program, BISTELLAR,
allowed to obtain many remarkable results on triangulations of manifolds, no
counter examples to $g$-conjecture were found. For more history of
$g$-theorem and related questions see~\cite{St4}, \cite{St5},
\cite[Lecture~8]{Zi}.

A simplicial complex $K$ is called a {\it simplicial manifold\/} (or {\it
triangulated manifold\/}) if the polyhedron $|K|$ is a topological manifold.
All manifolds considered in this review are compact, connected and closed
(unless otherwise stated). If $K^q$ is a simplicial manifold, then
$\link(I)$ has the homology of a $(q-\#I)$-sphere for each non-empty
simplex $I\in K^q$ (see Theorem~\ref{edwards} below). A $q$-dimensional {\it
$PL$ manifold\/} (or {\it combinatorial manifold\/}) is a simplicial complex
$K^q$ such that $\link(I)$ is a $PL$ sphere of dimension $(q-\#I)$ for each
non-empty simplex $I\in K^q$.

\begin{remark}
If $K^q$ is a $PL$ manifold, then for each vertex $\{i\}\in K^q$ the $PL$
$(q-1)$-sphere $\link\{i\}$ bounds the open neighbourhood $U_i$ which is
$PL$-homeomorphic to a $q$-ball. Since any point of $|K^q|$ is contained in
$U_i$ for some~$i$, this defines a $PL$-atlas on~$|K^q|$.
\end{remark}

\begin{example}[non-$PL$ simplicial 5-sphere]
\label{non-PL}
Let $S_H^3$ be any homology, but not topological, 3-sphere, i.e. a
non-simply-connected manifold with the same homology as the ordinary
3-sphere~$S^3$. The {\it Poincar\'e sphere\/} $\Sigma=SO(3)/A_5$ provides an
example of such a manifold. By Cannon's theorem~\cite{Ca}, the double
suspension $\Sigma^2S_H^3$ is homeomorphic to~$S^5$. However, $\Sigma^2S_H^3$
can not be $PL$, since $S_H^3$ appears as the link of some 1-simplex in
$\Sigma^2S_H^3$.
\end{example}

The following theorem gives a combinatorial characterisation of simplicial
complexes which are simplicial manifolds of dimension $\ge 5$ and generalises
the mentioned above result by Cannon.

\begin{theorem}[Edwards~\cite{Ed}]
\label{edwards}
  For $q\ge5$ the polyhedron of a simplicial complex $K^q$ is a topological
  $q$-manifold if and only if $\link I$ has the homology of a
  $(q-\#I)$-sphere for each non-empty simplex $I\in K^q$ and
  $\link\{i\}$ is simply connected for each vertex $\{i\}\in K$.
\end{theorem}

Which topological manifolds can be triangulated is the question of great
importance for combinatorial topologists. Any smooth manifold can be
triangulated by Whitney's theorem. All topological 2- and 3-dimensional
manifolds can be triangulated as well (for 2-manifolds this is almost
obvious, for 3-manifolds see~\cite{Mo}). Moreover, since the link of a vertex
in a simplicial 3-sphere is a 2-sphere (and 2-sphere is always $PL$), all 2-
and 3-manifolds are~$PL$. However, in dimension~4 there exist topological
manifolds that do not admit a $PL$-triangulation (e.g. Freedman's fake $\C
P^2$). Moreover, there exist topological 4-manifolds that {\it do not
admit any triangulation\/} (e.g. Freedman's topological 4-manifold with the
intersection form~$E_8$). Both examples can be found in~\cite{Ak}. In
dimensions $\ge5$ we have the famous combinatorial-topological problem:

\begin{problem}[Triangulation Conjecture]
  Is it true that every topological manifold of dimension $\ge5$ can be
  triangulated?
\end{problem}

Another well-known problem of $PL$-topology concerns the uniqueness of
a $PL$ structure on the topological sphere.

\begin{problem}
Is a $PL$ manifold homeomorphic to the topological {\rm 4}-sphere necessarily
a $PL$ sphere?
\end{problem}

Four is only dimension where the uniqueness of a $PL$ structure for the
topological sphere is open. For dimensions $\le 3$ the uniqueness was proved
by Moise~\cite{Mo}, and for dimensions $\ge 5$ it follows from the work of
Kirby and Siebenmann~\cite{KS}. In dimension~4 the category of $PL$ manifolds
is equivalent to the smooth category, hence, the above problem is equivalent
to if there exists an exotic 4-sphere.

More information about recent developments and open problems in combinatorial
and $PL$ topology can be found in \cite{No1},~\cite{Ra}.

\subsection{Simplicial complexes: commutative algebra}
\label{sim2}
The commutative algebra can be applied to combinatorics of
simplicial complexes and related objects. The main tool for translating
combinatorial results and problems to the algebraic language is the
Stanley--Reisner ring of simplicial complex. This approach was outlined by
R.~Stanley at the beginning of 1970's.

Remember that $\k[v_1,\ldots,v_m]$ denotes the graded polynomial algebra over
a commutative ring $\k$ with unit, $\deg(v_i)=2$.

\begin{definition}
\label{frsim}
The {\it face ring\/} (or the {\it Stanley--Reisner ring\/}) of a
simplicial complex $K$ with $m$ vertices is the quotient ring
$$
  \k(K)=\k[v_1,\ldots,v_m]/\mathcal I_K,
$$
where $\mathcal I_K$ is the homogeneous ideal generated by all square-free
monomials $v_{i_1}v_{i_2}\cdots v_{i_s}$, $i_1<\dots<i_s$, such that
$I=\{i_1,\ldots,i_s\}$ is not a simplex of~$K$.
\end{definition}

Note that the ideal $\mathcal I_K$ has basis consisting of square-free
monomials $v_{i_1}v_{i_2}\cdots v_{i_s}$ such that
$I=\{i_1,\ldots,i_s\}$ is a missing face of~$K$. Ideals in the polynomial
ring that admit a basis of momomials are called {\it monomial\/}.

\begin{proposition}\label{sfmi}
Any square-free monomial ideal of the polynomial ring has the form
$\mathcal I_K$ for some simplicial complex~$K$.
\end{proposition}
\begin{proof}
For any subset $I=\{i_1,\ldots,i_k\}\subset[m]$ denote by $v_I$ the
square-free monomial $v_{i_1}\cdots v_{i_k}$. Let $\mathcal I$ be a
square-free monomial ideal. Set
$$
  K=\{I\subset[m]\::\:v_{[m]\setminus I}\in\mathcal I\}.
$$
Then one easily checks that $K$ is a simplicial complex and $\mathcal
I=\mathcal I_K$.
\end{proof}

Let $P$ be a simple $n$-polytope, and $P^*$ its polar simplicial polytope.
Denote by $K_P$ the boundary of~$P^*$. Then $K_P$ is a polytopal simplicial
$(n-1)$-sphere. Obviously, the face ring of $P$ from Definition~\ref{frpol}
coincides with that of $K_P$ from Definition~\ref{frsim}: $\k(P)=\k(K_P)$.

Let $M=M^0\oplus M^1\oplus\dots$ be a graded $\k$-module. The series
$$
  F(M;t)=\sum_{i=0}^\infty(\dim_\k M^i)t^i
$$
is called the {\it Poincar\'e series\/} of~$M$.

\begin{remark}
In the algebraic literature the series $F(M;t)$ is known as the {\it Hilbert
series\/} or {\it Hilbert--Poincar\'e series\/}.
\end{remark}

The following lemma establishes the connection between two
combinatorial invariants of a simplicial complex: the face ring and the
$f$-vector (or the $h$-vector).

\begin{lemma}[Stanley~{\cite[Theorem~II.1.4]{St4}}]
\label{psfr}
  The Poincar\'e series of $\k(K^{n-1})$ can be calculated as
  $$
    F\bigl(\k(K^{n-1});t\bigr)=
    \sum_{i=-1}^{n-1}\frac{f_it^{2(i+1)}}{(1-t^2)^{i+1}}=
    \frac{h_0+h_1t^2+\dots+h_nt^{2n}}{(1-t^2)^n},
  $$
  where $(f_0,\ldots,f_{n-1})$ is the $f$-vector and $(h_0,\ldots,h_n)$ is
  the $h$-vector of~$K^{n-1}$.
\end{lemma}

Note that the second identity from Lemma~\ref{psfr} is an obvious corollary
of~(\ref{hf}).

\begin{example}\label{simbsim}
1. Let $K=\D^n$ (the $n$-simplex). Then $f_i=\bin{n+1}{i+1}$ for $-1\le i\le
n$, $h_0=1$, and $h_i=0$ for $i>0$. Since any subset of $[n+1]$ is a simplex
of $\D^n$, one has $\k(\D^n)=\k[v_1,\ldots,v_{n+1}]$. Then
$F(\k(\D^n);t)=\frac1{(1-t^2)^{n+1}}$, which agrees with Lemma~\ref{psfr}.

2. Let $K$ be the boundary of an $n$-simplex. Then $h_i=1$, $i=0,1,\ldots,n$,
and $\k(K)=\k[v_1,\ldots,v_{n+1}]/(v_1v_2\cdots v_{n+1})$. Hence,
$$
  F\bigl(\k(K);t\bigr)=\frac{1+t^2+\dots+t^{2n}}{(1-t^2)^n}.
$$
\end{example}

Now suppose $\k$ is a field. Let $A$ be a graded algebra over $\k$. The {\it
Krull dimension\/} of $A$ (denoted $\Kd A$) is the maximal number of
algebraically independent elements of~$A$. A sequence $\t_1,\ldots,\t_n$ of
$n=\Kd A$ homogeneous elements of $A$ is called a {\it hsop\/} (homogeneous
system of parameters) if the Krull dimension of the quotient
$A/(\t_1,\ldots,\t_n)$ is zero. Equivalently, $\t_1,\ldots,\t_n$ is a hsop if
$n=\Kd A$ and $A$ is a finitely-generated $\k[\t_1,\ldots,\t_n]$-module. The
elements of a hsop are algebraically independent.

\begin{lemma}[Noether normalisation lemma]
  For any finitely-generated graded algebra $A$ there exists a hsop. If $\k$
  is of zero characteristic and $A$ is generated by degree-two elements, then
  one can choose a degree-two hsop.
\end{lemma}

Below in this section we assume that $\k$ is of zero characteristic.
A sequence $\t_1,\ldots,\t_k$ of homogeneous elements of $A$ is called a
{\it regular sequence\/} if $\t_{i+1}$ is not a zero divisor in
$A/(\t_1,\ldots,\t_i)$ for $0\le i<k$ (i.e. the multiplication
by $\t_{i+1}$ is a monomorphism of $A/(\t_1,\ldots,\t_i)$ into itself).
Equivalently, $\t_1,\ldots,\t_k$ is a regular sequence if
$\t_1,\ldots,\t_k$ are algebraically independent and $A$ is a {\it free\/}
$\k[\t_1,\ldots,\t_k]$-module.

\begin{remark}
Regular sequences can be also defined in non-finitely-generated graded
algebras and in algebras over any integral domain.  Regular sequences in the
graded polynomial ring $R[a_1,a_2,\ldots,]$, $\deg a_i=-2i$, on infinitely
many generators, where $R$ is a subring of the field $\Q$ of rationals, are
used in the algebraic topology for constructing complex cobordism theories
with coefficients, see~\cite{La}.
\end{remark}

Any two maximal regular sequences have the same length, which is called the
{\it depth\/} of $A$ and denoted $\mathop{\rm depth}A$.
Obviously, $\mathop{\rm depth}A\le\Kd A$.

\begin{definition}
\label{CM}
  Algebra $A$ is called {\it Cohen--Macaulay\/} if it admits a
  regular sequence $\t_1,\ldots,\t_n$ of length $n=\Kd A$.
\end{definition}

Any regular sequence $\t_1,\ldots,\t_n$ of length $n=\Kd A$ is a hsop. It
follows that $A$ is Cohen--Macaulay if and only if there exists a hsop
$\t_1,\ldots,\t_n$ such that $A$ is a finitely-generated free
$\k[\t_1,\ldots,\t_n]$-module. If in addition $A$ is generated by degree-two
elements, then one can choose $\t_1,\ldots,\t_n$ to be of degree two. In this
case for the Poincar\'e series of $A$ holds
$$
  F(A;t)=\frac{F\bigl( A/(\t_1,\ldots,\t_n);t \bigr)}{(1-t^2)^n},
$$
where $F(A/(\t_1,\ldots,\t_n);t)=h_0+h_1t^2+\cdots$ is a polynomial. The
finite vector $(h_0,h_1,\ldots)$ is called the {\it $h$-vector\/} of~$A$.

A simplicial complex $K^{n-1}$ is called {\it Cohen--Macaulay\/} (over $\k$)
if its face ring $\k(K)$ is Cohen--Macaulay. Obviously, $\Kd\k(K)=n$.
Lemma~\ref{psfr} shows that the $h$-vector of $\k(K)$ coincides with the
$h$-vector of~$K$.

\begin{theorem}[Stanley]\label{stanmv}
  If $K^{n-1}$ is a Cohen--Macaulay simplicial complex, then
  $\mb h(K^{n-1})=(h_0,\ldots,h_n)$ is an $M$-vector (see
  section~{\rm\ref{sim1}}).
\end{theorem}
\begin{proof}
Let $\t_1,\ldots,\t_n$ be a regular sequence of degree-two elements
of~$\k(K)$. Then $A=\k(K)/(\t_1,\ldots,\t_n)$ is a graded algebra generated by
degree-two elements, and $\dim_\k A^{2i}=h_i$. Now the result follows from
Theorem~\ref{mvect}.
\end{proof}

The following fundamental theorem characterises Cohen--Macaulay complexes.

\begin{theorem}[Reisner~\cite{Re}]
\label{reisner}
  A simplicial complex $K$ is Cohen--Macaulay over $\k$ if and only if
  for any simplex $I\in K$ (including $I=\emptyset$) and $i<\dim(\link I)$,
  $\widetilde{H}_i(\link I;\k)=0$. (Here $\widetilde{H}_i(\:\cdot\:;\k)$
  denotes the $i$-th reduced homology group with coefficients in~$\k$.)
\end{theorem}

In particular, a simplicial sphere is a Cohen--Macaulay complex. Then
Theorem~\ref{stanmv} shows that the $h$-vector of a simplicial sphere is an
$M$-vector. This fact allowed Stanley to extend the UBT (Theorem~\ref{ubt})
to simplicial spheres.

\begin{corollary}[Stanley]
\label{ubtss}
  The Upper Bound Theorem holds for simplicial spheres. That is, if
  $(h_0,h_1,\ldots,h_n)$ is the $h$-vector of a simplicial $(n-1)$-sphere
  $K^{n-1}$ with $m$ vertices, then
  $$
    h_i(K^{n-1})\le\bin{m-n+i-1}i,\qquad 0\le i<{\textstyle\sbr n2}.
  $$
\end{corollary}
\begin{proof}
Since $\mb h(K^{n-1})$ is an $M$-vector, there exists a graded algebra
$A=A^0\oplus A^2\oplus\dots\oplus A^{2n}$ generated by degree-two elements
such that $\dim_\k A^{2i}=h_i$ (Theorem~\ref{mvect}). In particular, $\dim_\k
A^2=h_1=m-n$. Since $A$ is generated by~$A^2$, the number $h_i$ can not
exceed the total number of monomials of degree $i$ in $(m-n)$ variables. The
latter is exactly~$\bin{m-n+i-1}i$.
\end{proof}

\subsection{Homological properties of face rings
(Stanley--Reisner rings)}
\label{homo}
We start with reviewing some homology algebra. All modules in this section
are assumed to be finitely-generated graded $\k[v_1,\ldots,v_m]$-modules,
$\deg v_i=2$, unless otherwise stated.

A {\it finite free resolution\/} of a module $M$ is an exact sequence
\begin{equation}
\label{resol}
  0\to R^{-h} \stackrel d{\longrightarrow} R^{-h+1}
  \stackrel d{\longrightarrow}\cdots\longrightarrow R^{-1}
  \stackrel d{\longrightarrow} R^0\stackrel d{\longrightarrow} M\to 0,
\end{equation}
where the $R^{-i}$ are finitely-generated free modules and the maps $d$ are
degree-preserving. The minimal number $h$ for which a free
resolution~(\ref{resol}) exists is called the {\it homological dimension\/}
of $M$ and denoted $\hd M$. By the Hilbert syzygy theorem, $\hd M\le m$. A
resolution~(\ref{resol}) determines the free {\it bigraded differential
module\/} $[R,d]$, where $R=\bigoplus R^{-i,j}$, $R^{-i,j}:=(R^{-i})^j$ (the
$j$-th graded component of the free module $R^{-i}$). The cohomology of
$[R,d]$ is zero in non-zero dimensions, while $H^0[R,d]=M$. Conversely, a
free bigraded differential module $[R=\bigoplus_{i,j\ge0}
R^{-i,j},d:R^{-i,j}\to R^{-i+1,j}]$ with $H^0[R,d]=M$ and $H^{-i}[R,d]=0$ for
$i>0$ defines a free resolution~(\ref{resol}) with
$R^{-i}:=R^{-i,*}=\bigoplus_j R^{-i,j}$.

\begin{remark}
  For reasons specified below we numerate the terms of a free
  resolution by {\it non-positive\/} numbers, thus making it a {\it
  cochain\/} complex.
\end{remark}

The Poincar\'e series of $M$ can be calculated from any free
resolution~(\ref{resol}).

\begin{theorem}
\label{psresol}
  Suppose that the free $\k[v_1,\ldots,v_m]$-module $R^{-i}$
  in~{\rm(\ref{resol})} is generated by elements of degrees
  $d_{1i},\ldots,d_{q_ii}$, where $q_i=\mathop{\rm rank}R^{-i}$,
  $i=1,\ldots,h$. Then
  \begin{equation}
  \label{psresol1}
    F(M;t)=(1-t^2)^{-m}
    \sum_{i=0}^h(-1)^i(t^{d_{1i}}+\dots+t^{d_{q_ii}}).
  \end{equation}
\end{theorem}

\begin{example}[minimal resolution]
\label{minimal}
For different reasons it is convenient to have a resolution~(\ref{resol}) for
which each term $R^{-i}$ has the smallest possible rank. The following
definition is taken from Adams' paper~\cite{Ad}. Let $M$,~$M'$ be two
modules. Set $\mathcal J(M)=v_1M+v_2M+\dots+v_mM\subset M$. A map $f:M\to M'$
is called {\it minimal\/} if $\Ker f\subset\mathcal J(M)$. A
resolution~(\ref{resol}) is called {\it minimal\/} if all maps $d$ are
minimal.  Then it is easy to see that each $R^{-i}$ has the smallest
possible rank.

A minimal resolution can be constructed as follows.  Take a {\it minimal set
of generators\/} $a_1,\ldots,a_{k_0}$ for $M$ and define $R^0$ to be the free
$\k[v_1,\ldots,v_m]$-module with $k_0$ generators. Then take a minimal set of
generators $a_1,\ldots,a_{k_1}$ in the kernel of natural epimorphism $R^0\to
M$ and define $R^{-1}$ to be the free $\k[v_1,\ldots,v_m]$-module with $k_1$
generators, and so on. On the $i$-th step we take a minimal set of
generators in the kernel of the previously constructed map $d:R^{-i+1}\to
R^{-i+2}$ and define $R^{-i}$ to be the free module with the corresponding
generators. Note that a minimal resolution is unique up to an isomorphism.
\end{example}

\begin{example}[Koszul resolution]
\label{koszul}
Let $M=\k$. The $\k[v_1,\ldots,v_m]$-module structure on $\k$ is defined by
the map $\k[v_1,\ldots,v_m]\to\k$ that sends each $v_i$ to~0. Let
$\L[u_1,\ldots,u_m]$ denote the exterior algebra on $m$ generators.
The tensor product $R=\L[u_1,\ldots,u_m]\otimes\k[v_1,\ldots,v_m]$
(here and below $\otimes$ denotes $\otimes_\k$) becomes a
{\it differential bigraded algebra\/} by setting
\begin{gather}
  \bideg u_i=(-1,2),\quad\bideg v_i=(0,2),\notag\\
  \label{diff}
  du_i=v_i,\quad dv_i=0,
\end{gather}
and requiring that $d$ be a derivation of algebras. It can be shown
(see~\cite[\S 7.2]{Mac}) that $H^{-i}[R,d]=0$ for $i>0$ and
$H^0[R,d]=\k$. Since $\L[u_1,\ldots,u_m]\otimes\k[v_1,\ldots,v_m]$ is a
free $\k[v_1,\ldots,v_m]$-module, it determines a free resolution of~$\k$.
This resolution is called the {\it Koszul resolution\/}. More precisely, it
has the following form
\begin{multline*}
0\to\L^m[u_1,\ldots,u_m]\otimes\k[v_1,\ldots,v_m]\longrightarrow\cdots\\
\longrightarrow\L^1[u_1,\ldots,u_m]\otimes\k[v_1,\ldots,v_m]\longrightarrow
\k[v_1,\ldots,v_m]\longrightarrow\k\to0,
\end{multline*}
where $\L^i[u_1,\ldots,u_m]$ is the submodule of $\L[u_1,\ldots,u_m]$ spanned
by the monomials of length~$i$. Thus, in notations of~(\ref{resol}) we have
$R^{-i}=\L^i[u_1,\ldots,u_m]\otimes\k[v_1,\ldots,v_m]$.
\end{example}

Let $N$ be another module; then applying the functor
$\otimes_{\k[v_1,\ldots,v_m]}N$ to~(\ref{resol}) we obtain the following
cochain complex of graded modules:
$$
  0 \longrightarrow R^{-h}\otimes_{\k[v_1,\ldots,v_m]}N \longrightarrow
  \cdots \longrightarrow R^0\otimes_{\k[v_1,\ldots,v_m]}N
  \longrightarrow 0
$$
and the corresponding bigraded differential module $[R\otimes N,d]$.
The $(-i)$-th cohomology module of the above cochain complex is denoted
$\Tor^{-i}_{\k[v_1,\ldots,v_m]}(M,N)$, i.e.
\begin{multline*}
  \Tor^{-i}_{\k[v_1,\ldots,v_m]}(M,N):=H^{-i}[R\otimes_{\k[v_1,\ldots,v_m]}
  N,d]\\=\frac{\Ker[d:R^{-i}\otimes_{\k[v_1,\ldots,v_m]}
  N\to R^{-i+1}\otimes_{\k[v_1,\ldots,v_m]} N]}
  {d(R^{-i-1}\otimes_{\k[v_1,\ldots,v_m]}N)}.
\end{multline*}
Since $R^{-i}$ and $N$ are graded modules, one has
$$
  \Tor^{-i}_{\k[v_1,\ldots,v_m]}(M,N)=
  \bigoplus_j\Tor^{-i,j}_{\k[v_1,\ldots,v_m]}(M,N),
$$
where
$$
  \Tor^{-i,j}_{\k[v_1,\ldots,v_m]}(M,N)=
  \frac{\Ker\bigl[d:(R^{-i}\otimes_{\k[v_1,\ldots,v_m]}
  N)^j\to(R^{-i+1}\otimes_{\k[v_1,\ldots,v_m]} N)^j\bigl]}
  {d(R^{-i-1}\otimes_{\k[v_1,\ldots,v_m]} N)^j}.
$$
Thus, we have a {\it bigraded\/} $\k[v_1,\ldots,v_m]$-module
$$
  \Tor_{\k[v_1,\ldots,v_m]}(M,N)=
  \bigoplus_{i,j}\Tor^{-i,j}_{\k[v_1,\ldots,v_m]}(M,N).
$$

The following properties of $\Tor^{-i}_{\k[v_1,\ldots,v_m]}(M,N)$ are
well known (see e.g.~\cite{Mac}).
\begin{proposition}
\label{torprop}
{\rm (a)} The module $\Tor^{-i}_{\k[v_1,\ldots,v_m]}(M,N)$ does not depend,
up to isomorphism, on the choice of resolution~{\rm(\ref{resol})}.

{\rm (b)} Both $\Tor^{-i}_{\k[v_1,\ldots,v_m]}(\:\cdot\:,N)$ and
$\Tor^{-i}_{\k[v_1,\ldots,v_m]}(M,\:\cdot\:)$ are covariant functors.

{\rm (c)} $\Tor^0_{\k[v_1,\ldots,v_m]}(M,N)\cong
M\otimes_{\k[v_1,\ldots,v_m]}N$.

{\rm (d)} $\Tor^{-i}_{\k[v_1,\ldots,v_m]}(M,N)\cong
\Tor^{-i}_{\k[v_1,\ldots,v_m]}(N,M)$.
\end{proposition}

One can also define the $A$-modules $\Tor_A(M,N)$ for any finitely-generated
graded commutative algebra $A$ and (finitely-generated graded) $A$-modules
$M$,~$N$. Though an $A$-free resolution~(\ref{resol}) of $M$ may fail to
exist, there always exists a {\it projective\/} resolution of~$M$, which
allows to define $\Tor_A(M,N)$ in the same way as above. Note that projective
modules over polynomial algebra are free. This was known as the {\it Serre
problem\/}, now solved by Quillen and Suslin. However, for graded modules this
fact is not hard to prove. In our paper modules $\Tor_A(M,N)$ for algebras
$A$ different from the polynomial ring appear only in sections~\ref{eile}
and~\ref{diag}.

Now let $K^{n-1}$ be a simplicial complex on $m$ vertices, $M=\k(K)$ and
$N=\k$. Since $\deg v_i=2$, we have
$$
  \Tor_{\k[v_1,\ldots,v_m]}\bigl(\k(K),\k\bigr)=
  \bigoplus_{i,j=0}^m\Tor^{-i,2j}_{\k[v_1,\ldots,v_m]}\bigl(\k(K),\k\bigr)
$$
(i.e. non-zero elements of $\Tor_{\k[v_1,\ldots,v_m]}(\k(K),\k)$ have even
second degree). Define the {\it bigraded Betti numbers\/} of $\k(K)$ as
\begin{equation}
\label{bbnfr}
  \b^{-i,2j}\bigl(\k(K)\bigr):=
  \dim_\k\Tor^{-i,2j}_{\k[v_1,\ldots,v_m]}\bigl(\k(K),\k\bigr),\qquad
  0\le i,j\le m.
\end{equation}
Suppose that (\ref{resol}) is a {\it minimal\/} free resolution of $M=\k(K)$.
Then $R^0\cong\k[v_1,\ldots,v_m]$ is a free $\k[v_1,\ldots,v_m]$-module with
one generator of degree~0. The basis of $R^{-1}$ (the minimal generator set
for $\Ker[\k[v_1,\ldots,v_m]\to\k(K)]$) consists of elements
$v_{i_1,\ldots,i_k}$, $\deg v_{i_1,\ldots,i_k}=2k$, such that
$\{i_1,\ldots,i_k\}$ is a missing face of~$K$. The map $d:R^{-1}\to R^0$
takes $v_{i_1,\ldots,i_k}$ to $v_{i_1}\cdots v_{i_k}$.
Since the maps $d$ in~(\ref{resol}) are minimal, all differentials in
the cochain complex
$$
  0 \longrightarrow R^{-h}\otimes_{\k[v_1,\ldots,v_m]}\k \longrightarrow
  \cdots \longrightarrow R^0\otimes_{\k[v_1,\ldots,v_m]}\k
  \longrightarrow 0
$$
are trivial. Hence, for the minimal resolution of $\k(K)$ holds
\begin{gather}
\label{tormin}
  \Tor^{-i}_{\k[v_1,\ldots,v_m]}\bigl(\k(K),\k\bigr)\cong
  R^{-i}\otimes_{\k[v_1,\ldots,v_m]}\k,\\
  \b^{-i,2j}\bigl(\k(K)\bigr)=\mathop{\rm rank}R^{-i,2j}
  (=\dim_{\k[v_1,\ldots,v_m]}R^{-i,2j}).\notag
\end{gather}

The Betti numbers $\beta^{-i,2j}(\k(K))$ are important combinatorial
invariants of the simplicial complex~$K$. Some results describing these
numbers were obtained in~\cite{St4}. The following important theorem (which
was proved by combinatorial methods) reduces the calculation of numbers
$\beta^{-i,2j}(\k(K))$ to calculating the homology groups of
subcomplexes of~$K$.

\begin{theorem}[Hochster \cite{Ho}]\label{hoch}
The following formula holds for the Poincar\'e series of the module
$\Tor^{-i}_{\k[v_1,\ldots,v_m]}(\k(K),\k)$:
$$
  \sum_j\beta^{-i,2j}\bigl(\k(K)\bigr)t^{2j}=\sum_{I\subset[m]}
  \bigl( \dim_\k\widetilde{H}_{\#I-i-1}(K_I)\bigr)t^{2(\#I)},
$$
where $K_I$ is the subcomplex of $K$ consisting of all simplices with
vertices in~$I$.
\end{theorem}
\noindent We mention that calculations using this theorem become very
involved even for small~$K$. In chapter~4 we show that the numbers
$\beta^{-i,2j}(\k(K))$ equal the bigraded Betti numbers of the moment-angle
complex $\zk$ associated with the simplicial complex~$K$. This provides the
alternative (topological) way for calculating the numbers
$\beta^{-i,2j}(\k(K))$.

Now we look at the Koszul resolution (Example~\ref{koszul}).

\begin{lemma}
\label{koscom}
For any module $M$ holds
$$
  \Tor_{\k[v_1,\ldots,v_m]}(M,\k)\cong
  H\bigl[\L[u_1,\ldots,u_m]\otimes M,d\bigr],
$$
where $H[\L[u_1,\ldots,u_m]\otimes M,d]$ is the cohomology of the
bigraded differential module $\L[u_1,\ldots,u_m]\otimes M$ and $d$ is defined
as in~{\rm(\ref{diff})}.
\end{lemma}
\begin{proof}
  Using the Koszul resolution
  $[\L[u_1,\ldots,u_m]\otimes\k[v_1,\ldots,v_m],d]$ for the definition of
  $\Tor_{\k[v_1,\ldots,v_m]}(\k,M)$ we calculate
  \begin{multline*}
    \Tor_{\k[v_1,\ldots,v_m]}(M,\k)\cong
    \Tor_{\k[v_1,\ldots,v_m]}(\k,M)\\
    =H\bigl[ \L[u_1,\ldots,u_m]\otimes\k[v_1,\ldots,v_m]
    \otimes_{\k[v_1,\ldots,v_m]}M \bigr]\cong
    H\bigl[ \L[u_1,\ldots,u_m]\otimes M\bigr].
  \end{multline*}
\end{proof}

Suppose now that the $\k[v_1,\ldots,v_m]$-module $M$ is a $\k$-algebra. Then
the cohomology of $[\L[u_1,\ldots,u_m]\otimes M,d]$ is an algebra as well.
Lemma~\ref{koscom} allows to invest $\Tor_{\k[v_1,\ldots,v_m]}(M,\k)$ with
the canonical structure of a finite-dimensional bigraded {\it
$\k$-algebra\/}. We refer to this bigraded algebra as the {\it
$\Tor$-algebra\/} of~$M$. The Tor-algebra of a simplicial complex~$K$ is the
Tor-algebra of the face ring $\k(K)$.

\begin{remark}
For general $N\ne\k$ the module $\Tor_{\k[v_1,\ldots,v_m]}(M,N)$ has no
canonical structure of an algebra even if both $M$ and $N$ are algebras.
\end{remark}

\begin{construction}[multigraded structure in the Tor-algebra]
\label{mgrad}
Invest the polynomial ring $\k[v_1,\ldots,v_m]$ with the multigraded (more
precisely, $\N^m$-graded) structure by setting $\mathop{\rm mdeg}
v_i=(0,\ldots,0,2,0,\ldots,0)$, where 2 stands at the $i$-th place. Then the
multidegree of monomial $v_1^{i_1}\cdots v_m^{i_m}$ is $(2i_1,\ldots,2i_m)$.
Suppose that the algebra $M$ is the quotient of the polynomial ring by a
monomial ideal. Then the multigraded structure descends to~$M$ and to the
terms of resolution~(\ref{resol}). We may assume that the differentials in
the resolution preserve the multigraded structures. Then the module
$\Tor_{\k[v_1,\ldots,v_m]}(M,N)$ acquires the canonical
$\N\oplus\N^m$-grading, i.e.
$$
  \Tor_{\k[v_1,\ldots,v_m]}(M,\k)=\bigoplus_{i\ge0,\mb j\in\N^m}
  \Tor^{-i,\mb j}_{\k[v_1,\ldots,v_m]}(M,\k).
$$
In particular, the Tor-algebra of $K$ can be canonically made
an $\N\oplus\N^m$-graded algebra.
\end{construction}

\begin{remark}
According to our agreement the first degree in the Tor-algebra is
{\it non-positive\/}. (Recall that we numerate the terms of Koszul
$\k[v_1,\ldots,v_m]$-free resolution of $\k$ by non-positive integers.) In
such notations the Koszul complex $[M\otimes\Lambda[u_1,\ldots,u_m],d]$
becomes a {\it cochain\/} complex, and $\Tor_{\k[v_1,\ldots,v_m]}(M,\k)$
is its {\it cohomology\/}, not the homology as usually regarded. This is the
standard trick used for applying the Eilenberg--Moore spectral sequence, see
section~\ref{eile}. It also explains why we write
$\Tor^{*,*}_{\k[v_1,\ldots,v_m]}(M,\k)$ instead of usual
$\Tor_{*,*}^{\k[v_1,\ldots,v_m]}(M,\k)$.
\end{remark}

The upper bound $\hd M\le m$ from the Hilbert syzygy theorem
can be replaced by the following sharper result.

\begin{theorem}[Auslander and Buchsbaum]
  $\hd M=m-\mathop{\rm depth}M$.
\end{theorem}

From now on we assume that $M$ is generated by degree-two elements and the
$\k[v_1,\ldots,v_m]$-module structure on $M$ is defined by an epimorphism
$p:\k[v_1,\ldots,v_m]\to M$ (both assumptions are satisfied by definition
when $M=\k(K)$ for some~$K$). Suppose that $\t_1,\ldots,\t_k$ is a regular
sequence of degree-two elements of~$M$. Let $\mathcal
J:=(\t_1,\ldots,\t_k)\subset M$ be the ideal generated by $\t_1,\ldots,\t_k$.
Choose degree-two elements $t_i\subset\k[v_1,\ldots,v_m]$ such that
$p(t_i)=\t_i$, $i=1,\ldots,k$. The ideal in $\k[v_1,\ldots,v_m]$ generated by
$t_1,\ldots,t_k$ will be also denoted~$\mathcal J$. Then
$\k[v_1,\ldots,v_m]/\mathcal J\cong\k[w_1,\ldots,w_{m-k}]$. Under these
assumptions we have the following reduction lemma.

\begin{lemma}
\label{tortor}
  The following isomorphism of algebras holds:
  $$
    \Tor_{\k[v_1,\ldots,v_m]}(M,\k)=
    \Tor_{\k[v_1,\ldots,v_m]/\mathcal J}(M/\mathcal J,\k).
  $$
\end{lemma}

In order to prove the lemma we need the following fact from the homology
algebra.

\begin{theorem}[{\cite[p.349]{CE}}]
\label{change}
  Let $\Lambda$ be an algebra, $\Gamma$ its subalgebra, and
  $\Omega=\Lambda{/}\Gamma$ the quotient algebra. Suppose that $\Lambda$ is a
  free $\Gamma$-module and we are given an $\Omega$-module $A$ and a
  $\Lambda$-module~$C$. Then there exists a spectral sequence $\{E_r,d_r\}$
  such that
  $$
    E_r\Rightarrow \Tor_{\Lambda}(A,C),\quad
    E_2=\Tor_{\Omega}\bigl(A,\Tor_{\Gamma}(C,\k)\bigr).
  $$
\end{theorem}

\begin{proof}[Proof of Lemma~{\rm\ref{tortor}}]
Set $\L=\k[v_1,\ldots,v_m]$, $\Gamma=\k[t_1,\ldots,t_k]$, $A=\k$, $C=M$.
Then $\L$ is a free $\Gamma$-module and
$\O=\L{/}\Gamma=\k[v_1,\ldots,v_m]/\mathcal J$.  Therefore,
Theorem~\ref{change}, gives a spectral sequence
$$
  E_r\Rightarrow \Tor_{\k[v_1,\ldots,v_m]}(M,\k),\quad
  E_2=\Tor_{\O}\bigl(\Tor_{\Gamma}(M,\k),\k\bigr).
$$
Since $\t_1,\ldots,\t_k$ is a regular sequence, $M$ is a free $\Gamma$-module.
Therefore,
$$
  \Tor_{\Gamma}(M,\k)=M\otimes_{\Gamma}\k=M/\mathcal J\text{\quad and \quad}
  \Tor_{\Gamma}^q(M,\k)=0\text{\quad for }\;q\ne 0.
$$
It follows that $E_2^{p,q}=0$ for $q\ne 0$. Thus, the spectral sequence
collapses at the $E_2$ term, and we have
$$
  \Tor_{\k[v_1,\ldots,v_m]}(M,\k)=
  \Tor_{\O}\bigl(\Tor_{\Gamma}(M,\k),\k\bigr)=
  \Tor_{\k[v_1,\ldots,v_m]/\mathcal J}(M/\mathcal J,\k),
$$
which concludes the proof.
\end{proof}

It follows from Lemma~\ref{tortor} that if $M$ is Cohen--Macaulay of Krull
dimension~$n$, then $\mathop{\rm depth}M=n$, $\hd M=m-n$, and
$\Tor^{-i}_{\k[v_1,\ldots,v_m]}(M,\k)=0$ for $i>m-n$.

\begin{definition}
\label{goren}
  Suppose $M$ is a Cohen--Macaulay algebra of Krull dimension~$n$. Then $M$
  is called a {\it Gorenstein algebra\/} if
  $\Tor^{-(m-n)}_{\k[v_1,\ldots,v_m]}(M,\k)\cong\k$.
\end{definition}
Following Stanley~\cite{St4}, we call a simplicial complex $K$
{\it Gorenstein\/} if $\k(K)$ is a Gorenstein algebra. Further, $K$ is called
{\it Gorenstein*\/} if $\k(K)$ is Gorenstein and $K=\core K$ (see
section~\ref{sim1}). The following theorem characterises Gorenstein*
simplicial complexes.

\begin{theorem}[{\cite[\S II.5]{St4}}]
\label{gorencom}
  A simplicial complex $K$ is Gorenstein* over $\k$ if and only if
  for any simplex $I\in K$ (including $I=\emptyset$) the subcomplex
  $\link I$ has the homology of a sphere of dimension $\dim\,(\link I)$.
\end{theorem}

In particular, simplicial spheres and simplicial homology spheres (simplicial
manifolds with the homology of a sphere) are Gorenstein* complexes. Note,
however, that a Gorenstein* complex is not necessarily a simplicial
manifold (links of vertices are not necessarily simply connected, compare
with Theorem~\ref{edwards}).

\begin{theorem}[{\cite[\S II.5]{St4}}]
\label{tordual}
  Suppose $K^{n-1}$ is a Gorenstein* complex. Then the following identities
  hold for the Poincar\'e series
  of $\Tor^{-i}_{\k[v_1,\ldots,v_m]}(\k(K),\k)$, $0\le i\le m-n$:
  $$
    F\left(\Tor^{-i}_{\k[v_1,\ldots,v_m]}\bigl(\k(K),\k\bigr);\;t\right)=
    t^{2m}F\left(\Tor^{-(m-n)+i}_{\k[v_1,\ldots,v_m]}
    \bigl(\k(K),\k\bigr);\;{\textstyle\frac 1t}\right).
  $$
\end{theorem}

\begin{corollary}
\label{frdual}
  If $K^{n-1}$ is Gorenstein* then
  $$
    F\bigl(\k(K),t\bigr)=(-1)^nF\bigl(\k(K),{\textstyle\frac1t}\bigr).
  $$
\end{corollary}
\begin{proof}
We apply Theorem~\ref{psresol} to a minimal resolution of $\k(K)$. It
follows from~(\ref{tormin}) that the numerators of the summands in the right
hand side of~(\ref{psresol1}) are exactly
$F\bigl(\Tor^{-i}_{\k[v_1,\ldots,v_m]}(\k(K),\k);t\bigr)$, $i=1,\ldots,m-n$.
Hence,
$$
  F\bigl(\k(K);t\bigr)=(1-t^2)^{-m}\sum_{i=0}^{m-n}(-1)^i
  F\Bigl(\Tor^{-i}_{\k[v_1,\ldots,v_m]}\bigl(\k(K),\k\bigr);t\Bigr).
$$
Using Theorem~\ref{tordual} we get
\begin{multline*}
  F\bigl(\k(K);t\bigr)=(1-t^2)^{-m}\sum_{i=0}^{m-n}(-1)^i t^{2m}
  F\Bigl(\Tor^{-(m-n)+i}_{\k[v_1,\ldots,v_m]}\bigl(\k(K),\k\bigr);
  {\textstyle\frac1t}\Bigr)\\
  =\bigl(1-({\textstyle\frac1t})^2\bigr)^{-m}
  (-1)^m\sum_{j=0}^{m-n}(-1)^{m-n-j}
  F\Bigl(\Tor^{-j}_{\k[v_1,\ldots,v_m]}\bigl(\k(K),\k\bigr);
  {\textstyle\frac1t}\Bigr)\\
  =(-1)^nF\bigl(\k(K);{\textstyle\frac1t}\bigr).
\end{multline*}
\end{proof}

\begin{corollary}
\label{dsgor}
  The Dehn--Sommerville relations $h_i=h_{n-i}$, $0\le i\le n$, hold for
  any Gorenstein* complex $K^{n-1}$ (in particular, for any simplicial
  sphere).
\end{corollary}
\begin{proof}
This follows from Lemma~\ref{psfr} and Corollary~\ref{frdual}.
\end{proof}

As it was pointed out by Stanley in~\cite{St3}, Gorenstein* complexes are the
most general objects appropriate for the generalisation of the $g$-theorem.
(As we have seen, polytopal spheres, $PL$ spheres, simplicial spheres
and simplicial homology spheres are particular cases of Gorenstein*
complexes).

The Dehn--Sommerville equations can be generalized even beyond Gorenstein*
complexes. In~\cite{Kl} Klee reproved the Dehn--Sommerville equations in
the form~(\ref{DSf}) in the more general context of {\it Eulerian
manifolds\/}. In particular, this implies that equations~(\ref{DSf}) (with
the exception of $k=0$) hold for any simplicial manifold $K$ of
dimension~$n-1$. (In the case $k=0$ equation~(\ref{DSf}) expresses that
$\chi(K^{n-1})=\chi(S^{n-1})$.) Analogues of equations~(\ref{DSf}) were
obtained by Bayer and Billera~\cite{BB} (for {\it Eulerian posets\/}) and
Chen and Yan~\cite{CY} (for arbitrary polyhedra).

In section~\ref{coh3} we obtain (by topological methods) the following
form of the Dehn--Sommerville equations for simplicial manifolds:
$$
  h_{n-i}-h_i=(-1)^i\bigl(\chi(K^{n-1})-\chi(S^{n-1})\bigr)\bin ni,
  \quad i=0,1,\ldots,n.
$$
where $\chi(K^{n-1})=f_0-f_1+\ldots+(-1)^{n-1}f_{n-1}=1+(-1)^{n-1}h_n$ is the
Euler characteristic of $K^{n-1}$ and $\chi(S^{n-1})=1+(-1)^{n-1}$ is
that of sphere. Note that if $K$ is a simplicial sphere or has odd dimension,
then the above equations reduce to the classical $h_{n-i}=h_i$.

\subsection{Cubical complexes and cubical maps}
\label{cubi}
Define a $q$-dimensional {\it com\-bi\-na\-to\-rial-\-geo\-met\-ri\-cal
cube\/} as a polytope combinatorially equivalent to the standard
$q$-cube~(\ref{cube}). In this section cubes are
com\-bi\-na\-to\-rial-\-geo\-met\-ri\-cal cubes.

\begin{definition}
\label{cubcom}
A {\it cubical complex\/} is a subset $\mathcal C\subset\R^n$ represented as
a union of cubes of any dimensions in such a way that the intersection of any
two cubes is a face of each.
\end{definition}

\begin{remark}
The above definition of cubical complex is similar to
Definition~\ref{polyhed} of geometrical simplicial complex. One could also
define an {\it abstract cubical complex\/}, however this definition is more
subtle (and we do not need it).
\end{remark}

A {\it face\/} of a cubical complex $\mathcal C$ is a face of a cube
from~$\mathcal C$. The {\it dimension\/} of $\mathcal C$ is the maximal
dimension of its faces. The {\it $f$-vector\/} of a cubical complex $\mathcal
C$ is defined in the standard way ($f_i$ is the number of $i$-faces). Some
problems concerning the $f$-vectors of cubical complexes are discussed
in~\cite{St5}.

Obviously, the standard cube $I^q$ (together with all its faces) is a
$q$-dimensional cubical complex, which we also denote~$I^q$.
Any face of $I^q$ has the form
\begin{equation}
\label{ijface}
  C_{I\subset J}=\{(y_1,\ldots,y_q)\in I^q\: : \: y_i=0\text{ for }i\in
  I,\; y_i=1\text{ for }i\notin J\},
\end{equation}
were $I\subset J$ are two (possibly empty) subsets of~$[q]$. We denote
$C_J:=C_{\emptyset\subset J}$.

Unlike simplicial complexes (which are always subcomplexes of a
simplex), not any cubical complex can be realised as a subcomplex of
some~$I^q$. One example of a cubical complex not embedable as a
subcomplex in any $I^q$ is shown on Figure~1. Moreover, this complex is not
embedable into the standard cubical lattice in $\R^q$ (for any~$q$). The
authors are thankful to M.\,I.~Shtogrin for giving this example.
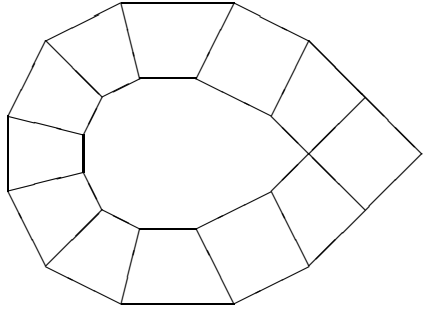
\begin{figure}
\begin{center}
\begin{picture}(55,40)
\put(0,15){\line(0,1){10}}
\put(0,25){\line(1,2){5}}
\put(5,35){\line(2,1){10}}
\put(15,40){\line(1,0){15}}
\put(30,40){\line(2,-1){10}}
\put(40,35){\line(1,-1){10}}
\put(40,35){\line(1,-1){15}}
\put(55,20){\line(-1,-1){15}}
\put(40,5){\line(-2,-1){10}}
\put(30,0){\line(-1,0){15}}
\put(15,0){\line(-2,1){10}}
\put(5,5){\line(-1,2){5}}
\put(10,17.5){\line(0,1){5}}
\put(10,22.5){\line(1,2){2.5}}
\put(12.5,27.5){\line(2,1){5}}
\put(17.5,30){\line(1,0){7.5}}
\put(25,30){\line(2,-1){10}}
\put(35,25){\line(1,-1){5}}
\put(40,20){\line(-1,-1){5}}
\put(35,15){\line(-2,-1){10}}
\put(25,10){\line(-1,0){7.5}}
\put(17.5,10){\line(-2,1){5}}
\put(12.5,12.5){\line(-1,2){2.5}}
\put(0,25){\line(4,-1){10}}
\put(5,35){\line(1,-1){7.5}}
\put(15,40){\line(1,-4){2.5}}
\put(30,40){\line(-1,-2){5}}
\put(40,35){\line(-1,-2){5}}
\put(47.5,27.5){\line(-1,-1){7.5}}
\put(47.5,12.5){\line(-1,1){7.5}}
\put(40,5){\line(-1,2){5}}
\put(30,0){\line(-1,2){5}}
\put(15,0){\line(1,4){2.5}}
\put(5,5){\line(1,1){7.5}}
\put(0,15){\line(4,1){10}}
\end{picture}
\caption{Cubical complex not embedable into cubical lattice}
\end{center}
\end{figure}

There is the following problem.
\begin{problem}[S.\,P.~Novikov]
\label{novikov}
Characterise $k$-dimensional cubical complexes
$\mathcal C$ (in particular, cubical manifolds) which admit

{\rm(a)} a (cubical) embedding into the standard cubical lattice
in~$\R^q${\rm;}

{\rm(b)} a map to the standard cubical lattice in~$\R^q$ whose restriction
to every $k$-dimensional cube identifies it with a certain
$k$-face of the lattice.
\end{problem}

In the case when $\mathcal C$ is homeomorphic to~$S^2$ the above problem was
solved in~\cite{DSS}. Problem~\ref{novikov} is an extension of the following
problem, formulated in~\cite{DSS}.

\begin{problem}[S.\,P.~Novikov]
\label{novikov1}
Suppose we are given a {\rm2}-dimensional cubical $\mod 2$ cycle
$\alpha$ in the standard cubical lattice in~$\R^3$. Describe all maps of
cubical subdivisions of {\rm2}-dimensional surfaces onto~$\alpha$ such that no
two different squares are mapped to the same square of~$\alpha$.
\end{problem}

As it was told to the authors by S.\,P.~Novikov, problem~\ref{novikov1}
have arisen in connection with the 3-dimensional Ising model during the
discussions with the well-known physicist A.\,M.~Polyakov.

Below we introduce some special cubical complexes, which will play a pivotal
r\^ole in our theory of moment-angle complexes. Each of these cubical
complexes admits a canonical cubical embedding into the standard cube. We
note that the problem of embedability into the cubical lattice is
closely connected with that of embedability into the standard cube. For
instance, as it was shown in~\cite{DSS}, if a cubical subdivision of a
2-dimensional surface is embedable into the standard cubical lattice
in~$\R^q$, then it also admits a cubical embedding into~$I^q$.

\begin{construction}[canonical simplicial subdivision of $I^{m}$]
\label{conbar}
Let $\D^{m-1}$ be the simplex on the set~$[m]$, i.e. $\D^{m-1}$ is the
collection of all subsets of~$[m]$. Assign to each subset
$I=\{i_1,\ldots,i_k\}\subset[m]$ the vertex $v_I:=C_{I\subset I}$ of~$I^m$.
That is, $v_I=(\e_1,\ldots,\e_m)$, where $\e_i=0$ if $i\in I$ and $\e_i=1$
otherwise. Regarding $I$ as a vertex of the barycentric
subdivision of $\D^{m-1}$, we can extend the mapping $I\mapsto v_I$ to the
piecewise linear embedding $i_c$ of the polyhedron $|\bs(\D^{m-1})|$ into the
(boundary complex of) standard cube~$I^m$. Under this embedding, the vertices
of $|\D^{m-1}|$ are mapped to the vertices $(1,\ldots,1,0,1,\ldots,1)\in
I^m$, while the barycentre of $|\D^{m-1}|$ is mapped to the vertex
$(0,\ldots,0)\in I^m$. The image $i_c(|\bs(\D^{m-1})|)$ is the union of $m$
facets of $I^m$ meeting at the vertex $(0,\ldots,0)$. For any pair $I$, $J$
of non-empty subsets of $[m]$ such that $I\subset J$ all simplices of
$\bs(\D^{m-1})$ of the form $I=I_1\subset I_2\subset\dots\subset I_k=J$ are
mapped to the same face $C_{I\subset J}\subset I^m$ (see~(\ref{ijface})).
The map $i_c:|\bs(\D^{m-1})|\to I^m$ extends to $|\cone(\bs(\D^{m-1}))|$ by
taking the vertex of the cone to $(1,\ldots,1)\in I^m$. The resulting map is
denoted by $\cone(i_c)$. Its image is the whole~$I^m$. Hence
$\cone(i_c):|\cone(\bs(\D^{m-1}))|\to I^m$ is a $PL$ homeomorphism linear on
the simplices of $|\cone(\bs(\D^{m-1}))|$. This defines the {\it canonical
triangulation\/} of~$I^m$. Thus, the canonical triangulation of
$I^m$ arises from the identification of $I^m$ with the cone over the
barycentric subdivision of~$\D^{m-1}$.
\end{construction}

\begin{construction}[cubical subdivision of a simple polytope]
\label{cubpol}
Let $P^n\subset\R^n$ be a simple polytope with $m$ facets
$F^{n-1}_1,\ldots,F^{n-1}_m$. Choose a point in the relative interior of
every face of $P^n$ (including the vertices and the polytope itself). We get
the set $\mathcal S$ of $1+f_0+f_1+\ldots+f_{n-1}$ points (here $\mb
f(P^n)=(f_0,f_1,\ldots,f_{n-1})$ is the $f$-vector of~$P^n$). For each vertex
$v\in P^n$ define the subset $\mathcal S_v\subset\mathcal S$ consisting of
the points chosen inside the faces containing~$v$. Since $P^n$ is simple, the
number of $k$-faces meeting at $v$ is $\binom nk$, $0\le k\le n$. Hence,
$\#\mathcal S_v=2^n$. The set $\mathcal S_v$ is the vertex set of an
$n$-cube, which we denote~$C_v^n$. The faces of $C^n_v$ can be described as
follows. Let $G^k_1$ and $G^l_2$ be two faces of $P^n$ such that $v\in
G^k_1\subset G^l_2$. Then there are exactly $2^{l-k}$ faces $G$ of $P^n$ such
that $G^k_1\subset G\subset G^l_2$. The corresponding $2^{l-k}$ points from
$\mathcal S$ form the vertex set of an $(l-k)$-face of~$C^n_v$. We denote this
face $C^{l-k}_{G_1\subset G_2}$. Every face of $C^n_v$ is
$C^i_{G_1\subset G_2}$ for some $G_1,G_2$ containing~$v$. The intersection of
any two cubes $C^n_v$, $C^n_{v'}$ is a face of each. Indeed, let $G^i\subset
P^n$ be the smallest face containing both vertices $v$ and~$v'$. Then
$C^n_v\cap C^n_{v'}=C^{n-i}_{G^i\subset P^n}$ is the face of both $I^n_v$
and~$I^n_{v'}$. Thus, we have constructed a cubical subdivision of $P^n$ with
$f_{n-1}(P^n)$ cubes of dimension~$n$. We denote this cubical complex
$\mathcal C(P^n)$.

There is an embedding of $\mathcal C(P^n)$ to $I^m$ constructed as follows.
Every $(n-k)$-face of $P^n$ is the intersection of $k$ facets:
$G^{n-k}=F^{n-1}_{i_1}\cap\ldots\cap F^{n-1}_{i_k}$. We map the corresponding
point of $\mathcal S$ to the vertex $(\e_1,\ldots,\e_m)\in I^m$, where
$\e_i=0$ if $i\in\{i_1,\ldots,i_k\}$ and $\e_i=1$ otherwise. This defines a
mapping from the vertex set $\mathcal S$ of $\mathcal C(P^n)$ to the vertex
set of~$I^m$. Using the canonical triangulation of $I^m$
(Construction~\ref{conbar}) we extend this mapping to the $PL$ embedding
$i_P:P^n\to I^m$. For each vertex
$v=F^{n-1}_{i_1}\cap\cdots\cap F^{n-1}_{i_n}\in P^n$ we have
\begin{equation}
\label{cubpolmap}
  i_P(C_v^n)=\bigl\{(y_1,\ldots,y_m)\in I^m\::\: y_j=1\text{ for }
  j\notin\{i_1,\ldots,i_n\}\bigr\},
\end{equation}
i.e. $i_P(C_v^n)=C_{\{i_1,\ldots,i_n\}}\subset I^m$ (in the notations
of~(\ref{ijface})). The embedding $i_P:P^n\to I^m$ for $n=2$, $m=3$ is shown
on Figure~2.
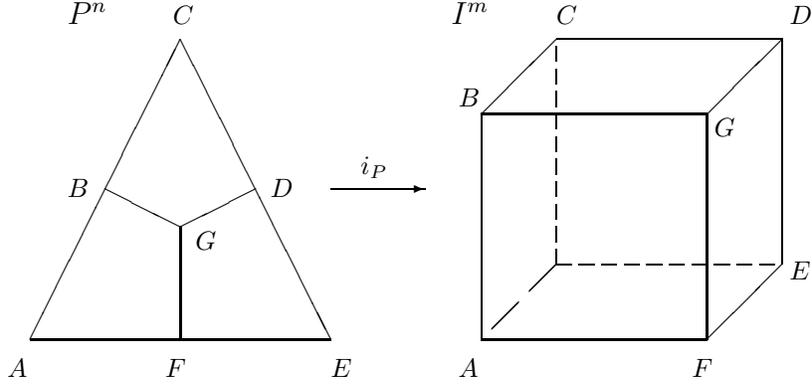
\begin{figure}
\begin{picture}(120,60)
  \put(10,10){\line(1,2){20}}
  \put(30,50){\line(1,-2){20}}
  \put(10,10){\line(1,0){40}}
  \put(20,30){\line(2,-1){10}}
  \put(40,30){\line(-2,-1){10}}
  \put(30,10){\line(0,1){15}}
  \put(50,30){\vector(1,0){12.5}}
  \put(70,10){\line(1,0){30}}
  \put(70,10){\line(0,1){30}}
  \put(70,40){\line(1,0){30}}
  \put(100,10){\line(0,1){30}}
  \put(70,40){\line(1,1){10}}
  \put(100,40){\line(1,1){10}}
  \put(100,10){\line(1,1){10}}
  \put(80,50){\line(1,0){30}}
  \put(110,20){\line(0,1){30}}
  \put(80,20){\line(-1,-1){3.6}}
  \put(75,15){\line(-1,-1){3.6}}
  \multiput(80,20)(3,0){10}{\line(1,0){2}}
  \multiput(80,20)(0,3){10}{\line(0,1){2}}
  \put(7,5){$A$}
  \put(28,5){$F$}
  \put(50,5){$E$}
  \put(32,22){$G$}
  \put(15,29){$B$}
  \put(42,29){$D$}
  \put(29,52){$C$}
  \put(15,52){\Large $P^n$}
  \put(54,32){$i_P$}
  \put(67,5){$A$}
  \put(98,5){$F$}
  \put(67,41){$B$}
  \put(101,37){$G$}
  \put(80,52){$C$}
  \put(111,52){$D$}
  \put(111,18){$E$}
  \put(66,52){\Large $I^m$}
\end{picture}
\caption{The embedding $i_P:P^n\to I^m$ for $n=2$, $m=3$.}
\end{figure}
\end{construction}

We summarise the facts from the above construction in the following statement

\begin{theorem}
\label{thcubpol}
  A simple polytope $P^n$ with $m$ facets can be split into cubes~$C^n_v$,
  one for each vertex $v\in P^n$. The resulting cubical complex
  $\mathcal C(P^n)$ embeds canonically
  into the boundary of~$I^m$, as described by~{\rm(\ref{cubpolmap})}.
\end{theorem}

\begin{lemma}
  The number of $k$-faces of the cubical complex $\mathcal C(P^n)$ is given by
  \begin{multline*}
    f_k\bigl(\mathcal C(P^n)\bigr)
    =\sum_{i=0}^{n-k}\bin{n-i}k f_{n-i-1}(P^n)\\
    =\bin nk f_{n-1}(P^n)+\bin{n-1}k f_{n-2}(P^n)+\dots+f_{k-1}(P^n),
    \quad k=0,\ldots,n.
  \end{multline*}
\end{lemma}
\begin{proof}
This follows from the fact that the $k$-faces of
$\mathcal C(P^n)$ are in one-to-one correspondence with pairs
$G_1^i,G_2^{i+k}$ of faces of $P^n$ such that $G_1^{i}\subset G_2^{i+k}$.
\end{proof}

\begin{construction}
\label{cubk}
Let $K^{n-1}$ be a simplicial complex on~$[m]$. Then $K$ is naturally a
subcomplex of~$\D^{m-1}$, and $\bs(K)$ is a subcomplex of $\bs(\D^{m-1})$. As
it follows from Construction~\ref{conbar}, there is a $PL$ embedding
$i_c|_{\bs(K)}:|\bs(K)|\to I^m$. The image $i_c(|\bs(K)|)$ is an
$(n-1)$-dimensional cubical subcomplex of~$I^m$,
which we denote $\cub(K)$. Then
\begin{equation}
\label{fcubk}
  \cub(K)=\bigcup_{\emptyset\ne I\subset J\in K}C_{I\subset
  J}\subset I^m,
\end{equation}
i.e. $\cub(K)$ is the union of faces $C_{I\subset J}\subset I^m$ over all
pairs $I\subset J$ of non-empty simplices of~$K$.
\end{construction}

\begin{construction}
\label{cck}
Since $\cone(\bs(K))$ is a subcomplex of $\cone(\bs(\D^{m-1}))$,
Construction~\ref{conbar} also provides a $PL$ embedding
$\cone(i_c)|_{\cone(\bs(K))}:|\cone(\bs(K))|\to I^m$. The image of this
embedding is an $n$-dimensional cubical subcomplex of~$I^m$, which we denote
$\cc(K)$. Then one easily obtains that
\begin{equation}
\label{fcck}
  \cc(K)=\bigcup_{J\in K}C_{I\subset J}\subset I^m.
\end{equation}
Since $C_{I\subset J}\subset C_{\emptyset\subset J}=C_J$, we also can write
$\cc(K)=\bigcup_{J\in K}C_J$.
\end{construction}

The following statement summarises the results of two previous constructions.
\begin{theorem}
\label{cubkcck}
  For any simplicial complex $K$ on the set $[m]$ there is a $PL$ embedding
  of the polyhedron $|K|$ into $I^m$ linear on the simplices of $\bs(K)$. The
  image of this embedding is the cubical subcomplex~{\rm(\ref{fcubk})}.
  Moreover, there is a $PL$ embedding of the polyhedron $|\cone(K)|$ into $I^m$
  linear on the simplices of $\cone(\bs(K))$. The image of this embedding is
  the cubical subcomplex~{\rm(\ref{fcck})}.
\end{theorem}

As in the case of simplicial complexes, a cubical complex $\mathcal C'$ is
called a {\it cubical subdivision\/} of cubical complex $\mathcal C$ if each
cube of $\mathcal C$ is a union of finitely many cubes of~$\mathcal C'$.

\begin{proposition}
  For every cubical subcomplex $\mathcal C$ there exists a cubical
  subdivision $\mathcal C'$ that can be realised as a subcomplex of
  some~$I^q$.
\end{proposition}
\begin{proof}
Subdividing each cube of $\mathcal C$ as described in
Construction~\ref{conbar} we obtain a simplicial complex, say~$K_{\mathcal
C}$. Then applying Construction~\ref{cubk} to $K_{\mathcal C}$ we get a
cubical complex that subdivides $|K_{\mathcal C}|=\mathcal C$ and embeds into
some $I^q$ (see Theorem~\ref{cubkcck}) as the subcomplex $\cub(K_{\mathcal
C})$.
\end{proof}

\begin{example}
\begin{figure}
  \begin{picture}(120,45)
  \put(15,5){\line(1,0){25}}
  \put(15,5){\line(0,1){25}}
  \put(40,5){\line(1,1){10}}
  \put(50,15){\line(0,1){25}}
  \put(50,40){\line(-1,0){25}}
  \put(25,40){\line(-1,-1){10}}
  \put(40,5){\circle*{2}}
  \put(15,30){\circle*{2}}
  \put(50,40){\circle*{2}}
  \put(40,30){\line(1,1){10}}
  \multiput(25,15)(5,0){5}{\line(1,0){3}}
  \multiput(25,15)(0,5){5}{\line(0,1){3}}
  \multiput(25,15)(-5,-5){2}{\line(-1,-1){4}}
  \put(26,16){$0$}
  \put(22,-2){(a)\ $K=\::\!\cdot$}
  \put(75,5){\line(1,0){25}}
  \put(75,5){\line(0,1){25}}
  \put(100,5){\line(1,1){10}}
  \put(110,15){\line(0,1){25}}
  \put(110,40){\line(-1,0){25}}
  \put(85,40){\line(-1,-1){10}}
  \put(100,5){\circle*{2}}
  \put(110,15){\circle*{2}}
  \put(75,30){\circle*{2}}
  \put(85,40){\circle*{2}}
  \put(75,5){\circle*{2}}
  \put(110,40){\circle*{2}}
  \put(100,30){\line(1,1){10}}
  \multiput(75,29.3)(0,0.1){16}{\line(1,1){10}}
  \multiput(100,4.3)(0,0.1){16}{\line(1,1){10}}
  \multiput(85,15)(5,0){5}{\line(1,0){3}}
  \multiput(85,15)(0,5){5}{\line(0,1){3}}
  \multiput(85,15)(-5,-5){2}{\line(-1,-1){4}}
  \put(86,16){$0$}
  \put(82,-2){(b)\ $K=\partial\D^2$}
  \put(15,30){\line(1,0){25}}
  \put(40,5){\line(0,1){25}}
  \put(75,30){\line(1,0){25}}
  \put(100,5){\line(0,1){25}}
  \linethickness{1mm}
  \put(75,5){\line(1,0){25}}
  \put(85,40){\line(1,0){25}}
  \put(75,5){\line(0,1){25}}
  \put(110,15){\line(0,1){25}}
  \end{picture}
  \caption{The cubical complex $\cub(K)$.}
  \end{figure}
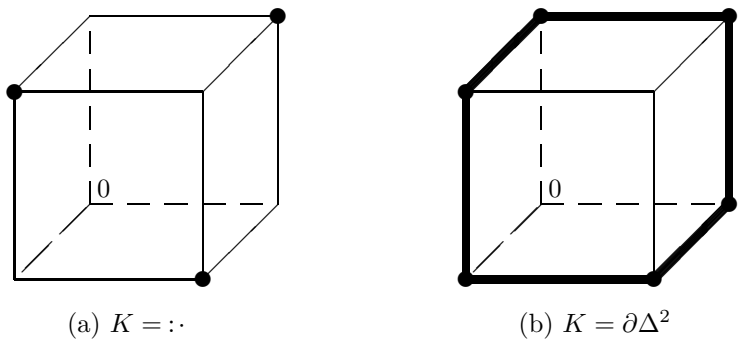
  \begin{figure}
  \begin{picture}(120,45)
  \put(15,5){\line(1,0){25}}
  \put(15,5){\line(0,1){25}}
  \put(40,5){\line(1,1){10}}
  \put(50,15){\line(0,1){25}}
  \put(50,40){\line(-1,0){25}}
  \put(25,40){\line(-1,-1){10}}
  \put(40,5){\circle*{2}}
  \put(40,30){\circle*{2}}
  \put(15,30){\circle*{2}}
  \put(50,40){\circle*{2}}
  \multiput(40,29.3)(0,0.1){16}{\line(1,1){10}}
  \multiput(25,15)(5,0){5}{\line(1,0){3}}
  \multiput(25,15)(0,5){5}{\line(0,1){3}}
  \multiput(25,15)(-5,-5){2}{\line(-1,-1){4}}
  \put(26,16){$0$}
  \put(22,-2){(a)\ $K=\::\!\cdot$}
  \put(75,5){\line(1,0){25}}
  \put(75,5){\line(0,1){25}}
  \put(100,5){\line(1,1){10}}
  \put(110,15){\line(0,1){25}}
  \put(110,40){\line(-1,0){25}}
  \put(85,40){\line(-1,-1){10}}
  \put(100,5){\circle*{2}}
  \put(110,15){\circle*{2}}
  \put(100,30){\circle*{2}}
  \put(75,30){\circle*{2}}
  \put(85,40){\circle*{2}}
  \put(75,5){\circle*{2}}
  \put(110,40){\circle*{2}}
  \multiput(100,29.3)(0,0.1){16}{\line(1,1){10}}
  \multiput(75,29.3)(0,0.1){16}{\line(1,1){10}}
  \multiput(100,4.3)(0,0.1){16}{\line(1,1){10}}
  \multiput(85,15)(5,0){5}{\line(1,0){3}}
  \multiput(85,15)(0,5){5}{\line(0,1){3}}
  \multiput(85,15)(-5,-5){2}{\line(-1,-1){4}}
  \multiput(78,5)(4,0){6}{\line(0,1){25}}
  \multiput(78,30)(4,0){6}{\line(1,1){10}}
  \multiput(100,8)(0,4){6}{\line(1,1){10}}
  \put(86,16){$0$}
  \put(82,-2){(b)\ $K=\partial\D^2$}
  \linethickness{1mm}
  \put(15,30){\line(1,0){25}}
  \put(40,5){\line(0,1){25}}
  \put(75,30){\line(1,0){25}}
  \put(75,5){\line(1,0){25}}
  \put(85,40){\line(1,0){25}}
  \put(100,5){\line(0,1){25}}
  \put(75,5){\line(0,1){25}}
  \put(110,15){\line(0,1){25}}
  \end{picture}
  \caption{The cubical complex $\cc(K)$.}
  \end{figure}
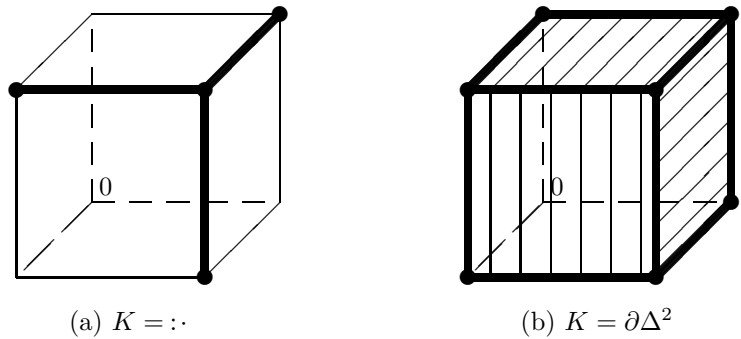
  Figure~3~(a) shows the cubical complex $\cub(K)$ for the case when $K$ is
  the disjoint union of 3 vertices ($K=\::\!\cdot$). Figure~3~(b) shows that
  for the case when $K$ is the boundary complex of a 2-simplex
  ($K=\partial\D^2$). The cubical complexes $\cc(K)$ in both cases are
  indicated on Figure~4~(a) and~(b).
\end{example}

\begin{remark}
As a topological space, $\cub(K)$ is homeomorphic to $|K|$, while
$\cc(K)$ is homeomorphic to $|\cone(K)|$. For simplicial complex $\cone(K)$
one may construct the cubical complex $\cub(\cone(K))$, which is also
homeomorphic to $|\cone(K)|$. However, as {\it cubical complexes\/}, $\cc(K)$
and $\cub(\cone(K))$ differ (since $\cone(\bs(K))\ne\bs(\cone(K))$).
\end{remark}

Let $P$ be a simple $n$-polytope and $K_P$ the corresponding simplicial
$(n-1)$-sphere (the boundary of the polar simplicial polytope~$P^*$). Then
$\cc(K_P)$ coincides with the cubical complex $\mathcal C(P)$ from
Construction~\ref{cubpol} (more precisely, $\cc(K_P)=i_P(\mathcal C(P))$).
Thus, the Construction~\ref{cubpol} is a particular case of
Construction~\ref{cck}.

\begin{remark}
Some versions of our previous constructions already appeared in the
literature. A version of Construction~\ref{cck} can be found
in~\cite[p.~434]{DJ} (it was used there for studying some torus actions; we
return to this in the next chapter). A version of our cubical subcomplex
$\cub(K)\subset I^m$ appeared in~\cite{SS} in connection with
Problem~\ref{novikov}.
\end{remark}

\smallskip

\section{Toric and quasitoric manifolds}
\subsection{Toric varieties}
\label{tori}
Toric varieties appeared in algebraic geometry in the beginning of 1970's in
connection with compactification problems for algebraic tori actions (see
below). Very quickly the geometry of toric varieties became one of the most
fascinating topics in algebraic geometry and found applications in many
mathematical sciences, which otherwise seem far from algebraic geometry.  We
have already mentioned the proof of necessity of $g$-theorem for simplicial
polytopes given by Stanley. Other remarkable applications include counting
lattice points and volumes of lattice polytopes; relations with Newton
polytopes and singularities (after Khovanskii and Kushnirenko);
discriminants, resultants and hypergeometric functions (after Gelfand,
Kapranov and Zelevinsky); reflexive polytopes and mirror symmetry for
Calabi--Yau toric hypersurfaces (after Batyrev). The standard references in
the toric geometry are Danilov's survey~\cite{Da} and books by
Fulton~\cite{Fu} and Oda~\cite{Od}. More recent survey article by
Cox~\cite{Co2} covers new applications, including mentioned above.  We are
not going to give a new survey of the toric geometry here. Instead, in this
section we stress upon some topological and combinatorial aspects of toric
varieties. We also give Stanley's argument for $g$-theorem.

Let $\C^*=\C\setminus\{0\}$ denote the multiplicative group of complex
numbers. The product $(\C^*)^n$ of $n$ copies of $\C^*$ is known as the torus
in the theory of algebraic groups. In topology, the {\it torus\/} $T^n$ is
the product of $n$ circles. We keep the topological notations, referring to
$(\C^*)^n$ as the {\it algebraic torus\/}. The torus $T^n$ is a subgroup of
the algebraic torus $(\C^{*})^n$ in the standard way:
\begin{equation}
\label{torus}
  T^n=\bigl\{\bigl(e^{2\pi i\f_1},\ldots,e^{2\pi i\f_n}\bigr)\in\C^n\bigr\},
\end{equation}
where $(\f_1,\ldots,\f_n)$ varies over~$\R^n$.

\begin{definition}
A {\it toric variety\/} is a normal algebraic variety $M$ containing the
algebraic torus $(\C^{*})^n$ as a Zariski open subset in such a way that the
natural action of $(\C^{*})^n$ on itself extends to an action on~$M$.
\end{definition}
Hence, $(\C^{*})^n$ acts on $M$ with a dense orbit.

Every toric variety is encoded by a set of combinatorial data, namely by a
(rational polyhedral) {\it fan\/} in some~$\R^n$.

Let $\R^n$ be the Euclidean space and $\Z^n\subset\R^n$ the integral lattice.
Given a finite set of vectors $\mb l_1,\ldots,\mb l_s\in\R^n$, define the {\it
convex polyhedral cone\/} $\sigma$ spanned by $\mb l_1,\ldots,\mb l_s$ as
$$
  \sigma=\{r_1\mb l_1+\dots+r_s\mb l_s\in\R^n\::\:r_i\ge0\}.
$$
Any convex polyhedral cone is a convex polyhedron in the sense of
Definition~\ref{pol2}. Hence, the {\it faces\/} of a convex polyhedral cone
are defined. A cone $\sigma$ is {\it rational\/} if its
generator vectors $\mb l_1,\ldots,\mb l_s$ can be taken from $\Z^n$ and is
{\it strongly convex\/} if it contains no line through the origin. All cones
considered below are strongly convex and rational.  A cone is called {\it
simplicial\/} (respectively {\it non-singular\/}) if it is generated by a
part of a basis of $\R^n$ (respectively~$\Z^n$). A {\it fan\/} is a set
$\Sigma$ of cones in some $\R^n$ such that each face of a cone in $\Sigma$ is
also a cone in~$\Sigma$, and the intersection of two cones in $\Sigma$ is a
face of each.  A fan $\Sigma$ is called {\it simplicial\/} (respectively {\it
non-singular\/}) if all cones of $\Sigma$ are simplicial (respectively
non-singular). A fan $\Sigma$ in $\R^n$ is called {\it complete\/} if the
union of all cones from $\Sigma$ is~$\R^n$.

A fan $\Sigma$ in $\R^n$ determines the toric variety $M_\Sigma$
of complex dimension~$n$, which orbit structure is described by the
combinatorics of~$\Sigma$. That is, the $k$-dimensional cones of $\Sigma$
correspond to the codimension-$k$ orbits of the algebraic torus action
on~$M^{2n}$. In particular, the $n$-dimensional cones correspond to the
fixed points, while the origin corresponds to the unique dense orbit. The
toric variety $M_\Sigma$ is compact if and only if $\Sigma$ is complete. If
$\Sigma$ is simplicial then $M_\Sigma$ is an {\it orbifold\/} (i.e. locally
homeomorphic to the quotient of $\R^{2n}$ by a finite group action). Finally,
if $\Sigma$ is non-singular, then, as one should expect, $M_\Sigma$ is
non-singular (smooth). Smooth toric varieties sometimes are called {\it toric
manifolds\/}.

Let $\Sigma$ be a simplicial fan in $\R^n$ with $m$ one-dimensional cones (or
{\it rays\/}). Choose generator vectors $\mb l_1,\ldots,\mb l_m$ for these $m$
rays to be integer and primitive (i.e. with relatively prime integer
coordinates). The fan $\Sigma$ defines the simplicial complex $K_\Sigma$ on
the vertex set~$[m]$, which is called the {\it underlying complex\/}
of~$\Sigma$.  By definition, $\{i_1,\ldots,i_k\}\subset[m]$ is a simplex of
$K_\Sigma$ if $\mb l_{i_1},\ldots,\mb l_{i_k}$ span a cone of~$\Sigma$. It is
easy to see that if $\Sigma$ is complete, then $K_\Sigma$ is a simplicial
$(n-1)$-sphere.

Denote $l_{ij}:=(\mb l_j)_i$, $1\le i\le n$, $1\le j\le m$. Here $\mb
l_j=(l_{1j},\ldots,l_{nj})^t\in\Z^n$. Assign to each vector $\mb l_j$ the
indeterminate $v_j$ of degree~2, and define linear forms
$$
  \t_i:=l_{i1}v_1+\dots+l_{im}v_m\in\Z[v_1,\ldots,v_m],\quad 1\le i\le n.
$$
Denote by $\mathcal J_\Sigma$ the ideal in $\Z[v_1,\ldots,v_m]$ spanned by
these linear forms, i.e. $\mathcal J_\Sigma=(\t_1,\ldots,\t_n)$. The images
of $\t_1,\ldots,\t_n$ and $\mathcal J_\Sigma$ in the Stanley--Reisner ring
$\Z(K_\Sigma)=\Z[v_1,\ldots,v_m]/\mathcal I_{K_\Sigma}$ (see
Definition~\ref{frsim}) will be denoted by the same symbols
$\t_1,\ldots,\t_n$ and~$\mathcal J_\Sigma$.

\begin{theorem}[Danilov and Jurkiewicz]
\label{danjur}
Let $\Sigma$ be a complete non-singular fan in~$\R^n$, and
$M_\Sigma$ the corresponding toric variety. Then

{\rm(a)} The Betti numbers (the ranks of homology groups) of $M_\Sigma$ vanish
in odd dimensions, while in even dimensions are given by
$$
  b_{2i}(M_\Sigma)=h_i(K_\Sigma),\quad i=0,1,\ldots,n,
$$
where $\mb h(K_\Sigma)=(h_0,\ldots,h_n)$ is the $h$-vector of~$K_\Sigma$.

{\rm(b)} The closures of orbits corresponding to the one-dimensional cones of
$\Sigma$ are codimension-{\rm2} submanifolds (divisors) $D_i$ of~$M_\Sigma$.
Let $v_i$, $1\le i\le m$, be the corresponding {\rm2}-dimensional cohomology
classes. Then the cohomology ring of $M_\Sigma$ is given by
$$
  H^*(M_\Sigma;\Z)\cong\Z[v_1,\ldots,v_m]/(\mathcal I_{K_\Sigma}+
  \mathcal J_\Sigma)=\Z(K_\Sigma)/\mathcal J_\Sigma.
$$
Moreover, $\t_1,\ldots,\t_n$ is a regular sequence in $\Z(K_\Sigma)$.
\end{theorem}
\noindent This theorem was proved by Jurkiewicz for projective smooth toric
varieties and by Danilov~\cite[Theorem~10.8]{Da} in the general case. Note
that the first part of Theorem~\ref{danjur} follows from the second part and
Lemma~\ref{psfr}.

Theorem~\ref{danjur} shows that the cohomology of $M_\Sigma$ is identified
with the {\it Chow ring\/}~\cite[\S~5.1]{Fu} of $M_\Sigma$ and is generated
by two-dimensional classes. Note that the ideal $\mathcal I_{K_\Sigma}$
depends only on simplicial complex $K_\Sigma$ (i.e. on the intersection
lattice of fan), while $\mathcal J_\Sigma$ depends on fan $\Sigma$ itself.

\begin{remark}
As it was shown by Danilov, Theorem~\ref{danjur} remains true for simplicial
fans and corresponding toric varieties if one replaces the coefficient ring
$\Z$ by any field of zero characteristic (e.g.~$\Q$).
\end{remark}

\begin{construction}[Normal fan and toric varieties arising from polytopes]
\label{nf}
Suppose we are given an $n$-polytope~(\ref{ptope}) with vertices in the
integer lattice $\Z^n\subset\R^n$. (Such polytopes are called integral (or
{\it lattice\/}).) Then the vectors $\mb l_i$ in~(\ref{ptope}), $1\le
i\le m$, can be chosen integer and primitive, and the numbers $a_i$ can be
chosen integer. Note that each $\mb l_i$ is the normal vector of the facet
$F_i\subset P^n$ and $\mb l_i$ is pointing inside the polytope~$P$. Define
the complete fan $\Sigma(P)$ whose cones are generated by sets of normal
vectors $\mb l_{i_1},\ldots,\mb l_{i_k}$ such that the corresponding facets
$F_{i_1},\ldots,F_{i_k}$ have non-empty intersection in~$P$. The fan
$\Sigma(P)$ is called the {\it normal fan\/} of~$P$. Alternatively, if $0\in
P$ then the normal fan consists of cones over the faces of the polar
polytope~$P^*$. Define the toric variety $M_P:=M_{\Sigma(P)}$. The variety
$M_P$ is smooth if and only if $P$ is simple and the normal vectors $\mb
l_{i_1},\ldots,\mb l_{i_n}$ of any set of $n$ facets
$F_{i_1},\ldots,F_{i_n}$ meeting at the same vertex form a basis of~$\Z^n$.
\end{construction}

\begin{remark}
Any combinatorial simple polytope is {\it rational\/}, that is, admits a
geometrical realisation with rational (or, equivalently, integer) vertex
coordinates. Indeed, there is a small perturbation of defining inequalities
in~(\ref{ptope}) that makes all of them rational but does not change the
combinatorial type (since the half-spaces defined by the inequalities are in
general position). As a result, one gets a
simple polytope $P'$ of the same combinatorial type with rational vertex
coordinates. To obtain a realisation with integral vertex coordinates one
should take the magnified polytope $kP'$ for appropriate $k\in\Z$. Quite
surprising thing is that there exist {\it non-rational\/} convex polytopes
(non-simple and non-simplicial), see~\cite[Example~6.21]{Zi}. Returning to
simple polytopes, we note that different realisations of a given
combinatorial simple polytope as lattice polytopes may produce different
(even topologically) toric varieties~$M_P$. At the same time there exist
combinatorial simple polytopes that do not admit {\it any\/} geometrical
realisation with smooth~$M_P$. We discuss one such example in the next
section (see Example~\ref{nonqt}).
\end{remark}

In the rest of this section all polytopes are assumed to be simple.
Construction~\ref{nf} allows to define the simplicial fan $\Sigma(P)$ and the
toric variety $M_P$ from any lattice simple polytope~$P$. However, the
polytope $P$ contains more information than the fan~$\Sigma(P)$. Indeed,
besides the normal vectors~$\mb l_i$, we also have numbers $a_i\in\Z$, $1\le
i\le m$, (see~(\ref{ptope})). The linear combination $D=a_1D_1+\dots+a_mD_m$
(see Theorem~\ref{danjur}) is an {\it ample divisor\/}. It defines a
projective embedding $M_P\subset\C P^r$ for some $r$ (which can be taken to
be the number of vertices of~$P$). Thus, all toric varieties from polytopes are
projective. Conversely, given a smooth projective toric variety $M\subset\C
P^r$, one gets a very ample divisor (line bundle) $D$ of a hyperplane section
whose zero cohomology is generated by the sections corresponding to lattice
points in a certain lattice simple polytope~$P$. For this $P$ one
has~$M=M_P$. Let $\omega:=a_1v_1+\dots+a_mv_m\in H^2(M_P;\Q)$ be the
cohomology class of~$D$.

\begin{theorem}[Hard Lefschetz theorem for toric varieties]
\label{hlth}
Let $P^n$ be a lattice simple polytope~{\rm(\ref{ptope})}, $M_P$ the
toric variety defined by~$P$, and $\omega=a_1v_1+\dots+a_mv_m\in
H^2(M_P;\Q)$. Then the maps
$$
\begin{CD}
  H^{n-i}(M_P;\Q) @>{\cdot\:\omega^i}>> H^{n+i}(M_P;\Q),\qquad 1\le i\le n,
\end{CD}
$$
are isomorphisms.
\end{theorem}
\noindent It follows from the projectivity that if $M_P$ is smooth then it is
K\"ahler, and $\omega$ is the class of the K\"ahler 2-form.

\begin{remark}
If one replaces the ordinary cohomology by the {\it intersection
cohomology\/}, then Theorem~\ref{hlth} holds for any projective toric
variety, not necessarily arising from a simple lattice polytope (see the
discussion in~\cite[\S~5.2]{Fu}).
\end{remark}

\begin{example}
\label{cpn}
The complex projective space $\C P^n=\{(z_0:z_1:\dots:z_n),z_i\in\C\}$ is a
toric variety. $(\C^*)^n$ acts on $\C P^n$ by $(t_1,\ldots,t_n)\cdot
(z_0:z_1:\dots:z_n)=(z_0:t_1z_1:\dots:t_nz_n)$. Obviously,
$(\C^*)^n\subset\C^n\subset\C P^n$ is a dense open subset. The fan defining
$\C P^n$ consists of cones generated by all proper subsets of $(n+1)$ vectors
$e_1,\ldots,e_n,-e_1-\dots-e_n$ in~$\R^n$. Theorem~\ref{danjur} identifies
the cohomology ring $H^*(\C P^n;\Z)=\Z[u]/(u^{n+1})$, $\dim u=2$, with the
quotient $\Z[v_1,\ldots,v_{n+1}]/(v_1\cdots
v_{n+1},v_1-v_{n+1},\ldots,v_n-v_{n+1})$. The toric variety $\C P^n$ arises
from a polytope: $\C P^n=M_P$, where $P$ is the standard
$n$-simplex~(\ref{stsim}). The class $\omega\in H^2(\C P^n;\Q)$ from
Theorem~\ref{hlth} in this case is $\omega=v_{n+1}$.
\end{example}

Now we are ready to give Stanley's argument for the necessity of
$g$-theorem for simple polytopes.
\begin{proof}[Proof of necessity of Theorem~{\rm\ref{gth}}]
Realise the simple polytope as a lattice polytope $P^n\subset\R^n$. Let
$M_P$ be the corresponding toric variety. Part~(a) is already proved
(Theorem~\ref{ds}). It follows from Theorem~\ref{hlth} that the
multiplication by $\omega\in H^2(M_P;\Q)$ is a monomorphism
$H^{2i-2}(M_P;\Q)\to H^{2i}(M_P;\Q)$ for $i\le \sbr n2$. This together with
the part~(b) of Theorem~\ref{danjur} gives $h_{i-1}\le h_i$, $0\le i\le
\sbr n2$, thus proving~(b). To prove~(c), define the graded commutative
$\Q$-algebra $A:=H^*(M_P;\Q)/(\omega)$. Then $A^0=\Q$,
$A^{2i}=H^{2i}(M_P;\Q)/\omega\cdot H^{2i-2}(M_P;\Q)$ for $1\le i\le\sbr n2$,
and $A$ is generated by degree-two elements (since so is $H^*(M_P;\Q)$). It
follows from Theorem~\ref{mvect} that the numbers $\dim A^{2i}=h_i-h_{i-1}$,
$0\le i\le\sbr n2$, are the components of an $M$-vector, thus proving~(c) and
the whole theorem.
\end{proof}

\begin{remark}
  The Dehn--Sommerville equations now follow from the Poincar\'e duality for
  $M_P$ (which holds also for singular $M_P$ provided that $P$ is simple).
\end{remark}

Now we consider the action of the torus $T^n\subset(\C^*)^n$ on a non-singular
compact toric variety~$M$. This action is locally equivalent to the standard
action of $T^n$ on $\C^n$ (see the next section for the precise definition).
The orbit space $M/T^n$ is homeomorphic to an $n$-ball, invested with the
topological structure of {\it manifold with corners\/} by the fixed point
sets of appropriate subtori, see~\cite[\S~4.1]{Fu}. (Roughly speaking, a
manifold with corners is a space that is locally modelled by open subsets
of the positive cone~$\R^n_+$~(\ref{pcone}). From this description it is
easy to deduce the strict definition, see~\cite{Ja}, which we omit here.)

\begin{construction}
\label{corners}
Let $P^n$ be a simple polytope. For any vertex $v\in P^n$ denote by $U_v$ the
open subset of $P^n$ obtained by deleting all faces not containing~$v$.
Obviously, $U_v$ is diffeomorphic to $\R^n_+$ (and even affinely isomorphic
to an open set of $\R^n_+$ containing~0). It follows that $P^n$ is a manifold
with corners, with atlas $\{U_v\}$.
\end{construction}

If $M=M_P$ arises from some (simple) polytope~$P^n$, then the orbit space
$M/T^n$ is diffeomorphic, as a manifold with corners, to~$P^n$. Furthermore,
there exists an explicit map $M_P\to\R^n$ (the {\it moment map\/}) with image
$P^n$ and $T^n$-orbits as fibres, see~\cite[\S4.2]{Fu}. (We consider
relationships with moment maps and some aspects of symplectic geometry in
more details in section~\ref{coor}.) Under this map, the interior of a
codimension-$k$ face of $P^n$ is identified with the set of orbits having the
same $k$-dimensional isotropy subgroup. In particular, the action is free
over the interior of the polytope. Regarded as a smooth manifold, $M_P$ can
be identified with quotient space $T^n\times P^n/{\sim}$ for some equivalence
relation~$\sim$. Such description of the torus action on a non-singular toric
variety motivated the appearance of an important topological analogue of
toric varieties, namely, the theory of {\it quasitoric manifolds\/}.

\subsection{Quasitoric manifolds}
\label{quas}
A quasitoric manifold is a smooth manifold with a torus action whose
properties are similar to that of the (compact) torus action on a
non-singular projective toric variety. This notion appeared in~\cite{DJ}
under the name ``toric manifold". We use the term ``quasitoric manifold",
since ``toric manifold" is occupied in the algebraic geometry for
``non-singular toric variety". In the consequent definitions we
follow~\cite{DJ}, taking into account some adjustments from~\cite{BR2}.

As in the previous section, we regard the torus $T^n$ as the standard
subgroup~(\ref{torus}) in~$\C^*$, therefore specifying the orientation and
the coordinate subgroups $T_i\cong S^1$, $i=1,\ldots,n$, in~$T^n$. We refer
to the representation of $T^n$ by diagonal matrices in $U(n)$ as the {\it
standard\/} action on~$\C^n$. The orbit space of this action is the positive
cone~$\R^n_+$. The canonical projection
$$
  T^n\times\R^n_+\to\C^n\::\:(t_1,\ldots,t_n)\times(x_1,\ldots,x_n)\to
  (t_1x_1,\ldots,t_nx_n)
$$
allows to identify $\C^n$ with the quotient space
$T^n\times \R^n_+/{\sim}$ for some equivalence relation~$\sim$, which will
play an important r\^ole in our future considerations.

Let $M^{2n}$ be a $2n$-dimensional manifold with an action of the
torus~$T^n$. Say that the $T^n$-action is {\it locally standard\/} if every
point $x\in M^{2n}$ lies in some $T^n$-invariant neighbourhood $U(x)$ for
which there exists a $\psi$-equivariant homeomorphism $f:U(x)\to W$ with
some ($T^n$-stable) open subset $W\subset\C^n$. That is, there is an
automorphism $\psi:T^n\to T^n$ such that $f(t\cdot y)=\psi(t)f(y)$ for
all $t\in T^n$, $y\in U(x)$. The orbit space for a locally standard action of
$T^n$ on $M^{2n}$ is an $n$-dimensional manifold with corners. Quasitoric
manifolds correspond to the important case when this orbit space is
diffeomorphic, as manifold with corners, to a simple polytope~$P^n$.

\begin{definition}
\label{qtm}
Given a simple polytope $P^n$, a manifold $M^{2n}$ with locally standard
$T^n$-action is said to be a {\it quasitoric manifold over\/} $P^n$ if there
is a projection map $\pi:M^{2n}\to P^n$ whose fibres are the orbits of the
action.
\end{definition}

\noindent Under projection $\pi$, points that have same isotropy subgroup of
codimension $k$ are taken to the interior of a certain $k$-face of~$P^n$. In
particular, the action of $T^n$ is free over the interior~of $P^n$, while the
vertices of $P^n$ correspond to the $T^n$-fixed points of~$M^{2n}$.

\begin{remark}
Two simple polytopes are combinatorially equivalent if and only if they are
diffeomorphic as manifolds with corners.
\end{remark}

Suppose $P^n$ has $m$ facets $F_1,\ldots,F_m$. For every facet $F_i$ the
pre-image $\pi^{-1}(F_i)$ is a submanifold $M_i^{2(n-1)}\subset M^{2n}$ with
isotropy subgroup a circle $T(F_i)$ in~$T^n$. Then
\begin{equation}
\label{fisotr}
  T(F_i)=\bigl\{\bigl(e^{2\pi i\l_{1i}\f},\ldots,e^{2\pi
  i\l_{ni}\f}\bigr)\in T^n\bigr\},
\end{equation}
where $\f\in\R$ and $\bl_i=(\l_{1i},\ldots,\l_{ni})^t\in\Z^n$ is a primitive
vector. This $\bl_i$ is determined by $T(F_i)$ only up to a sign. The choice
of this sign specifies an orientation of~$T(F_i)$. For now we
do not care about this sign and choose it arbitrary. More detailed treatment
of these signs is the subject of the next section. We refer to $\bl_i$ as
the {\it facet vector\/} corresponding to~$F_i$.
The action of $T^n/T(F_i)$ on $M_i$ describes it as a quasitoric manifold
over~$F_i$. The correspondence
\begin{equation}
\label{charf}
  \ell:F_i\mapsto T(F_i)
\end{equation}
is called the {\it characteristic map\/} of~$M^{2n}$. Suppose we have a
codimension-$k$ face $G^{n-k}$ written as the intersection of $k$ facets:
$G^{n-k}=F_{i_1}\cap\dots\cap F_{i_k}$. Then the submanifolds
$M_{i_1},\ldots,M_{i_k}$ intersect transversally in a {\it facial
submanifold\/} $M(G)^{2(n-k)}$. The map $T(F_{i_1})\times\dots\times
T(F_{i_k})\to T^n$ is injective since $T(F_{i_1})\times\dots\times
T(F_{i_k})$ is identified with the $k$-dimensional isotropy subgroup of
$M(G)^{2(n-k)}$. It follows that the vectors $\bl_{i_1},\ldots,\bl_{i_k}$
form a part of integral basis of~$\Z^n$.

Let $\Lambda$ be the integer $(n\times m)$-matrix whose $i$-th column is
formed by the coordinates of the facet vector~$\bl_i$, $i=1,\ldots,m$. Each
vertex $v\in P^n$ is the in intersection of $n$ facets:
$v=F_{i_1}\cap\cdots\cap F_{i_n}$. Let
$\Lambda_{(v)}:=\Lambda_{(i_1,\ldots,i_n)}$ be the maximal minor of $\Lambda$
formed by the columns $i_1,\ldots,i_n$. Then
\begin{equation}
\label{L}
  \det\Lambda_{(v)}=\pm1.
\end{equation}

The correspondence
$$
  G^{n-k}\mapsto\text{ isotropy subgroup of }M(G)^{2(n-k)}
$$
extends the characteristic map~(\ref{charf}) to a map from the face lattice
of $P^n$ to the lattice of subtori of~$T^n$.

Like in the case of standard action of $T^n$ on~$\C^n$, there is a projection
$T^n\times P^n\to M^{2n}$ whose fibre over $x\in M^{2n}$ has the form
(isotropy subgroup of $x)\times($orbit of~$x$). This argument can be used for
reconstructing the quasitoric manifold from any given pair $(P^n,\ell)$, where
$P^n$ is a (combinatorial) simple polytope and $\ell$ is a map from
facets of $P^n$ to one-dimensional subgroups of $T^n$ such that
$\ell(F_{i_1})\times\dots\times\ell(F_{i_k})\to T^n$ is injective whenever
$F_{i_1}\cap\cdots\cap F_{i_k}\ne\emptyset$. Such $(P^n,\ell)$ is
called a {\it characteristic pair\/}. The map $\ell$ directly extends to a map
from the face lattice of $P^n$ to the lattice of subtori of~$T^n$.

\begin{construction}[Quasitoric manifold from a characteristic pair]
\label{der}
Note that each point $q$ of $P^n$ lies in the relative interior of a unique
face~$G(q)$. Now construct the identification space $(T^n\times P^n)/{\sim}$,
where $(t_1,q)\sim(t_2,q)$ if and only if $t_1t_2^{-1}$ lies in the subtorus
$\ell(G(q))$. The free action of $T^n$ on $T^n\times P^n$ obviously descends
to an action on $(T^n\times P^n)/{\sim}$ with quotient~$P^n$. The latter
action is free over the interior of $P^n$ and has a fixed point for each
vertex of~$P^n$. Just as $P^n$ is covered by the open sets $U_v$, based on
the vertices and diffeomorphic to $\R^n_+$ (see Construction~\ref{corners}),
so the space $(T^n\times P^n)/{\sim}$ is covered by open sets $(T^n\times
U_v)/{\sim}$ homeomorphic to $(T^n\times\R^n_+)/{\sim}$, and therefore
to~$\C^n$. This implies that the $T^n$-action on $(T^n\times P^n)/{\sim}$ is
locally standard, and therefore $M^{2n}(\ell):=(T^n\times P^n)/{\sim}$ is a
quasitoric manifold.
\end{construction}

Given an automorphism $\psi:T^n\to T^n$, say that two quasitoric manifolds
$M^{2n}_1$, $M^{2n}_2$ over the same $P^n$ are {\it $\psi$-equivariantly
diffeomorphic\/} if there is a diffeomorphism $f:M^{2n}_1\to M^{2n}_2$
such that $f(t\cdot x)=\psi(t)f(x)$ for all $t\in T^n$, $x\in M^{2n}_1$.
The automorphism $\psi$ induces an automorphism $\psi_*$ of
the lattice of subtori of~$T^n$. Any such automorphism descends to a
{\it $\psi$-translation\/} of characteristic pairs, in which the two
characteristic maps differ by~$\psi_*$. The following proposition
is proved in~\cite[Proposition~2.6]{BR2} and generalises Proposition~1.8
of~\cite{DJ}.

\begin{proposition}\label{equivar}
For any automorphism $\psi$, Construction~{\rm\ref{der}} defines
a bijection between $\psi$-equivariant diffeomorphism classes of quasitoric
manifolds and $\psi$-translations of pairs $(P^n,\ell)$.
\end{proposition}

When $\psi$ is the identity, we deduce that two quasitoric manifolds are
equivariantly diffeomorphic if and only if their characteristic maps are the
same.

Now we are going to construct a cellular decomposition of $M^{2n}$ with only
even dimensional cells and calculate the Betti numbers accordingly,
following~\cite{DJ}.

\begin{construction}
We recall ``Morse-theoretical arguments" from the proof of Dehn--Sommerville
relations (Theorem~\ref{ds}). There we turned the 1-skeleton of $P^n$ into a
directed graph and defined the index $\ind(v)$ of a vertex $v\in P^n$ as the
number of incident edges that point towards~$v$. These inward edges define a
face $G_v$ of dimension $\ind(v)$. Denote by $\widehat{G}_v$ the subset of
$G_v$ obtained by deleting all faces not containing~$v$. Obviously,
$\widehat{G}_v$ is diffeomorphic to $\R^{\ind(v)}_+$ and is contained in
the open set $U_v\subset P^n$ from Construction~\ref{corners}. Then
$e_v:=\pi^{-1}\widehat{G}_v$ is identified with $\C^{\ind(v)}$, and the union
of the $e_v$ over all vertices of $P^n$ define a cell decomposition
of~$M^{2n}$. Note that all cells are even-dimensional and the closure of
the cell $e_v$ is the facial submanifold $M(G_v)^{2\ind(v)}\subset M^{2n}$.
This construction was earlier used by Khovanskii~\cite{Kh} for constructing
cellular decompositions of toric varieties.
\end{construction}

\begin{proposition}
\label{qtbn}
The Betti numbers of $M^{2n}$ vanish
in odd dimensions, while in even dimensions are given by
$$
  b_{2i}(M^{2n})=h_i(P^n),\quad i=0,1,\ldots,n,
$$
where $\mb h(P^n)=(h_0,\ldots,h_n)$ is the $h$-vector of~$P^n$.
\end{proposition}
\begin{proof}
The $2i$-th Betti number equals the number of $2i$-dimensional cells in the
cellular decomposition constructed above. This number equals the number of
vertices of index $i$, which is $h_i(P^n)$ by the argument from the proof of
Theorem~\ref{ds}.
\end{proof}

Given a quasitoric manifold $M^{2n}$ with characteristic map~(\ref{charf})
and facet vectors $\bl_i=(\l_{1i},\ldots,\l_{ni})^t\in\Z^n$, $i=1,\ldots,m$,
define linear forms
\begin{equation}
\label{theta}
  \t_i:=\l_{i1}v_1+\dots+\l_{im}v_m\in\Z[v_1,\ldots,v_m],\quad 1\le i\le n.
\end{equation}
The images of these linear forms in the face ring $\Z(P^n)$ will be denoted
by the same letters.

\begin{lemma}[Davis and Januszkiewicz]
\label{crs}
  For any quasitoric manifold $M^{2n}$ over~$P^n$, the sequence
  $\t_1,\ldots,\t_n$ is a (degree-two) regular sequence in $\Z(P^n)$.
\end{lemma}

Let $\mathcal J_\ell$ denote the ideal in $\Z(P^n)$ generated by
$\t_1,\ldots,\t_n$.

\begin{theorem}[Davis and Januszkiewicz]
\label{qtcoh}
Let $v_i$, $1\le i\le m$, be the {\rm2}-dimensional cohomology classes
dual to the submanifolds $M_i^{2(n-1)}\subset M^{2n}$.
Then the cohomology ring of $M^{2n}$ is given by
$$
  H^*(M^{2n};\Z)\cong\Z[v_1,\ldots,v_m]/(\mathcal I_P+
  \mathcal J_\ell)=\Z(P^n)/\mathcal J_\ell.
$$
\end{theorem}

We give proofs for the above two statements in section~\ref{hom1}.

\begin{remark}
Change of sign of vector $\bl_i$ corresponds to passing from $v_i$ to $-v_i$
in the description of the cohomology ring from Theorem~\ref{qtcoh}. This fact
will be crucial in the next section.
\end{remark}

\begin{example}
\label{tcf}
A non-singular projective toric variety $M_P$ arising from a lattice simple
polytope $P^n$ is a quasitoric manifold over~$P^n$. The corresponding
characteristic map $\ell:F_i\mapsto T(F_i)$ is defined by~(\ref{fisotr}),
where $\bl_i=(\l_{1i},\ldots,\l_{ni})^t$ are the normal vectors $\mb l_i$ of
facets of $P^n$, $i=1,\ldots,m$ (see~(\ref{ptope})). The corresponding
characteristic $n\times m$-matrix $\L$ is the matrix $L$ from
Construction~\ref{dist}. In particular, if $P^n$ is the standard simplex
$\D^n$~(\ref{stsim}) then $M_P$ is $\C P^n$ (Example~\ref{cpn}) and
$\L=(\mathrm E\,|-\!\mathbf 1)$, where $\mathrm E$ is the unit $n\times
n$-matrix and $\mathbf 1$ is the column of units. See also
Example~\ref{cpnoo} below.
\end{example}

Generally, a non-singular toric variety is not necessarily a quasitoric
manifold: although the orbit space (for the action of~$T^n$) is a manifold
with corners (see section~\ref{tori}), it may fail to be diffeomorphic
(or combinatorially equivalent) to a simple polytope. The authors are
thankful to N.~Strickland for drawing our attention to this fact. However, we
do not know any such example. In~\cite[p.~71]{Fu} one can find the example of
a complete non-singular fan $\Sigma$ in $\R^3$ which can not be obtained by
taking the cones with vertex~0 over the faces of a geometrical simplicial
polytope. Nevertheless, since the corresponding simplicial complex $K_\Sigma$
is a simplicial 2-sphere, it is combinatorially equivalent to a polytopal
2-sphere. This means that the non-singular toric variety~$M_\Sigma$, although
being non-projective, is still a quasitoric manifold.

\begin{problem}
Give an example of a non-singular toric variety which is not a quasitoric
manifold.
\end{problem}

This problem seems to be not very hard and reduces to constructing a complete
non-singular fan~$\Sigma$ whose associated simplicial complex $K_\Sigma$ is
the Barnette sphere or any other non-polytopal $PL$-sphere (see
section~\ref{sim1}). As it was pointed out by N.~Strickland, the definition
of quasitoric manifold can be modified in such a way that non-singular toric
varieties become quasitoric manifolds. For this purpose, the simple polytope
in Definition~\ref{qtm} is to be replaced by the polytopal complex dual to a
simplicial sphere. The corresponding constructions are currently being
developed.

On the other hand, it is easy to construct a quasitoric manifold which is not
a toric variety. The simplest example is the manifold $\C P^2\cs\C P^2$ (the
connected sum of two copies of~$\C P^2$). It is a quasitoric manifold over
the square $I^2$ (this follows from the construction of equivariant connected
sum, see~\cite[1.11]{DJ}, section~\ref{stab} and corollary~\ref{4dtopclas}
below). However, $\C P^2\cs\C P^2$ do not admit even an almost complex
structure (i.e. a complex structure in the tangent bundle). The following
problem arises.

\begin{problem}\label{acs}
Let $P^n$ be simple polytope with $m$ facets, $\ell$ a characteristic
map~{\rm(\ref{charf})}, and $M^{2n}(\ell)$ the derived quasitoric
manifold (Construction~{\rm\ref{der}}). Find conditions on $P^n$ and $\ell$
so that $M^{2n}(\ell)$ admits a $T^n$-invariant complex (respectively almost
complex) structure.
\end{problem}

The almost complex case of the above problem was formulated
in~\cite[Problem~7.6]{DJ}. Example~\ref{tcf} (characteristic functions of
lattice polytopes) provides a sufficient condition for Problem~\ref{acs}
(since a non-singular toric variety is a complex manifold). However, this
condition is obviously non-necessary even for the existence of a complex
structure.  Indeed, there are non-singular (non-projective) toric varieties
that do not arise from any lattice simple polytope (see the already mentioned
example from~\cite[p.~71]{Fu}). At the same time, we do not know any example
of {\it non-toric\/} complex quasitoric manifold.

\begin{problem}\label{ntcs}
Find an example of a non-toric quasitoric manifold that admits a
$T^n$-invariant complex structure.
\end{problem}

Although a general quasitoric manifold may fail to be complex or almost
complex, it always admits a $T^n$-invariant complex structure on the {\it
stable\/} tangent bundle. The corresponding constructions are the subject of
the next section.

Another class of problems arises in connection with the classification of
quasitoric manifolds over a given combinatorial simple polytope. The general
setting of this problem is the subject of section~\ref{prob}.
Example~\ref{nonqt} below shows that there are combinatorial simple polytopes
that do not admit a characteristic map (and therefore can not arise as orbit
spaces for quasitoric manifolds).

\begin{problem}\label{ecf}
Give a combinatorial description of the class of polytopes $P^n$ that
admit a characteristic map~{\rm(\ref{charf})}.
\end{problem}
A generalisation of this problem is given in section~\ref{part}
(Problem~\ref{sp}).

A characteristic map is determined by an integer $n\times m$-matrix~$\L$
which satisfies~(\ref{L}) for every vertex $v\in P^n$. The equation
$(\det\L_{(v)})^2=1$ defines a hypersurface in the space $\mathcal M(n,m;\Z)$
of integer $n\times m$-matrices. Problem~\ref{ecf} can be reformulated in the
following way.

\begin{proposition}
\label{intersect}
The set of characteristic matrices coincides with the intersection
\begin{equation}\label{fintersect}
  \bigcap_{v\in P^n}\bigl\{(\det\L_{(v)})^2=1\bigr\}
\end{equation}
of hypersurfaces in the space $\mathcal M(n,m;\Z)$, where $v$ varies over the
vertices of the polytope~$P^n$.
\end{proposition}

\begin{example}[{polytope that do not admit a characteristic
function,~\cite[Example~1.22]{DJ}}]
\label{nonqt}
Suppose $P^n$ is a 2-neighbourly simple polytope with $m\ge2^n$ facets (e.g.
the polar of the cyclic polytope $C^n(m)$ (Example~\ref{cyclic}) with $n\ge4$
and $m\ge2^n$). Then $P^n$ does not admit a characteristic map, and therefore
can not appear as the quotient space for a quasitoric manifold.  Indeed, by
Proposition~\ref{intersect}, it is sufficient to show that
intersection~(\ref{fintersect}) is empty. Since $m\ge2^n$, any matrix
$\L\in\mathcal M(n,m;\Z)$ (without zero columns) contains two columns, say
$i$-th and $j$-th, which coincide modulo~2. Since $P^n$ is 2-neighbourly, the
corresponding facets $F_i$ and $F_j$ have non-empty intersection in~$P^n$.
Hence, columns $i$ and $j$ of~$\L$ enter some minor of the form~$\L_{(v)}$.
This implies that the determinant of this minor is even and
intersection~(\ref{fintersect}) is empty.
\end{example}

\subsection{Stably complex structures and quasitoric representatives in
cobordism classes}\label{stab}
This section is the review of results obtained by N.~Ray and the first
author in~\cite{BR1} and~\cite{BR2} supplied with some additional comments.

A {\it stably complex structure\/} on a (smooth) manifold $M$ is defined by a
complex structure on the vector bundle $\tau(M)\oplus\R^{k}$ for some~$k$,
where $\tau(M)$ is the {\it tangent bundle\/} of~$M$ and $\R^k$ denotes the
trivial real $k$-bundle over~$M$. A {\it stably complex manifold\/} (in other
notations, {\it weakly almost complex manifold\/} or $U$-{\it manifold\/}) is
a manifold with fixed stably complex structure, that is, a pair $(M,\xi)$,
where $\xi$ is a complex bundle isomorphic, as a real bundle, to
$\tau(M)\oplus\R^{k}$ for some~$k$. If $M$ itself is a complex manifold, then
it possesses the {\it canonical\/} stably complex structure $(M,\tau(M))$.
The operations of disjoint union and product endow the set of cobordism
classes $[M,\xi]$ of stably complex manifolds with the structure of a graded
ring, called the {\it complex cobordism ring\/}~$\Omega^U$. By the theorem of
Milnor and Novikov, $\O^U\cong\Z[a_1,a_2,\ldots]$, $\deg a_i=2i$
(see~\cite{No1},~\cite{St}). The ring $\O^U$ serves as the coefficient ring
for a generalised (co)homology theory known as the {\it complex
(co)bordisms\/}.

Stably complex manifolds were the main subject of F.~Hirzebruch's talk at the
1958 International Congress of mathematicians (see~\cite{icm58}). Using {\it
Milnor hypersurfaces\/} (Example~\ref{milnor}) and the Milnor--Novikov
theorem it was shown that every complex cobordism class contains a
non-singular algebraic variety (not necessarily connected).  The following
problem is still open.

\begin{problem}[Hirzebruch]
\label{hirz}
Which cobordism classes in $\O^U$ contain {\it connected\/} non-singular
algebraic varieties?
\end{problem}

\begin{example}
The 2-dimensional cobordism group $\O^U_2\cong\Z$ is generated by the class
of $[\C P^1]$ (Riemannian sphere). Every cobordism class $k[\C P^1]\in\O^U_2$
contains a non-singular algebraic variety, namely, the disjoint union of $k$
copies of $\C P^1$ for $k>0$ and the disjoint union of $k$ copies of a genus
2 Riemannian surface for~$k<0$.  However, a {\it connected\/} algebraic
variety is contained only in cobordism classes $k[\C P^1]$ for~$k\le1$.
\end{example}

The problem of choice of appropriate generators for the ring $\O^U$ plays a
pivotal r\^ole in the cobordism theory and its applications. In this section
we give a solution of the quasitoric analogue of problem~\ref{hirz}, recently
obtained in~\cite{BR1} and~\cite{BR2}. This solution relies upon an important
additional structure on a quasitoric manifold, namely, the
{\it omniorientation\/}, which provides a canonical stably complex structure
described in the combinatorial terms.

Let $\pi:M^{2n}\to P^n$ be a quasitoric manifold with characteristic
map~$\ell$. Since the torus $T^n$~(\ref{torus}) is oriented, a choice of
orientation for $P^n$ is equivalent to a choice of orientation for~$M^{2n}$.
(An orientation of $P^n$ is specified by orienting the ambient space~$\R^n$.)

\begin{definition}\label{omnior}
An {\it omniorientation\/} of a quasitoric manifold $M^{2n}$ consists of a
choice of an orientation for $M^{2n}$ and for every submanifold
$M_i^{2(n-1)}=\pi^{-1}(F_i)$, $i=1,\ldots,m$.
\end{definition}
There are therefore $2^{m+1}$ omniorientations in all for given~$M^{2n}$.

An omniorientation of $M^{2n}$ determines an orientation for every normal
bundle $\nu_i:=\nu(M_i\subset M^{2n})$, $i=1,\ldots,m$. Since every
$\nu_i$ is a real 2-plane bundle, an orientation of $\nu_i$ allows to
interpret it as a complex line bundle. The isotropy subgroup $T(F_i)$ of the
submanifold $M_i^{2(n-1)}=\pi^{-1}(F_i)$ acts on the normal bundle
$\nu_i$, $i=1,\ldots,m$. Thus, we have the following statement.

\begin{proposition}\label{oosigns}
A choice of omniorientation for $M^{2n}$ is equivalent to a choice of
orientation for $P^n$ together with an unambiguous choice of
facet vectors $\bl_i$, $i=1,\ldots,m$ in~{\rm(\ref{fisotr})}.
\end{proposition}

We refer to a characteristic map $\ell$ as {\it directed\/} if all circles
$\ell(F_i)$, $i=1,\ldots,m$, are oriented. This implies that
signs of the facet vectors $\bl_i=(\l_{1i},\ldots,\l_{ni})^t$,
$i=1,\ldots,m$, are determined unambiguously. In the previous section we
organised the facet vectors into the integer $n\times m$ matrix~$\Lambda$.
This matrix satisfies~(\ref{L}). Due to~(\ref{fisotr}), knowing a matrix
$\Lambda$ is equivalent to knowing a directed characteristic map. Let $\Z^\F$
denote the $m$-dimensional free $\Z$-module spanned by the set $\F$ of facets
of~$P^n$. Then $\Lambda$ defines an epimorphism $\l:\Z^\F\to\Z^n$ by
$\l(F_i)=\bl_i$ and an epimorphism $T^\F\to T^m$, which we denote by the same
letter~$\l$. In the sequel we write $\Z^m$ for $\Z^\F$ and $T^m$ for $T^\F$
with the agreement that the vector $\mb e_i$ of the standard basis of $\Z^m$
corresponds to $F_i\in\Z^\F$, $i=1,\ldots,m$ (and the same for~$T^m$).
A {\it directed characteristic pair\/} $(P^n,\Lambda)$ consists of a simple
polytope $P^n$ and an integer matrix $\Lambda$ (or, equivalently, an
epimorphism $\l:\Z^m\to\Z^n$) that satisfies~(\ref{L}).

Proposition~\ref{oosigns} shows that the characteristic pair of an
omnioriented quasitoric manifold is directed. On the other hand, the
quasitoric manifold derived from a directed characteristic pair by
Construction~\ref{der} is omnioriented.

\begin{construction}
The orientation of the bundle $\nu_i$ over $M_i$ defines an integral Thom
class in the cohomology group $H^2(M(\nu_i))$, represented by a complex line
bundle over the Thom complex $M(\nu_i)$. We pull this back along the
Pontryagin--Thom collapse $M^{2n}\to M(\nu_i)$, and denote the resulting
bundle~$\rho_i$. The restriction of $\rho_i$ to $M_i\subset M^{2n}$ is~$\nu_i$.
\end{construction}

\begin{theorem}[{\cite[Theorem~3.8]{BR2}}]
\label{stcomst}
  Every omniorientation of a quasitoric manifold $M^{2n}$ determines a stably
  complex structure on it by means of the following isomorphism of real
  $2m$-bundles:
  $$
    \tau(M^{2n})\oplus\R^{2(m-n)}\cong\rho_1\oplus\cdots\oplus\rho_m.
  $$
\end{theorem}
It follows that a directed characteristic pair $(P^n,\Lambda)$ determines a
complex cobordism class $[M^{2n},\rho_1\oplus\cdots\oplus\rho_m]\in\O^U$. At
the same time, the above constructions can be directly applied to computing
the complex cobordism ring $\O_U^*(M^{2n})$ of an omnioriented quasitoric
manifold.

\begin{theorem}[{\cite[Proposition~5.3]{BR2}}]
\label{qtcob}
Let $v_i$ denote the first cobordism Chern class $c_1(\rho_i)\in
\O_U^2(M^{2n})$ of the bundle $\rho_i$, $1\le i\le m$.
Then the cobordism ring of $M^{2n}$ is given by
$$
  \O_U^*(M^{2n})\cong\O^U[v_1,\ldots,v_m]/(\mathcal I_P+
  \mathcal J_\Lambda),
$$
where the ideals $\mathcal I_P$ and $\mathcal J_\Lambda$ are defined in the
same way as in Theorem~{\rm\ref{qtcoh}}.
\end{theorem}
\noindent Note that the Chern class $c_1(\rho_i)$ is Poincar\'e dual to the
inclusion $M_i^{2(n-1)}\subset M^{2n}$ by construction of~$\rho_i$. This
highlights the remarkable fact that the complex bordism groups
$\O_*^U(M^{2n})$ are spanned by {\it embedded\/} submanifolds. By definition,
the fundamental cobordism class $\langle M^{2n}\rangle\in\O_U^{2n}(M^{2n})$
is dual to the bordism class of a point.  Thus, $\langle
M^{2n}\rangle=v_{i_1}\cdots v_{i_n}$ for any set $i_1,\ldots,i_n$ such that
$F_{i_1}\cap\dots\cap F_{i_n}$ is a vertex of~$P^n$.

\begin{example}[bounded flag manifold~\cite{BR0}]
\label{bfm}
A {\it bounded flag\/} in $\C^{n+1}$ is a complete flag $U=\{U_1\subset
U_2\subset\cdots\subset U_{n+1}=\C^{n+1}\}$ for which $U_k$ contains the
coordinate subspace $\C^{k-1}$ (spanned by the first $k-1$ standard basis
vectors) for $2\le k\le n$. As it is shown in~\cite[Example~2.8]{BR2}, the
$2n$-dimensional manifold $B_n$ of all bounded flags in $\C^{n+1}$ is a
quasitoric manifold over the cube $I^n$ with respect to the action induced by
$t\cdot z=(t_1z_1,\ldots,t_nz_n,z_{n+1})$ on $\C^{n+1}$ (here $t\in T^n$).
\end{example}

\begin{example}
A family of manifolds $B_{i,j}$, $0\le i\le j$, is introduced in~\cite{BR1}.
The manifold $B_{i,j}$ consists of pairs $(U,W)$, where $U$ is a bounded flag
in $\C^{i+1}$ (see Example~\ref{bfm}) and $W$ is a line in
$U_1^\bot\oplus\C^{j-i}$. So $B_{i,j}$ is a smooth $\C P^{j-1}$-bundle
over~$B_i$. It is shown in~\cite[Example~2.9]{BR2} that $B_{i,j}$ is a
quasitoric manifold over the product $I^i\times D^{j-1}$.
\end{example}

The canonical stably complex structures and omniorientations on the manifolds
$B_n$ and $B_{i,j}$ are described in~\cite[examples~4.3, 4.5]{BR2}.

\begin{remark}
The product of two quasitoric manifolds $M_1^{2n_1}$ and $M_2^{2n_2}$ over
polytopes $P^{n_1}_1$ and $P^{n_2}_2$ is a quasitoric manifold over
$P^{n_1}_1\times P^{n_2}_2$. This construction extends to omnioriented
quasitoric manifolds and is compatible with stably complex structures
(details can be found in~\cite[Proposition~4.7]{BR2}).
\end{remark}

It is shown in \cite{BR1} that the cobordism classes of $B_{i,j}$
multiplicatively generate the ring~$\O^U$. So every $2n$-dimensional
complex cobordism class may be represented by a disjoint union of products
\begin{equation}
\label{prodbij}
  B_{i_1,j_1}\times B_{i_2,j_2}\times\dots\times B_{i_k,j_k},
\end{equation}
where $\sum_{q=1}^k(i_q+j_q)-2k=n$. Each such component is a quasitoric
manifold, under the product quasitoric structure. This result is the
substance of \cite{BR1}. The stably complex structures of
products~(\ref{prodbij}) are induced by omniorientations, and are therefore
also preserved by the torus action.

\begin{example}\label{milnor}
The standard set of multiplicative generators for $\O^U$ consists of
projective spaces $\C P^i$, $i\ge0$, and {\it Milnor hypersurfaces\/}
$H_{i,j}\subset\C P^{i}\times\C P^{j}$, $1\le i\le j$. The hypersurfaces
$H_{i,j}$ are defined
by
$$
  H_{i,j}=\Bigl\{ (z_0:\cdots:z_i)\times(w_0:\cdots:w_j)\in
  \C P^{i}\times\C P^{j}:\sum_{q=0}^iz_qw_q=0\Bigr\}.
$$
However, as it was shown in~\cite{BR1}, the hypersurfaces $H_{i,j}$ are {\it
not\/} quasitoric manifolds for~$i\ge2$.
\end{example}

To give genuinely toric representatives (which are, by definition, connected)
for each cobordism class of dimension~$>2$, it remains only to replace the
disjoint union of products~(\ref{prodbij}) with their connected sum. This is
done in~\cite[\S6]{BR2} using Construction~\ref{consum} and its extension to
omnioriented quasitoric manifolds.

\begin{theorem}[Buchstaber and Ray]
\label{cobrep}
In dimensions $>2$, every complex cobordism class contains a quasitoric
manifold, necessarily connected, whose stably complex structure is
induced by an omniorientation, and is therefore compatible with the
action of the torus.
\end{theorem}

This theorem gives a solution to the quasitoric analogue of
Problem~\ref{hirz}.

\subsection{Combinatorial formulae for Hirzebruch genera of quasitoric
manifolds}\label{comb}
The constructions from the previous section open the way to evaluation of
cobordism invariants (Chern numbers, Hirzebruch genera etc.) on omnioriented
quasitoric manifolds in terms of the combinatorics of the quotient. In this
section we expose the results obtained in this direction by the second author
in \cite{Pan1},~\cite{Pan2}. Namely, using arguments similar to that from the
proof of Theorem~\ref{ds} we construct a circle action with only isolated
fixed points on any quasitoric manifold~$M^{2n}$. If $M^{2n}$ is
omnioriented, then this action preserves the stably complex structure and its
local representations near fixed points are retrieved from the characteristic
matrix~$\L$. This allows to calculate Hirzebruch's $\chi_y$-genus as the sum
of contributions corresponding to the vertices of polytope. These
contributions depend only on the ``local combinatorics" near the vertex.
After some adjustments the formula also allows to calculate the signature and
the Todd genus of~$M^{2n}$.

\begin{definition}
The {\it Hirzebruch genus\/} associated with the series
$$
  Q(x)=1+\sum q_kx^k,\quad q_k\in\Q,
$$
is the ring homomorphism $\varphi_Q:\O^U\to\Q$ that to each cobordism class
$[M^{2n}]\in\O^U_{2n}$ assigns the value given by the formula
$$
  \f_Q[M^{2n}]=\Bigl( \prod_{i=1}^nQ(x_i),\langle\M^{2n}\rangle\Bigr).
$$
Here $M^{2n}$ is a smooth manifold whose stable tangent bundle
$\tau(M^{2n})$ is a complex bundle with complete Chern class in cohomology
$$
  c(\tau)=1+c_1(\tau)+\dots+c_n(\tau)=\prod_{i=1}^n(1+x_i),
$$
and $\langle M^{2n}\rangle$ is a fundamental class in homology.
\end{definition}

The {\it $\chi_y$-genus\/} is the Hirzebruch genus associated with the series
$$
  Q(x)=\frac{x(1+ye^{-x(1+y)})}{1-e^{-x(1+y)}},
$$
where $y\in\R$ is a parameter. In particular cases $y=-1,0,1$ we obtain
the $n$-th Chern number, the {\it Todd genus\/} and the {\it
signature\/} of the manifold $M^{2n}$ correspondingly.

Provided that $M^{2n}$ is a complex manifold, the value $\chi_y(M^{2n})$ can
be calculated in terms of Euler characteristics of the {\it Dolbeault
complexes\/} on~$M^{2n}$. The general information about Hirzebruch genera can
be found in~\cite{HBY}.

In this section we assume that we are given an omnioriented quasitoric
manifold $M^{2n}$ over some~$P^n$ with characteristic matrix~$\Lambda$.
This specifies a stably complex structure on~$M^{2n}$, as described in the
previous section. The orientation of $M^{2n}$ defines the fundamental class
$\langle M^{2n}\rangle\in H_{2n}(M^{2n};\Z)$.

\begin{construction}
\label{order}
Suppose $v$ is a vertex of $P^n$ expressed as the intersection of $n$ facets:
\begin{equation}
\label{vert}
  v=F_{i_1}\cap\cdots\cap F_{i_n}.
\end{equation}
Assign to each facet $F_{i_k}$ the edge $E_k:=\bigcap_{j\ne k}F_{i_j}$ (that
is, $E_k$ contains $v$ and is opposite to~$F_{i_k}$). Let $\mb e_k$ be a
vector along $E_k$ beginning at~$v$. Then $\mb e_1,\ldots,\mb e_n$ is a basis
of~$\R^n$, which may be either positively or negatively oriented depending on
the ordering of facets in~(\ref{vert}). Throughout this section this ordering
is assumed to be so that $\mb e_1,\ldots,\mb e_n$ is a positively oriented
basis.
\end{construction}

Once we specified an ordering of facets in~(\ref{vert}), the facet vectors
$\bl_{i_1},\ldots,\bl_{i_n}$ at $v$ may in turn constitute either positively
or negatively oriented basis depending on the sign of the determinant of
$\Lambda_{(v)}=(\bl_{i_1},\ldots,\bl_{i_n})$ (see~(\ref{L})).

\begin{definition}
\label{sign}
The {\it sign\/} of a vertex $v=F_{i_1}\cap\cdots\cap F_{i_n}$ of $P^n$ is
$$
  \sigma(v):=\det\L_{(v)}.
$$
\end{definition}

The collection of signs of vertices of $P^n$ provides an important invariant
of an oriented omnioriented quasitoric manifold. Note that reversing the
orientation of $M^{2n}$ changes all signs $\sigma(v)$ to the opposite.
At the same time changing the direction of a facet vector reverses the signs
only for the vertices contained in the corresponding facet.

Let $E$ be an edge of $P^n$. The isotropy subgroup of the 2-dimensional
submanifold $\pi^{-1}(E)\subset M^{2n}$ is an $(n-1)$-dimensional subtorus,
which we denote by~$T(E)$. Then we can write
\begin{equation}
\label{eisotr}
  T(E)=\bigl\{\bigl(e^{2\pi i\f_1},\ldots,e^{2\pi i\f_n}\bigr)\in T^n\::
  \:\mu_1\f_1+\ldots+\mu_n\f_n=0\bigl\}
\end{equation}
for some integers $\mu_1,\ldots,\mu_n$. We refer to
$\bmu:=(\mu_1,\ldots,\mu_n)^t$ as the {\it edge vector\/} corresponding
to~$E$. This $\bmu$ is a primitive vector in the dual lattice~$(\Z^n)^*$; it
is determined by $E$ only up to a sign. There is no canonical way to fix
these signs simultaneously for all edges. However, the following lemma shows
that the omniorientation of $M^{2n}$ provides a canonical way to choose signs
of edge vectors ``locally" at each vertex.

\begin{lemma}
\label{lm}
  For any vertex $v\in P^n$, signs of edge vectors $\bmu_1,\ldots,\bmu_n$
  meeting at $v$ can be chosen in such way that the $n\times n$-matrix
  $\M_{(v)}:=(\bmu_1,\ldots,\bmu_n)$ satisfies the identity
  $$
    \M_{(v)}^t\cdot\L_{(v)}={\mathrm E},
  $$
  where $\mathrm E$ is the identity matrix. In other words,
  $\bmu_1,\ldots,\bmu_n$ and $\bl_{i_1},\ldots,\bl_{i_n}$ are conjugate
  bases.
\end{lemma}
\begin{proof}
At the beginning we choose signs of the edge vectors at $v$ arbitrary, and
express $v$ as in~(\ref{vert}). Then $\bmu_k$ is the edge vector
corresponding to the edge $E_k$ opposite to $F_{i_k}$, $k=1,\ldots,n$.  It
follows that $E_k\subset F_{i_l}$ and $T(F_{i_l})\subset T(E_k)$ for $l\ne
k$. Hence,
\begin{equation}
\label{mlne}
  \<\bmu_k,\bl_{i_l}\>=0, \quad l\ne k,
\end{equation}
(see~(\ref{fisotr}) and~(\ref{eisotr})). Since $\bmu_k$ is a primitive
vector, it follows from~(\ref{mlne}) that $\<\bmu_k,\bl_{i_k}\>=\pm1$.
Changing the sign of $\bmu_{i_k}$ if necessary, we obtain
$$
  \<\bmu_k,\bl_{i_k}\>=1,
$$
which together with~(\ref{mlne}) gives $\M_{(v)}^t\cdot\L_{(v)}=\mathrm E$,
as needed.
\end{proof}

In the sequel we assume that signs of edge vectors at each $v$ are chosen as
in the above lemma. It follows that the edge vectors $\bmu_1,\ldots,\bmu_n$
meeting at $v$ constitute an integer basis of $\Z^n$ and
\begin{equation}
\label{M}
  \det\M_{(v)}=\sigma(v).
\end{equation}

Suppose $M^{2n}=M_P$ is a smooth toric variety arising from a lattice simple
polytope $P$ defined by~(\ref{ptope}). Then $\bl_i=\mb l_i$, $i=1,\ldots,m$
(see Example~\ref{tcf}), while the edge vectors at $v\in P^n$ are the
primitive integer vectors $\mb e_1,\ldots,\mb e_n$ along the edges
beginning at~$v$. It follows from Construction~\ref{order} that $\sigma(v)=1$
for any~$v$. Lemma~\ref{lm} in this case expresses the fact that $\mb
e_1,\ldots,\mb e_n$ and $\mb l_{i_1},\ldots,\mb l_{i_n}$ are conjugate bases
of~$\Z^n$. Similarly, the following statement holds.

\begin{proposition}\label{acsigns}
Suppose that the omniorientation of the quasitoric manifold $M^{2n}$ arises
from a $T^m$-invariant almost complex structure (i.e. complex structure in
the tangent bundle $\tau(M^{2n})$). Then $\sigma(v)=1$ for any vertex $v\in
P^n$.
\end{proposition}
\begin{proof}
The orientation of $P^n$ is determined by the canonical orientation of the
almost complex manifold $M^{2n}$ and the orientation of the
torus~(\ref{torus}). Since the almost complex structure on $M^{2n}$ is
$T^n$-invariant, it induces almost complex structures on the $T(F_i)$-fixed
submanifolds $M^{2(n-1)}_i$. It follows that for any vertex~(\ref{vert}) the
vectors $\bl_{i_1},\ldots,\bl_{i_n}$ constitute a positively oriented basis
of~$\R^n$.
\end{proof}

Proposition \ref{acsigns} provides a necessary condition for the existence of
a $T^n$-invariant almost complex structure on $M^{2n}$ (see
Problem~\ref{acs}).

\begin{remark}
Globally Lemma~\ref{lm} provides two directions (signs) for an edge vector,
one for each of its ends. These signs are always different provided that
$M^{2n}$ is a complex manifold (e.g. a smooth toric variety), but in general
this fails to be true.
\end{remark}

Let $\bnu=(\nu_1,\ldots,\nu_n)^t\in\Z^n$ be a primitive vector such that
\begin{equation}
\label{genvect}
  \<\bmu,\bnu\>\ne0\quad\text{ for any edge vector }\bmu.
\end{equation}
The vector $\bnu$ defines the one-dimensional oriented subtorus
$$
  T_{\bnu}:=\bigl\{\bigl(e^{2\pi i\nu_1\f},\ldots,e^{2\pi
  i\nu_n\f}\bigr)\in T^n\::\:\f\in\R\bigr\}.
$$

\begin{lemma}[{\cite[Theorem~2.1]{Pan2}}]
\label{sa}
  For any $\bnu$ satisfying {\rm(\ref{genvect})} the circle $T_{\bnu}$ acts on
  $M^{2n}$ with only isolated fixed points (corresponding to the vertices
  of~$P^n$). For each vertex $v=F_{i_1}\cap\cdots\cap F_{i_n}$ the action of
  $T_{\bnu}$ induces a representation of $S^1$ in the tangent space
  $T_vM^{2n}$ with weights $\<\bmu_1,\bnu\>,\ldots,\<\bmu_n,\bnu\>$.
\end{lemma}

\begin{remark}
  If $M^{2n}=M_P$ is a smooth toric variety, then the genericity
  condition~(\ref{genvect}) is equivalent to that from the proof of
  Theorem~\ref{ds}.
\end{remark}

\begin{definition}
\label{ind}
Suppose we are given a primitive vector $\bnu$ satisfying~(\ref{genvect}).
Define the {\it index\/} of a vertex $v\in P^n$ as the number of negative
weights of the $S^1$-representation in $T_vM^{2n}$ from Lemma~\ref{sa}. That
is, if $v=F_{i_1}\cap\cdots\cap F_{i_n}$, then
$$
  \ind_{\bnu}(v)=\{\# k:\<\bmu_k,\bnu\><0\}.
$$
\end{definition}

\begin{remark}
The index of a vertex $v$ can be also defined in terms of the facet vectors
at~$v$. Indeed, Lemma~\ref{lm} shows that if
$v=F_{i_1}\cap\cdots\cap F_{i_n}$, then
$$
  \bnu=\<\bmu_1,\bnu\>\bl_{i_1}+\dots+\<\bmu_n,\bnu\>\bl_{i_n}.
$$
Hence, $\ind_{\bnu}(v)$ equals the number of negative coefficients in
the representation of $\bnu$ as a linear combination of basis vectors
$\bl_{i_1},\ldots,\bl_{i_n}$.
\end{remark}

\begin{theorem}[{\cite[Theorem~6]{Pan1}, \cite[Theorem~3.1]{Pan2}}]
\label{chi}
  For any $\bnu$ satisfying~{\rm(\ref{genvect})}, the $\chi_y$-genus of
  $M^{2n}$ can be calculated as
  $$
    \chi_y(M^{2n})=\sum_{v\in P^n}(-y)^{\ind_\nu(v)}\sigma(v).
  $$
\end{theorem}
This theorem is proved by applying the Atiyah--Hirzebruch
formula~\cite{AH} to the circle action defined in Lemma~\ref{sa}.

The value of the $\chi_y$-genus $\chi_y(M^{2n})$ at $y=-1$ equals the $n$-th
Chern number $c_n(\xi)\langle M^{2n}\rangle$ for any $2n$-dimensional stably
complex manifold $[M^{2n},\xi]$. In quasitoric case Theorem~\ref{chi} gives
\begin{equation}
\label{cn}
  c_n[M^{2n}]=\sum_{v\in P^n}\sigma(v).
\end{equation}
If $M^{2n}$ is a complex manifold (e.g. a smooth toric variety), then
$\sigma(v)=1$ for all vertices $v\in P^n$ and $c_n[M^{2n}]$ equals the Euler
characteristic $e(M^{2n})$. Hence, for complex $M^{2n}$ the Euler
characteristic equals the number of vertices of $P^n$ (which is well known
for toric varieties). For general quasitoric $M^{2n}$ the Euler
characteristic is also equal to the number of vertices of $P^n$ (since the
Euler characteristic of any $S^1$-manifold equals the sum of Euler
characteristics of fixed submanifolds), however the latter number may differ
from $c_n[M^{2n}]$ (see Example~\ref{cpquas} below).

The value of the $\chi_y$-genus at $y=1$ is the {\it signature\/} (or the
{\it $L$-genus\/}). Theorem~\ref{chi} gives in this case
\begin{corollary}
\label{signat}
 The signature of an omnioriented quasitoric manifold $M^{2n}$
 can be calculated as
 $$
   \sign(M^{2n})=\sum_{v\in P^n}(-1)^{\ind_\nu(v)}\sigma(v).
 $$
\end{corollary}

Being an invariant of an {\it oriented\/} cobordism class, the signature of
$M^{2n}$ does not depend on a stably complex structure (i.e. on an
omniorientation) and is determined only by an orientation of $M^{2n}$
(or~$P^n$). The following modification of Corollary~\ref{signat} provides a
formula for $\sign(M^{2n})$ that does not depend on an omniorientation.

\begin{corollary}[{\cite[Corollary~3.3]{Pan2}}]
\label{signor}
  The signature of an oriented quasitoric manifold $M^{2n}$ can be calculated
  as
  $$
  \sign(M^{2n})=\sum_{v\in P^n}\det(\widetilde{\bmu}_1,\ldots,
  \widetilde{\bmu}_n),
  $$
  where $\widetilde{\bmu}_k$, $k=1,\ldots,n$, are the edge vectors at $v$
  oriented in such way that ${\<\widetilde{\bmu}_k,\bnu\>>0}$.
\end{corollary}

If $M^{2n}=M_P$ is a smooth toric variety, then $\sigma(v)=1$ for any
$v\in P^n$, and Corollary~\ref{signat} gives
$$
  \sign(M_P)=\sum_{v\in P^n}(-1)^{\ind_\nu(v)}.
$$
Since in this case the index $\ind_\nu(v)$ is the same as the index from the
proof of Theorem~\ref{ds}, we obtain
\begin{equation}
\label{toricsign}
  \sign(M_P)=\sum_{k=1}^n(-1)^kh_k(P).
\end{equation}
Note that if $n$ is odd, then the right hand side of the above formula
vanishes due to the Dehn--Sommerville equations. The
formula~(\ref{toricsign}) appears in a more general context in recent
work of Leung and Reiner~\cite{LR}. The quantity in the right hand side
of~(\ref{toricsign}) arises in the following well known combinatorial
conjecture.
\begin{problem}[Charney--Davis conjecture]
  Let $K$ be a $(2q-1)$-dimen\-si\-o\-nal Gorenstein* flag complex with
  $h$-vector $(h_0,h_1,\ldots,h_{2q})$. Is it true that
  $$
    (-1)^q(h_0-h_1+\dots+h_{2q})\ge0?
  $$
\end{problem}
\noindent This conjecture was made in~\cite[Conjecture~D]{CD} for flag
simplicial homology spheres. Stanley~\cite[Problem~4]{St5} proposed to extend
it to Gorenstein* complexes. The Charney--Davis conjecture is closely
connected with the following differential-geometrical conjecture.

\begin{problem}[Hopf conjecture]
  Let $M^{2q}$ be a Riemannian manifold of non-positive sectional curvature.
  Is it true that the Euler characteristic $\chi(M^{2n})$ satisfies the
  inequality
  $$
    (-1)^q\chi(M^{2q})\ge0?
  $$
\end{problem}
Both above conjectures are known to be true for $q=1,2$ and for some special
cases. More details can be found in~\cite{CD}. For more relations of the
two problems with the signature of a toric variety see~\cite{LR}.

Now we turn again to the $\chi_y$-genus of an omnioriented quasitoric
manifold. The next important particular case is the {\it Todd genus\/}
corresponding to~$y=0$. In this case the summands in the formula from
Theorem~\ref{chi} are not defined for the vertices of index 0, so it
requires some additional analysis.

\begin{theorem}[{\cite[Theorem~7]{Pan1}, \cite[Theorem~3.4]{Pan2}}]
\label{todd}
  The Todd genus of an omnioriented quasitoric manifold
  can be calculated as
  $$
    \td(M^{2n})=\sum_{v\in P^n:\ind_\nu(v)=0}\sigma(v)
  $$
  (the sum is taken over all vertices of index~{\rm0}).
\end{theorem}

In the case of smooth toric variety there is only one vertex of index~0. This
is the ``bottom" vertex of~$P^n$, which has all incident edges pointing out
(in the notations used in the proof of Theorem~\ref{ds}). Since $\sigma(v)=1$
for every $v\in P^n$, Theorem~\ref{todd} gives $\td(M_P)=1$, which is well
known (see e.g.~\cite[\S5.3]{Fu}; for algebraic varieties the Todd genus
equals the {\it arithmetic genus\/}).

If $M^{2n}$ is an almost complex manifold, then $\td(M^{2n})>0$ by
Proposition~\ref{acsigns} and Theorem~\ref{todd}.

\begin{example}
\label{cpnoo}
The projective space $\C P^2$, regarded as a toric variety, arises from the
standard lattice 2-simplex $\D^2$ (with vertices $v_1=(0,0)$, $v_2=(1,0)$,
$v_3=(0,1)$). The orientation is standard (determined by the complex
structure). The omniorientation is determined by the facet vectors
$\bl_1,\bl_2,\bl_3$, which in this case are the primitive normal vectors
pointing inside the polytope. The edge vectors are the primitive vectors
along edges pointing out of the vertex. This can be seen on Figure~5. The
corresponding stably complex structure is the standard one, that is,
determined by the isomorphism of bundles $\tau(\C
P^2)\oplus\R^2\cong\bar\eta\oplus\bar\eta\oplus\bar\eta$, where $\eta$ is the
canonical Hopf line bundle. Let us calculate the Todd genus and the signature
from Corollary~\ref{signat} and Theorem~\ref{todd}. We have
$\sigma(v_1)=\sigma(v_2)=\sigma(v_3)=1$. Take $\bnu=(1,2)$, then
$\ind_\nu(v_1)=0$, $\ind_\nu(v_2)=1$, $\ind_\nu(v_3)=2$ (recall that the
index is the number of negative scalar products of edge vectors with~$\bnu$).
Thus, $\sign(\C P^2)= \sign[\C P^2,\bar\eta\oplus\bar\eta\oplus\bar\eta]=1$,
$\td(\C P^2)=\td[\C P^2,\bar\eta\oplus\bar\eta\oplus\bar\eta]=1$.
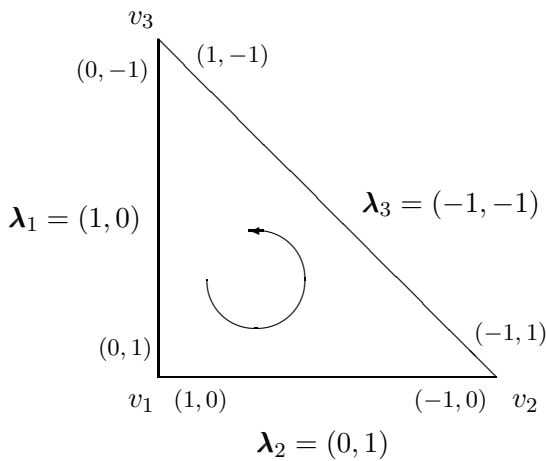
\begin{figure}
\begin{center}
\begin{picture}(100,60)
  \put(20,10){\line(0,1){45}}
  \put(20,10){\line(1,0){45}}
  \put(20,55){\line(1,-1){45}}
  \put(33,23){\oval(13,13)[b]}
  \put(33,23){\oval(13,13)[tr]}
  \put(34,29.5){\vector(-1,0){2}}
  \put(16,6){{\large $v_1$}}
  \put(22,6){\small $(1,0)$}
  \put(54,6){\small $(-1,0)$}
  \put(67,6){{\large $v_2$}}
  \put(33,0){{\large $\bl_2=(0,1)$}}
  \put(12,13){\small $(0,1)$}
  \put(0,30){{\large $\bl_1=(1,0)$}}
  \put(9,50){\small $(0,-1)$}
  \put(16,57){{\large $v_3$}}
  \put(25,52){\small $(1,-1)$}
  \put(62,15){\small $(-1,1)$}
  \put(47,32){{\large $\bl_3=(-1,-1)$}}
\end{picture}%
\caption{$\tau(\C P^2)\oplus\C\simeq\bar\eta\oplus\bar\eta\oplus\bar\eta$}
\end{center}
\end{figure}
\end{example}

\begin{example}\label{cpquas}
Now consider $\C P^2$ with the omniorientation determined by the three facet
vectors $\bl_1,\bl_2,\bl_3$ shown on Figure~6. This omniorientation differs
from the previous example by the sign of~$\bl_3$. The corresponding stably
complex structure is determined by the isomorphism $\tau(\C
P^2)\oplus\R^2\cong\bar\eta\oplus\bar\eta\oplus\eta$. Using~(\ref{M}) we
calculate
$$
  \sigma(v_1)=\begin{vmatrix} 1&0\\0&1 \end{vmatrix}=1,\quad
  \sigma(v_2)=\begin{vmatrix} -1&1\\1&0 \end{vmatrix}=-1,\quad
  \sigma(v_3)=\begin{vmatrix} 0&1\\1&-1 \end{vmatrix}=-1.
$$
For $\bnu=(1,2)$ we find $\ind_\nu(v_1)=0$, $\ind_\nu(v_2)=0$,
$\ind_\nu(v_3)=1$.  Thus, $\sign[\C P^2,\bar\eta\oplus\bar\eta\oplus\eta]=1$,
$\td[\C P^2,\bar\eta\oplus\bar\eta\oplus\eta]=0$. Note that in this case
formula~(\ref{cn}) gives $c_n[\C P^2,\bar\eta\oplus\bar\eta\oplus\eta]=
\sigma(v_1)+\sigma(v_2)+\sigma(v_3)=-1$, while the Euler number of $\C P^2$
is $c_n[\C P^2,\bar\eta\oplus\bar\eta\oplus\bar\eta]=3$.
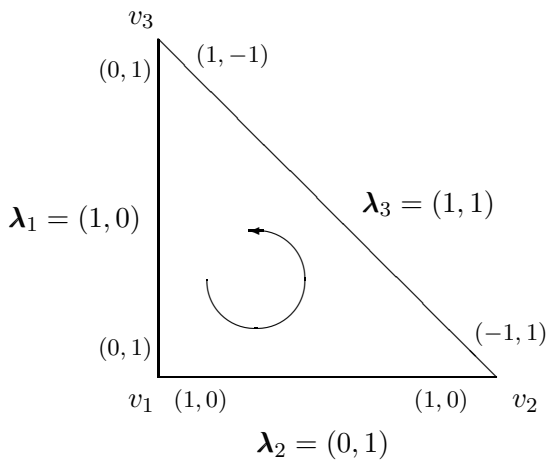
\begin{figure}
\begin{center}
\begin{picture}(100,60)
  \put(20,10){\line(0,1){45}}
  \put(20,10){\line(1,0){45}}
  \put(20,55){\line(1,-1){45}}
  \put(33,23){\oval(13,13)[b]}
  \put(33,23){\oval(13,13)[tr]}
  \put(34,29.5){\vector(-1,0){2}}
  \put(16,6){{\large $v_1$}}
  \put(22,6){\small $(1,0)$}
  \put(54,6){\small $(1,0)$}
  \put(67,6){{\large $v_2$}}
  \put(33,0){{\large $\bl_2=(0,1)$}}
  \put(12,13){\small $(0,1)$}
  \put(0,30){{\large $\bl_1=(1,0)$}}
  \put(12,50){\small $(0,1)$}
  \put(16,57){{\large $v_3$}}
  \put(25,52){\small $(1,-1)$}
  \put(62,15){\small $(-1,1)$}
  \put(47,32){{\large $\bl_3=(1,1)$}}
\end{picture}%
\caption{$\tau(\C P^2)\oplus\C\simeq\bar\eta\oplus\bar\eta\oplus\eta$}
\end{center}
\end{figure}
\end{example}

$T^n$-equivariant stably complex and almost complex manifolds were considered
in works of Hattori~\cite{Ha} and Masuda~\cite{Mas} as a separate
generalisation (called the {\it unitary toric manifolds\/}) of toric
varieties. Instead of Davis and Januszkiewicz's characteristic maps, Masuda
in~\cite{Mas} uses the notion of {\it multi-fan\/} to describe the
combinatorial structure of the orbit space. The multi-fan is a collection of
cones which may overlap unlike a usual fan. The Todd genus of a unitary toric
manifold was calculated in~\cite{Mas} via the degree of the overlap of cones
in the multi-fan. This result is equivalent to our Theorem~\ref{todd} in the
case of quasitoric manifolds.

\subsection{The classification problem for quasitoric manifolds
over a given simple polytope}
\label{prob}
More precisely, there are two classification problems: the equivariant
(i.e. up to an equivariant diffeomorphism) and the topological (i.e. up
to a diffeomorphism). Due to Proposition~\ref{equivar}, the equivariant
classification problem reduces to the description of all characteristic maps
for given simple polytope~$P^n$. The topological classification problem
usually requires an additional analysis.

Let $M^{2n}$ be a quasitoric manifold over~$P^n$ with characteristic
map~$\ell$. We suppose here that the first $n$ facets $F_1,\ldots,F_n$ share
a common vertex.

\begin{lemma}\label{ucf}
Up to $\psi$-equivalence, we may assume that $\ell(F_i)$ is the $i$-th
coordinate subtorus $T_i\subset T^n$, $i=1,\ldots,n$.
\end{lemma}
\begin{proof}
Since the one-dimensional subtori $\ell(F_i)$, $i=1,\ldots,n$,
generate~$T^n$, we may define $\psi$ as any automorphism of $T^n$ that
maps $\ell(F_i)$ to~$T_i$, $i=1,\ldots,n$.
\end{proof}
It follows that $M^{2n}$ admits such an omniorientation that the
corresponding characteristic $n\times m$-matrix $\L$ has the form $(\mathrm
E\,|\,{*})$, where $\mathrm E$ is the identity matrix and ${*}$ is a certain
integer $n\times(m-n)$-matrix.

In the simplest case $P^n=\D^n$ the equivariant (and topological)
classification of quasitoric manifolds reduces to the following easy result.

\begin{proposition}\label{clsim}
Any quasitoric manifold over the simplex $\D^n$ is equivariantly
diffeomorphic to $\C P^n$ (regarded as a toric variety, see
examples~{\rm\ref{cpn}} and~{\rm\ref{tcf}}).
\end{proposition}
\begin{proof}
The characteristic map for $\C P^n$ has the form
$$
  \ell_{\C P^n}(F_i)=T_i,\quad i=1,\ldots,n,\quad
  \ell_{\C P^n}(F_{n+1})=S_d,
$$
where $S^1_d:=\{(e^{2\pi i\f},\ldots,e^{2\pi i\f})\in T^n\}$, $\f\in\R$, is
the diagonal subgroup in~$T^n$. Let $M^{2n}$ be a quasitoric manifold over
$\D^n$ with characteristic map~$\ell_M$. We may assume that $\ell_M(F_i)=T_i$,
$i=1,\ldots,n$, by Lemma~\ref{ucf}. Then it easily follows from~(\ref{L})
that
$$
  \ell_M(F_{n+1})=\bigl\{\bigl(e^{2\pi i\e_1\f},\ldots,e^{2\pi
  i\e_n\f}\bigr)\in T^n\bigr\}, \quad\f\in\R,
$$
where $\e_i=\pm1$, $i=1,\ldots,n$. Now define the automorphism
$\psi:T^n\to T^n$ by
$$
  \psi\bigl(e^{2\pi i\f_1},\ldots,e^{2\pi i\f_n}\bigr)=
  \bigl(e^{2\pi i\e_1\f_1},\ldots,e^{2\pi i\e_n\f_n}\bigr).
$$
It can be readily seen that $\psi\cdot\ell_M=\ell_{\C P^n}$, which together
with Proposition~\ref{equivar} completes the proof.
\end{proof}

Note that the equivariant diffeomorphism provided by Proposition~\ref{clsim}
not necessarily preserves the orientation of the quasitoric
manifold~$M^{2n}$.

The problem of equivariant and topological classification also admits a
complete solution in the case $n=2$ (i.e. for quasitoric manifolds over
polygons).

\begin{example}
The {\it Hirzebruch surface\/} is the 2-dimensional complex manifold $H_p=\C
P(\zeta_p\oplus\C)$, where $\zeta_p$ is the complex line bundle over $\C P^1$
with first Chern class $p$, $\C$ is the trivial complex line bundle and
$\C P(\cdot)$ denotes the projectivisation of the complex bundle.  Hence,
there is a bundle $H_p\to\C P^1$ with fibre~$\C P^1$. The surface $H_p$ is
homeomorphic to $S^2\times S^2$ for even $p$ and to $\C P^2\cs\overline{\C
P}{}^2$ for odd $p$, where $\overline{\C P}{}^2$ denotes the space $\C P^2$
with reversed orientation.  The Hirzebruch surfaces are non-singular
projective toric varieties (see~\cite[p.~8]{Fu}). The orbit space for $H_p$
(regarded as a quasitoric manifold) is a combinatorial square; the
corresponding characteristic maps can be described using Example~\ref{tcf}
(see also~\cite[Example~1.19]{DJ}).
\end{example}

\begin{theorem}[{\cite[p.~553]{OR}}]
A quasitoric manifold of dimension {\rm4} is equivariantly diffeomorphic to
the equivariant connected sum of several copies of $\C P^2$ and Hirzebruch
surfaces~$H_p$.
\end{theorem}

\begin{corollary}\label{4dtopclas}
A quasitoric manifold of dimension {\rm4} is diffeomorphic to the connected
sum of several copies of $\C P^2$, $\overline{\C P}{}^2$ and $S^2\times S^2$.
\end{corollary}

The classification problem for quasitoric manifolds over a given simple
polytope can be considered as a generalisation of the corresponding problem
for non-singular toric varieties. In~\cite{Od} to every toric variety over a
simple 3-polytope $P^3$ was assigned two integer {\it weights\/} on every
edge of the dual simplicial complex~$K_P$. Using the special ``monodromy
conditions" for weights, the complete classification of toric varieties over
simple 3-polytopes with~$\le8$ facets was obtained in~\cite{Od}. A similar
construction was used in~\cite{Kls} to obtain the classification of toric
varieties over $P^n$ with $m=n+2$ facets (note that any such simple polytope
is the product of two simplices).

In \cite{Do} the construction of weights from~\cite{Od} was generalised to
the case of quasitoric manifolds. This allowed to obtain a
criterion~\cite[Theorem~3]{Do} for the existence of a quasitoric manifold
with given weight set and signs of vertices $\sigma(v)$ (see
Definition~\ref{sign}; note that our signs of vertices correspond to the
two-paint colouring used in~\cite{Do}). As an application, the complete
classification of characteristic maps for the cube $I^3$ is obtained
in~\cite{Do}, along with a series of results on the classification of
quasitoric manifolds over the product of arbitrary number of simplices.

\medskip

\section{Moment-angle complexes}
\subsection{Moment-angle manifolds $\zp$ defined by simple polytopes}
\label{mom1}
For any combinatorial simple polytope $P^n$ with $m$ facets Davis and
Januszkiewicz introduced in~\cite{DJ} a space $\zp$ with an action
of the torus $T^m$ and the orbit space~$P^n$. This space is universal for all
quasitoric manifolds over $P^n$ in the sense that for every quasitoric
manifold $\pi:M^{2n}\to P^n$ there is a principal $T^{m-n}$-bundle $\zp\to
M^{2n}$ whose composite map with $\pi$ is the orbit map for~$\zp$.  The space
$\zp$ and some its generalisations turn to be very important and effective
tool for studying different combinatorial objects such as Stanley--Reisner
rings, subspace arrangements, cubical complexes etc. In this section we
reproduce the original definition of $\zp$ and adjust it in the way
convenient for subsequent generalisations.

Let $\F=\{F_1,\ldots,F_m\}$ be the set of facets of~$P^n$. For any facet
$F_i\in\F$ denote by $T_{F_i}$ the one-dimensional coordinate subgroup of
$T^\F\cong T^m$ corresponding to~$F_i$. Then assign to every face $G$ the
coordinate subtorus
$$
  T_G=\bigoplus_{F_i\supset G}T_{F_i}\subset T^\F.
$$
Note that $\dim T_G=\mathop{\rm codim}G$. Recall that for any point $q\in
P^n$ we denoted by $G(q)$ the unique face containing $q$ in the relative
interior.

\begin{definition}
\label{zp}
For any combinatorial simple polytope $P^n$ introduce the identification
space
$$
  \zp=(T^\F\times P^n)/{\sim},
$$
where $(t_1,q)\sim(t_2,q)$ if and only if $t_1t_2^{-1}\in T_{G(q)}$.
\end{definition}
The free action of $T^m$ on $T^\F\times P^n$ obviously descends to an action
on $\zp$, with quotient~$P^n$. Let $\rho:\zp\to P^n$ be the orbit map. The
action of $T^m$ on $\zp$ is free over the interior of~$P^n$, while each
vertex $v\in P^n$ represents the orbit $\rho^{-1}(v)$ with maximal isotropy
subgroup of dimension~$n$.

\begin{lemma}
\label{zpman}
The space $\zp$ is a smooth manifold of dimension~$m+n$.
\end{lemma}

In this paper we provide several different proofs of this lemma,
each of which arises from an equivalent definition of~$\zp$. To give our
first proof we need the following simple topological fact.

\begin{proposition}
\label{toremb}
The torus $T^k$ admits a smooth embedding into~$\R^{k+1}$.
\end{proposition}
\begin{proof}
In order to construct a required embedding we provide a smooth function
$g_{k+1}(x_1,\ldots,x_{k+1})$ such that the equation $g_{k+1}=0$ defines a
hypersurface diffeomorphic to~$T^k$. For $k=1$ we take
$g_2(x_1,x_2)=x_1^2+(x_2-2)^2-1$; then we proceed by induction on~$k$.
Suppose we have a function $g_i(x_1,\ldots,x_i)$ such that $\{g_i=0\}$
defines a smooth embedding $T^{i-1}\hookrightarrow\R^i$ and $x_i>0$ for any
$(x_1,\ldots,x_i)$ satisfying $g_i=0$. Then set
$$
  g_{i+1}(x_1,\ldots,x_i,x_{i+1}):=
  g_i\bigl(x_1,\ldots,x_{i-1},\sqrt{x_i^2+x_{i+1}^2}\bigr).
$$
The hypersurface $\{g_{i+1}=0\}\subset\R^{i+1}$ is easily seen to be
diffeomorphic to~$T^i$.
\end{proof}

\begin{proof}[Proof of Lemma~{\rm\ref{zpman}}]
Construction~\ref{corners} provides the atlas $\{U_v\}$ for $P^n$ as a
manifold with corners. The set $U_v$ is based on the vertex $v$ and is
diffeomorphic to~$\R^n_+$. Then $\rho^{-1}(U_v)\cong T^{m-n}\times\R^{2n}$.
We claim that $T^{m-n}\times\R^{2n}$ can be realised as an open set
in~$\R^{m+n}$, thus providing a chart for~$\zp$. To see this we embed
$T^{m-n}$ into $\R^{m-n+1}$ as a closed hypersurface~$H$
(Proposition~\ref{toremb}). Since the normal bundle is trivial, the small
neighbourhood of $H\subset\R^{m-n+1}$ is homeomorphic to $T^{m-n}\times\R$.
Taking the cartesian product with $\R^{2n-1}$ we obtain an open set in
$\R^{m+n}$ homeomorphic to $T^{m-n}\times\R^{2n}$.
\end{proof}

The following statement follows easily from the definition of~$\zp$.

\begin{proposition}
\label{zpprod}
If $P=P_1\times P_2$ for some simple polytopes $P_1$,~$P_2$, then
$\zp=\mathcal Z_{P_1}\times\mathcal Z_{P_2}$. If $G\subset P$ is a face, then
$\mathcal Z_G$ is a submanifold of~$\zp$.
\end{proposition}

Suppose now that we are given a characteristic map $\ell$ on $P^n$ and
$M^{2n}(\ell)$ is the derived quasitoric manifold (Construction~\ref{der}).
Choosing an omniorientation in any way we obtain a directed characteristic
map $\l:T^\F\to T^n$. Denote its kernel by $H(\ell)$ (it depends only
on~$\ell$); then $H(\ell)$ is an $(m-n)$-dimensional subtorus of~$T^\F$.

\begin{proposition}
\label{zpbun}
The subtorus $H(\ell)$ acts freely on~$\zp$, thus defining a principal
$T^{m-n}$-bundle $\zp\to M^{2n}(\ell)$.
\end{proposition}
\begin{proof}
It follows from~(\ref{L}) that $H(\ell)$ meets every isotropy subgroup
only at the unit. This implies that the action of $H(\ell)$ on $\zp$ is
free. By definitions of $\zp$ and $M^{2n}(\ell)$, the projection
$\l\times\id:T^\F\times P^n\to T^n\times P^n$ descends to the projection
$$
  (T^\F\times P^n)/{\sim}\longrightarrow (T^n\times P^n)/{\sim},
$$
which displays $\zp$ as a principal $T^{m-n}$-bundle over~$M^{2n}(\ell)$.
\end{proof}

To simplify notations we would write $T^m$, $\C^m$ etc. instead of
$T^\F$, $\C^\F$ etc.

Define the unit poly-disk $(D^2)^m$ in the complex space as
$$
  (D^2)^m=\bigl\{ (z_1,\ldots,z_m)\in\C^m:\: |z_i|\le1,\quad i=1,\ldots,m
  \bigr\}.
$$
Then $(D^2)^m$ is stable under the standard action of $T^m$ on~$\C^m$, and
the quotient is the standard cube $I^m\subset\R_+^m$.

\begin{lemma}\label{ie}
  The cubical embedding $i_P:P^n\to I^m$ from Construction~{\rm\ref{cubpol}}
  is covered by an equivariant embedding $i_e:\zp\to(D^2)^m$.
\end{lemma}
\begin{proof}
Recall that the cubical decomposition of $P^n$ consists of the cubes $C^n_v$
based on the vertices $v\in P^n$. Note that $C^n_v$ is contained in the open
set $U_v\subset P^n$ (see Construction~\ref{corners}). The inclusion
$C^n_v\subset U_v$ is covered by an equivariant inclusion $B_v\subset\C^m$,
where $B_v=\rho^{-1}(C^n_v)$ is a closed subset homeomorphic to
$(D^2)^m\times T^{m-n}$. Since $\zp=\bigcup_{v\in P^n}B_v$ and $B_v$ is
stable under the $T^m$-action, the resulting embedding $\zp\to(D^2)^m$ is
equivariant.
\end{proof}

It follows from the proof that the manifold $\zp$ is represented as a union
of $f_{n-1}(P)$ closed $T^m$-invariant charts~$B_v$. In section~\ref{cell}
we provide two different ways to construct a cellular decomposition of
each~$B_v$, thus describing $\zp$ as a cellular complex. For now, we mention
that if $v=F_{i_1}\cap\dots\cap F_{i_n}$, then
$$
  i_e(B_v)=(D^2)^n_{i_1,\ldots,i_n}\times
  T^{m-n}_{[m]\setminus\{i_1,\ldots,i_n\}}\subset(D^2)^m,
$$
or, more precisely,
$$
  i_e(B_v)=\bigl\{(z_1,\ldots,z_m)\in
  (D^2)^m\: : \:|z_i|=1\text{ for }i\notin\{i_1,\ldots,i_n\}\bigr\}.
$$
Recalling that the vertices of $P^n$ correspond to the maximal simplices of
the polytopal sphere $K_P$ (boundary of the polar polytope~$P^*$), we can
write
\begin{equation}
\label{mazp}
  i_e(\zp)=\bigcup_{I\in K_P}(D^2)_I\times T_{[m]\setminus I}\subset(D^2)^m.
\end{equation}
This can be regarded as an alternative definition of~$\zp$. Introducing the
polar coordinates in $(D^2)^m$ we see that $i_e(B_v)$ is parametrised by $n$
radial (or moment) and $m$ angle coordinates. That is why we refer to
$\zp$ as the {\it moment-angle manifold defined by\/}~$P^n$.

\begin{example}\label{zpsphere}
Let $P^n=\D^n$ (the $n$-simplex). Then $\zp$ is homeomorphic to the
$(2n+1)$-sphere~$S^{2n+1}$. The cubical complex $\mathcal C(\D^n)$ (see
Construction~\ref{cubpol}) consists of $(n+1)$ cubes~$C_v^n$. Each subset
$B_v=\rho^{-1}(C_v^n)$ is homeomorphic to $(D^2)^n\times S^1$. In particular,
for $n=1$ we obtain the well-known representation of the 3-sphere $S^3$ as a
union of two solid tori $D^2\times S^1$.
\end{example}

Another way to construct an equivariant embedding of $\zp$ into $\C^m$ can be
derived from Construction~\ref{dist}.

\begin{construction}
Formula (\ref{af}) defines the (affine) embedding
$A_P:P^n\hookrightarrow\R^m_+$. This embedding is covered by an equivariant
embedding $\zp\hookrightarrow\C^m$. A choice of matrix $W$ in
Construction~\ref{dist} defines a basis in the $(m-n)$-dimensional subspace
orthogonal to the $n$-plane that contain $A_P(P^n)$ (see~(\ref{af})). The
following statement follows.

\begin{corollary}[{see also~\cite[\S3]{BR2}}] The embedding
$\zp\hookrightarrow\C^m$ has the trivial normal bundle.
\end{corollary}
\end{construction}

\subsection{General moment-angle complexes}
\label{mom2}
In this section we extend the construction of $\zp$ to the case of general
simplicial complex~$K$. The resulting space is not a manifold for
arbitrary~$K$, but is so when $K$ is a simplicial sphere.

Here we denote by $\rho$ the canonical projection $(D^2)^m\to I^m$, as well
as any of its restriction to a closed $T^m$-stable subset of~$(D^2)^m$.
For each face $C_{I\subset J}$ of $I^m$ (see~(\ref{ijface})) define
\begin{multline}
\label{bij}
  B_{I\subset J}:=\rho^{-1}(C_{I\subset J})\\ =\{(z_1,\ldots,z_m)\in
  (D^2)^m\: : \: z_i=0\text{ for }i\in I,\; |z_i|=1\text{ for }i\notin J\}.
\end{multline}
It follows that if $\#I=i$, $\#J=j$ then $B_{I\subset
J}\cong(D^2)^{j-i}\times T^{m-j}$, where the disk factors
$D^2\subset(D^2)^{j-i}$ are parametrised by $J\setminus I$, while
the circle factors $S^1\subset T^{m-j}$ are parametrised by $[m]\setminus J$.

\begin{definition}
\label{ma}
  Let $\mathcal C$ be a cubical subcomplex of~$I^m$. The {\it moment-angle
  complex\/} $\ma(\mathcal C)$ corresponding to $\mathcal C$ is the
  $T^m$-invariant decomposition of $\rho^{-1}(\mathcal C)$ into the
  ``moment-angle" blocks $B_{I\subset J}$~(\ref{bij}) corresponding
  to the faces $C_{I\subset J}$ of~$\mathcal C$. Hence, $\ma(\mathcal C)$ is
  defined from the commutative diagram
  $$
  \begin{CD}
    \ma(\mathcal C) @>>> (D^2)^m\\
    @VVV @VV\rho V\\
    \mathcal C @>>> I^m
  \end{CD}.
  $$
\end{definition}
\noindent The torus $T^m$ acts on $\ma(\mathcal C)$ with orbit space~$\mathcal
C$.

In section \ref{cubi} two canonical cubical subcomplexes of~$I^m$, namely
$\cub(K)$~(\ref{fcubk}) and $\cc(K)$~(\ref{fcck}), were associated to every
simplicial complex $K^{n-1}$ on $m$ vertices. We denote the corresponding
moment-angle complexes by $\wk$ and $\zk$ respectively. Thus, we have
\begin{align}
\label{zkwk}
  \begin{CD}
    \wk @>>> (D^2)^m\\
    @V\rho VV @VV\rho V\\
    \cub(K) @>>> I^m
  \end{CD} &&
  \text{and} &&
  \begin{CD}
    \zk @>>> (D^2)^m\\
    @V\rho VV @VV\rho V\\
    \cc(K) @>>> I^m
  \end{CD}, &
\end{align}
where the horizontal arrows are embeddings, while the vertical ones are orbit
maps for $T^m$-actions. Note that $\dim\zk=m+n$ and $\dim\wk=m+n-1$.

\begin{remark}
Suppose that $K=K_P$ for some simple polytope~$P$. Then it follows
from~(\ref{mazp}) that $\zk$ is identified with $\zp$ (or, more precisely,
with~$i_e(\zp)$).
\end{remark}

\begin{lemma}
\label{maman}
If $K$ is a simplicial $(n-1)$-sphere, then $\zk$ is an $(m+n)$-dimensional
(closed) manifold.
\end{lemma}
\begin{proof}
In this proof we identify the polyhedrons $|K|$ and $|\cone(K)|$ with their
images $\cub(K)\subset I^m$ and $\cc(K)\subset I^m$ under the map
$|\cone(K)|\to I^m$, see Theorem~\ref{cubkcck}. For each vertex $\{i\}\in K$
denote by $F_i$ the union of $(n-1)$-cubes of $\cub(K)$ that contain~$\{i\}$.
Alternatively, $F_i$ is $|\star_{\bs(K)}\{i\}|$. These $F_1,\ldots,F_m$ are
analogues of facets of a simple polytope. Moreover, if $K=K_P$ for some~$P$,
then $F_i$ is the image of a facet of $P$ under the map $i_P:\mathcal C(P)\to
I^m$ (Construction~\ref{cubpol}). As in the case of simple polytopes, we
define ``faces" of $\cc(K)$ as non-empty intersections of ``facets"
$F_1,\ldots,F_m$. Then the ``vertices" (i.e. non-empty intersections of $n$
``facets") are the barycentres of $(n-1)$-simplices of~$|K|$. For every such
barycentre $b$ denote by $U_b$ the open subset of $\cc(K)$ obtained by
deleting all ``faces" not containing~$b$.  Then $U_b$ is identified
with~$\R^n_+$, while $\rho^{-1}(U_b)$ is homeomorphic to
$T^{m-n}\times\R^{2n}$. This defines a structure of manifold with corners on
the $n$-ball $\cc(K)=|\cone(K)|$, with atlas~$\{U_b\}$. At the same time we
see that $\zk=\rho^{-1}(\cc(K))$ is a manifold, with atlas
$\{\rho^{-1}(U_b)\}$.
\end{proof}

\begin{problem}
\label{zkmanprob}
Characterise simplicial complexes $K$ for which $\zk$ is a manifold.
\end{problem}

As we will see below (Theorem~\ref{cohom1}), for homological reasons, if
$\zk$ is a manifold, then $K$ is a Gorenstein* complex. So, the answer to
the above problem is somewhere between ``simplicial spheres" and
``Gorenstein* complexes".

\subsection{Cellular structures on moment-angle complexes}
\label{cell}
Here we consider two cellular decompositions of~$(D^2)^m$, which provide
cellular decompositions for moment-angle complexes. The first one has $5^m$
cells and define a cellular complex structure (with 5 types of cells) for any
moment-angle complex $\ma(\mathcal C)\subset(D^2)^m$. The second cellular
decomposition of $(D^2)^m$ has only $3^m$ cells, but is appropriate only for
defining a cellular structure (with 3 types of cells) on moment-angle
complexes~$\zk$.

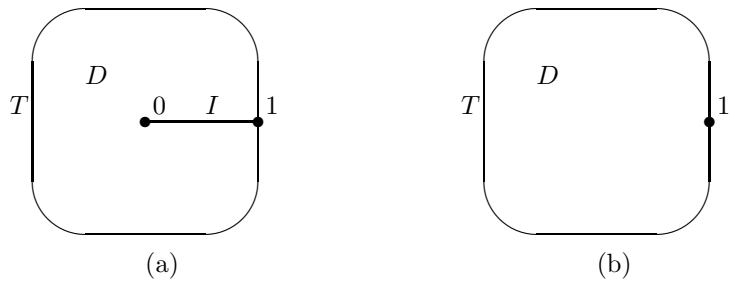
\begin{figure}
  \begin{picture}(120,35)
  \put(30,20){\oval(30,30)}
  \put(30,20){\line(1,0){15}}
  \put(30,20){\circle*{1.5}}
  \put(45,20){\circle*{1.5}}
  \put(31,21){0}
  \put(46,21){1}
  \put(38,21){$I$}
  \put(12,21){$T$}
  \put(22,25){$D$}
  \put(30,0){(a)}
  \put(90,20){\oval(30,30)}
  \put(105,20){\circle*{1.5}}
  \put(106,21){1}
  \put(72,21){$T$}
  \put(82,25){$D$}
  \put(90,0){(b)}
  \end{picture}
  \caption{Cellular decompositions of $D^2$.}
\end{figure}
Let us consider the cellular decomposition of $D^2$ with one 2-cell~$D$, two
1-cells $I$,~$T$, and two 0-cells 0,~1, see Figure~7~(a). It defines a
cellular decomposition of the poly-disk $(D^2)^m$ with $5^m$
cells. Each cell of this cellular complex is the product of cells of 5
different types: $D_i$, $I_i$, $0_i$, $T_i$ and $1_i$, $i=1,\ldots,m$. We
encode the cells of $(D^2)^m$ by ``words" of type $D_II_J0_LT_P1_Q$, where
$I,J,L,P,Q$ are pairwise disjoint subsets of $[m]$ such that $I\cup J\cup
L\cup P\cup Q=[m]$. Sometimes we would drop the last factor~$1_Q$, so in our
notations $D_II_J0_LT_P=D_II_J0_LT_P1_{[m]\setminus I\cup J\cup L\cup P}$.
The closure of $D_II_J0_LT_P1_Q$ is homeomorphic to the product
of $\#I$ disks, $\#J$ segments, and $\#P$ circles. The constructed cellular
decomposition of $(D^2)^m$ allows to display moment-angle complexes as
certain cellular subcomplexes in~$(D^2)^m$.

\begin{lemma}
\label{macell}
  For any cubical subcomplex $\mathcal C$ of~$I^m$ the corresponding
  moment-angle complex $\ma(\mathcal C)$ is a cellular
  subcomplex of~$(D^2)^m$.
\end{lemma}
\begin{proof}
Indeed, $\ma(\mathcal C)$ is a union of ``moment-angle" blocks $B_{I\subset
J}$~(\ref{bij}), and each $B_{I\subset J}$ is the closure of the cell
$D_{J\setminus I}I_\emptyset0_IT_{[m]\setminus J}1_\emptyset$.
\end{proof}

Now we concentrate on the moment-angle complex $\zk$ corresponding to the
cubical complex $\cc(K)\subset I^m$ (see~(\ref{zkwk})). By definition, $\zk$
is the union of moment-angle blocks $B_{I\subset J}\subset(D^2)^m$ with $J\in
K$. Denote
\begin{equation}
\label{bj}
  B_J:=B_{\emptyset\subset J}=\bigl\{(z_1,\ldots,z_m)\in
  (D^2)^m\: : \: |z_j|=1\text{ for }j\notin J\bigl\}.
\end{equation}
Then $B_J=\rho^{-1}(C_J)$ (recall that $C_J:=C_{\emptyset\subset J}$) and
$B_{I\subset J}\subset B_J$ for any $I\subset J$. It follows that
\begin{equation}
\label{zk=bigcu}
  \zk=\bigcup_{J\in K}B_J
\end{equation}
(compare this with the note after~(\ref{fcck})).

\begin{remark}
If $K=K_P$ for a simple polytope $P$ and $\#J=n$, then $B_J$ is $i_e(B_v)$
for $v=\bigcap_{j\in J}F_j$. Hence,~(\ref{zk=bigcu}) reduces to~(\ref{mazp})
in this case.
\end{remark}

Note that $B_J\cap B_{J'}=B_{J\cap J'}$.  This allows to simplify the
cellular decomposition from Lemma~\ref{macell} in the case $\ma(\mathcal
C)=\zk$. To do this we replace the union of cells $0$, $I$, $D$ (see
Figure~7~(a)) by one 2-dimensional cell (which we also denote~$D$). The
resulting cellular decomposition of $D^2$ with 3 cells is shown on
Figure~7~(b). It defines a cellular decomposition of $(D^2)^m$
with $3^m$ cells, each of which is the product of 3 different types of cells:
$D_i$, $T_i$ and $1_i$, $i=1,\ldots,m$. We encode these cells of $(D^2)^m$
as $D_IT_P1_Q$, where $I,P,Q$ are pairwise disjoint subsets of $[m]$ such
that $I\cup P\cup Q=[m]$. We use the notation $D_IT_P:=D_IT_P1_{[m]\setminus
I\cup P}$. The closure of $D_IT_P$ is the product of $\#I$ disks and $\#P$
circles.

\begin{lemma}
\label{zkcell}
  The moment-angle complex $\zk$ is a cellular
  subcomplex of $(D^2)^m$ with respect to the $3^m$-cell decomposition (see
  Figure~{\rm7~(b)}). Each cell of $\zk$ has the form $D_IT_P$, $I\in K$.
\end{lemma}
\begin{proof}
Since $B_J=B_{\emptyset\subset J}$ is the closure of the cell
$D_JT_{[m]\setminus J}1_\emptyset$, the statement follows
from~(\ref{zk=bigcu}).
\end{proof}

\begin{remark}
Note that for general $\mathcal C$ the moment-angle complex $\ma(\mathcal C)$
is {\it not\/} a cellular subcomplex for the $3^m$-cell decomposition
of~$(D^2)^m$.
\end{remark}

It follows from Construction~\ref{cck} that the cubical complex $\cc(K)$
always contains the vertex $(1,\ldots,1)\in I^m$. Hence, the torus
$T^m=\rho^{-1}(1,\ldots,1)$ is contained in~$\zk$.

\begin{lemma}
\label{the incl}
  The inclusion $T^m=\rho^{-1}(1,\ldots,1)\hookrightarrow\zk$ is a cellular
  map homotopical to the map to a point, i.e. the torus
  $T^m=\rho^{-1}(1,\ldots,1)$ is a contractible cellular subcomplex of~$\zk$.
\end{lemma}
\begin{proof}
To prove that $T^m=\rho^{-1}(1,\ldots,1)$ is a cellular subcomplex of $\zk$
we just mention that it is the closure of the $m$-cell $D_\emptyset
T_{[m]}\subset\zk$. So, it remains to prove that $T^m$ is contractible
within~$\zk$. To do this we show that the embedding $T^m\subset(D^2)^m$ is
homotopic to the map to the point $(1,\ldots,1)\in T^m\subset(D^2)^m$.  On
the first step we note that $\zk$ contains the cell $D_1T_{2,\ldots,m}$,
whose closure contains $T^m$ and is homeomorphic to $D^2\times T^{m-1}$.
Hence, our $T^m$ can be contracted to $1\times T^{m-1}$ within~$\zk$. On the
second step we note that $\zk$ contains the cell $D_2T_{3,\ldots,m}$, whose
closure contains $1\times T^{m-1}$ and is homeomorphic to $D^2\times
T^{m-2}$. Hence, $1\times T^{m-1}$ can be contracted to $1\times 1\times
T^{m-2}$ within~$\zk$, and so on. On the $k$-th step we note that $\zk$
contains the cell $D_kT_{k+1,\ldots,m}$, whose closure contains
$1\times\cdots\times1\times T^{m-k+1}$ and is homeomorphic to $D^2\times
T^{m-k}$. Hence, $1\times\cdots\times1\times T^{m-k+1}$ can be contracted to
$1\times\cdots\times1\times T^{m-k}$ within~$\zk$. We end up at the point
$1\times\cdots\times1$ to which the whole torus $T^m$ can be contracted.
\end{proof}

\begin{corollary}
  For any simplicial complex $K$ the moment-angle complex $\zk$ is simply
  connected.
\end{corollary}
\begin{proof}
  Indeed, the 1-skeleton of our cellular decomposition of $\zk$ is contained
  in the torus $T^m=\rho^{-1}(1,\ldots,1)$.
\end{proof}

\subsection{Borel construction and Stanley--Reisner space}
\label{hom1}
Let $ET^m$ be the contractible space of the universal principal $T^m$-bundle
over the classifying space~$BT^m$. It is well known that $BT^m$ is (homotopy
equivalent to) the product of $m$ copies of infinite-dimensional projective
space~$\C P^\infty$. The cellular decomposition of $\C P^\infty$ with one
cell in every even dimension determines the {\it canonical\/} cellular
decomposition of~$BT^m$. It follows that the cohomology of $BT^m$ (with
coefficients in $\k$) is the polynomial ring $\k[v_1,\ldots,v_m]$, $\deg
v_i=2$.

\begin{definition}
\label{borel}
Let $X$ be a space with an action of the torus~$T^m$. The {\it Borel
construction\/} (alternatively, {\it homotopy quotient\/} or {\it associated
bundle\/}) is the identification space
$$
  ET^m\times_{T^m} X:=ET^m\times X/{\sim},
$$
where $(e,x)\sim(eg,g^{-1}x)$ for any $e\in ET^m$, $x\in X$, $g\in T^m$.
\end{definition}

The projection $(e,x)\to e$ displays $ET^m\times_{T^m} X$ as the total space
of a bundle $ET^m\times_{T^m} X\to BT^m$ with fibre $X$ and structure
group~$T^m$. At the same time, there is a principal $T^m$-bundle $ET^m\times
X\to ET^m\times_{T^m} X$.

In the sequel we denote the Borel construction $ET^m\times_{T^m} X$
corresponding to a $T^m$-space $X$ by~$B_TX$. In particular, for any
simplicial complex $K$ on $m$ vertices there are defined the Borel
construction $\bk$ and the bundle $p:\bk\to BT^m$ with fibre~$\zk$.

For each $i=1,\ldots,m$ denote by $BT_i$ the $i$-th coordinate subspace of
$BT^m=(\C P^\infty)^m$. For any subset $I\subset[m]$, we denote by
$BT_I$ the product of $BT_i$ with $i\in I$. Obviously, $BT_I$ is a cellular
subcomplex of~$BT^m$, and $BT_I\cong BT^k$ provided that $\#I=k$.

\begin{definition}
\label{srspace}
Let $K$ be a simplicial complex. We refer to the cellular subcomplex
$$
  \bigcup_{I\in K}BT_I\subset BT^m
$$
as the {\it Stanley--Reisner space\/}, and denote it~$SR(K)$.
\end{definition}
\noindent The name refers to the following statement, which is an immediate
corollary of the definition of Stanley--Reisner ring $\k(K)$
(Definition~\ref{frsim}).

\begin{proposition}
\label{homsrs}
The cellular cochain algebra $C^*(SR(K))$ and the cohomology algebra
$H^*(SR(K))$ are isomorphic to the face ring~$\k(K)$. The cellular inclusion
$i:SR(K)\hookrightarrow BT^m$ induces the quotient epimorphism
$i^*:\k[v_1,\ldots,v_m]\to\k(K)=\k[v_1,\ldots,v_m]/\mathcal I_K$ in the
cohomology.
\end{proposition}

\begin{theorem}
\label{homeq1}
  The fibration $p:\bk\to BT^m$ is homotopy equivalent to the cellular
  inclusion $i:SR(K)\hookrightarrow BT^m$. More precisely, there is a
  deformation retraction $\bk\to SR(K)$ such that the diagram
  $$
  \begin{CD}
    \bk @>p>> BT^m\\
    @VVV @|\\
    SR(K) @>i>> BT^m
  \end{CD}
  $$
  is commutative.
\end{theorem}
\begin{proof}
Consider the decomposition~(\ref{zk=bigcu}). Since each $B_J\subset\zk$ is
$T^m$-stable, the Borel construction $\bk=ET^m\times_{T^m}\zk$ is patched
from the Borel constructions $ET^m\times_{T^m}B_J$, $J\in K$.  Suppose
$\#J=j$; then $B_J\cong(D^2)^j\times T^{m-j}$ (see~(\ref{bj})). By definition
of Borel construction,
$ET^m\times_{T^m}B_J\cong(ET^j\times_{T^j}(D^2)^j)\times ET^{m-j}$. The space
$ET^j\times_{T^j}(D^2)^j$ is the total space of a $(D^2)^j$-bundle
over~$BT^j$. It follows that there is a deformation retraction
$ET^m\times_{T^m}B_J\to BT_J$, which defines a homotopy equivalence between
the restriction of $p:\bk\to BT^m$ to $ET^m\times_{T^m}B_J$ and the cellular
inclusion $BT_J\hookrightarrow BT^m$. These homotopy equivalences
corresponding to different simplices $J\in K$ fit together to yield a
required homotopy equivalence between $p:\bk\to BT^m$ and
$i:SR(K)\hookrightarrow BT^m$.
\end{proof}

\begin{corollary}
\label{homfib}
  The moment-angle complex $\zk$ is a homotopy fibre of the cellular
  inclusion $i:SR(K)\hookrightarrow BT^m$.
\end{corollary}

\begin{corollary}
\label{cohombk}
The cohomology algebra $H^*(\bk)$ is isomorphic to the face ring~$\k(K)$. The
projection $p:\bk\to BT^m$ induces the quotient epimorphism
$p^*:\k[v_1,\ldots,v_m]\to\k(K)=\k[v_1,\ldots,v_m]/\mathcal I_K$ in the
cohomology.
\end{corollary}

\begin{remark}
The above statement was proved in~\cite[Theorem~4.8]{DJ} in the polytopal case
($\zk=\zp$) by other methods.
\end{remark}

The following information about the homotopy groups of $\zk$ can be retrieved
from the above constructions.

\begin{theorem}
\label{homgr}
  {\rm(a)}  The complex $\zk$ is {\rm2}-connected (i.e.
  $\pi_1(\zk)=\pi_2(\zk)=0$), and $\pi_i(\zk)=\pi_i(\bk)=\pi_i(SR(K))$ for
  $i\ge 3$.

  {\rm(b)}  If $K=K_P$ and $P$ is $q$-neighbourly, then $\pi_i(\zk)=0$ for
  $i<2q+1$, and $\pi_{2q+1}(\zp)$ is a free Abelian group whose
  generators correspond to the $(q+1)$-element missing faces of $K_P$ (or
  equivalently, to degree-$(2q+2)$ generators of the ideal $\mathcal I_P$,
  see Definition~{\rm\ref{frpol}}).
\end{theorem}
\begin{proof}
Note that $BT^m=K(\Z^m,2)$ and the 3-skeleton of $SR(K)$ coincides with that
of~$BT^m$. If $P$ is $q$-neighbourly, then it follows from
Definition~\ref{srspace} that the $(2q+1)$-skeleton of $SR(K_P)$ coincides
with that of~$BT^m$. Now, both statements follow easily from the exact
homotopy sequence of the map $i:SR(K)\to BT^m$ with homotopy fibre $\zk$
(Corollary~\ref{homfib}).
\end{proof}

\begin{remark}
Say that a simplicial complex $K$ on the set $[m]$ is {\it $k$-neighbourly\/}
if any $k$-element subset of $[m]$ is a simplex of~$K$. (This definition is
an obvious extension of the notion of $k$-neighbourly simplicial polytope to
arbitrary simplicial complexes). Then the second part of Theorem~\ref{homgr}
holds for arbitrary $q$-neighbourly simplicial complex.
\end{remark}

Suppose now that $K=K_P$ for some simple $n$-polytope $P$ and $M^{2n}$ is a
quasitoric manifold over $P$ with characteristic function~$\ell$. Then we
have the subgroup $H(\ell)\subset T^m$ acting freely on $\zp$ and
the principal $T^{m-n}$-bundle $\zp\to M^{2n}$ (Proposition~\ref{zpbun}).

\begin{proposition}
\label{qtbc}
The Borel construction $ET^n\times_{T^n}M^{2n}$ is homotopy equivalent
to~$\bp$.
\end{proposition}
\begin{proof}
Since $H(\ell)$ acts freely on $\zp$, we have
\begin{multline*}
  \bp=ET^m\times_{T^m}\zp\\
  =EH(\ell)\times
  \Bigl(E\bigl(T^m/H(\ell)\bigr)\times_{T^m/H(\ell)}\zp/H(\ell)\Bigr)
  \simeq ET^n\times_{T^n}M^{2n}.
\end{multline*}
\end{proof}

\begin{theorem}[{\cite[Theorem~4.12]{DJ}}]
\label{qtdeg}
The Leray--Serre spectral sequence of the bundle
\begin{equation}
\label{qtbun1}
  ET^n\times_{T^n}M^{2n}\to BT^n
\end{equation}
with fibre $M^{2n}$ collapses at the $E_2$ term, i.e.
$E_2^{p,q}=E_\infty^{p,q}$.
\end{theorem}
\begin{proof}
The differentials of the spectral sequence are all trivial by dimensional
reasons, since both $BT^n$ and $M^{2n}$ have cells only in even dimensions
(see Proposition~\ref{qtbn}).
\end{proof}

\begin{corollary}
\label{qtmaps}
The projection~{\rm(\ref{qtbun1})} induces a monomorphism
$\k[t_1,\ldots,t_n]\to\k(P)$ in the cohomology. The inclusion of the fibre
$M^{2n}\hookrightarrow ET^n\times_{T^n}M^{2n}$ induces an epimorphism
$\k(P)\to H^*(M^{2n})$.
\end{corollary}

Now we are ready to prove the statements from section~\ref{quas}.

\begin{proof}[Proof of Lemma~{\rm\ref{crs}} and Theorem~{\rm\ref{qtcoh}}]
The monomorphism
$$
  H^*(BT^n)=\k[t_1,\ldots,t_n]\to\k(P)=H^*(ET^n\times_{T^n}M^{2n})
$$
takes $t_i$ to $\theta_i$, $i=1,\ldots,n$. Since $\k(P)$ is a free
$\k[t_1,\ldots,t_n]$-module (note that this follows from Theorem~\ref{qtdeg},
so we do not need to use Theorem~\ref{reisner}), $\t_1,\ldots,\t_n$ is a
regular sequence. It follows that the kernel of $\k(P)\to H^*(M^{2n})$ is
exactly $\mathcal J_\ell=(\t_1,\ldots,\t_n)$.
\end{proof}

\subsection{Generalisations, analogues and additional comments}
\label{gen2}
Many important constructions from our survey (namely, the cubical
complex~$\cc(K)$, the moment-angle complex~$\zk$, the Borel
construction~$\bk$, the Stanley--Reisner space~$SR(K)$, and also the
complement $U(K)$ of a coordinate subspace arrangement from
section~\ref{coor}) admit the unifying combinatorial interpretation in terms
of the following construction, which was proposed by N.~Strickland.

\begin{construction}\label{nsc}
Let $X$ be a space, and $W$ a subspace of~$X$. Let
$K$ be a simplicial complex on the vertex set~$[m]$. Define the following
subset of the product of $m$ copies of~$X$:
$$
  K_\bullet(X,W)=\bigcup_{I\in K}\Bigl(\prod_{i\in I}X\times
  \prod_{i\notin I}W \Bigl).
$$
\end{construction}

\begin{example}
1. $\cc(K)=K_\bullet(I^1,1)$ (see~(\ref{fcck})).

2. $\zk=K_\bullet(D^2,S^1)$ (see~(\ref{zk=bigcu})).

3. $SR(K)=K_\bullet(\C P^\infty,{*})$ (see Definition~\ref{srspace}).

4. $\bk=K_\bullet(ES^1\times_{S^1}D^2,ES^1\times_{S^1}D^2)$
(see the proof of Theorem~\ref{homeq1}).
\end{example}

Note that Construction~\ref{nsc} is obviously extended to the set of
$m$ pairs $(X_1,W_1),\ldots,(X_m,W_m)$; another generalisation is obtained by
replacing the cartesian product by the fibred product.

Using construction \ref{nsc}, the part of Theorem~\ref{cubkcck} (dealing with
the complex~$\cc(K)$) and Theorem~\ref{he1} were obtained independently by
N.~Strickland.

Almost all constructions of our paper incorporate some action of the torus,
i.e. the product of circles~$S^1$. These constructions admit natural
$\Z/2$-analogues. To see this we replace the torus $T^m$ by its ``real
analogue", that is, the group~$(\Z/2)^m$. Then the standard cube
$I^m=[0,1]^m$ is the orbit space for the action of $(\Z/2)^m$
on the bigger cube $[-1,1]^m$, which can be regarded as a ``real analogue" of
the poly-disk $(D^2)^m\subset\C^m$. Now, given a cubical subcomplex $\mathcal
C\subset I^m$, we can construct a $(\Z/2)^m$-symmetrical cubical complex
embedded into $[-1,1]^m$ just in the same way as we did it in
Definition~\ref{ma}. In particular, for any simplicial complex $K$ on the
vertex set~$[m]$ we can introduce the cubical complexes $\R\zk$ and $\R\wk$,
the ``real analogues" of the moment-angle complexes $\zk$ and
$\wk$~(\ref{zkwk}).  In the notations of Construction~\ref{nsc} we have
$$
  \R\zk=K_\bullet\bigl([-1,1],\{-1,1\}\bigr).
$$
If $K$ is a simplicial $(n-1)$-sphere, then $\R\zk$ is an $n$-dimensional
manifold (the proof is identical to that of Lemma~\ref{maman}).  Thus, for
any simplicial sphere $K^{n-1}$ with $m$ vertices we get a
$(\Z/2)^m$-symmetric $n$-manifold with a $(\Z/2)^m$-invariant cubical
subdivision. This class of cubical manifolds may be useful for
the cubical analogue for the combinatorial theory of $f$-vectors of simplicial
complexes (see also~\cite{St5}). Finally, the real analogue $\R\zp$ of the
manifold $\zp$ (corresponding to the case of polytopal simplicial
sphere~$K^{n-1}$) is known as the {\it universal Abelian cover\/} of the
polytope~$P^n$ regarded as an orbifold (or manifold with corners),
see~\cite[\S4.5]{Gr}. In~\cite{Iz1} manifolds $\R\zp$ and $\zp$ are
interpreted as the configuration spaces of {\it hinge mechanisms\/} in $\R^2$
and~$\R^3$.

Passage from $T^n$ to $(\Z/2)^n$ in Definition~\ref{qtm} leads to real
analogues of quasitoric manifolds, which were introduced in~\cite{DJ} under
the name {\it small covers\/}. Thus, every small cover of a simple polytope
$P^n$ is a manifold $M^{n}$ with an action of $(\Z/2)^n$ and quotient~$P^n$.
The name refers to the fact that any branched cover of $P^n$ by a smooth
manifold have at least $2^n$ sheets. Small covers were studied in~\cite{DJ}
along with quasitoric manifolds, and many results on quasitoric manifolds
cited from~\cite{DJ} in section~\ref{quas} have analogues in the small cover
case. On the other hand, every small cover is the quotient of the universal
cover $\R\zp$ by a certain free action of the group~$(\Z/2)^{m-n}$.

An important class of small covers over 3-dimensional simple polytopes
$P^3$ was considered in~\cite{Iz2}. It can be shown that a simple polytope
$P^3$ admits a 3-paint colouring of its facets such that any two adjacent
facets have different colour if and only if every facet have even number of
edges. Every such colouring $\varrho$ defines the quasitoric manifold
$M^{6}(\varrho)$ and the small cover~$M^3(\varrho)$. It was shown
in~\cite{Iz2} that every manifold $M^3(\varrho)$ admits an equivariant
embedding into $\R^4=\R^3\times\R$ with the standard action of
$(\Z/2)^3$ on $\R^3$ and the trivial action on~$\R$. It was also shown there
that all manifolds $M^3(\varrho)$ can be obtained from the set of
3-dimensional tori by applying the operations of equivariant connected sum and
{\it equivariant Dehn twist\/}.

The quaternionic analogue of moment-angle complexes can be constructed by
replacing $T^n$ by the quaternionic torus $Sp(1)^n\cong(S^3)^n$. Developing
the quaternionic analogues of toric and quasitoric manifolds is a problem of
particular interest. This is acknowledged, in particular, by the important
results of~\cite{BGMR}.

At the end we give one application of the above described constructions to
the case of general group~$G$.

\begin{example}[classifying space for group~$G$]
Let $K$ be a simplicial complex on the vertex set~$[m]$. Set
$\zk(G):=K_\bullet(\cone(G),G)$ (see Construction~\ref{nsc}), where
$\cone(G)$ is the cone over~$G$. By the construction, the group $G^m$ acts on
$\zk(G)$ with quotient~$\cone(K)$. The diagonal embedding $G\hookrightarrow
G^m$ defines a free action of $G$ on~$\zk(G)$.

Suppose now that $K_0\subset K_1\subset\dots\subset K_i\subset\cdots$ is a
sequence of embedded simplicial complexes such that $K_i$ is
$q_i$-neighbourly and $q_i\to\infty$ as $i\to\infty$. Such a sequence can be
constructed, for instance, by taking $K_{i+1}:=K_i*K$ (see
Construction~\ref{join}), where $K_0$ and $K$ are arbitrary simplicial
complexes. The group $G$ acts freely on the space
$\lim\limits_{\longrightarrow}\mathcal Z_{K_i}(G)$, and the corresponding
quotient provides a realisation of the classifying space~$BG$. Thus, we have
the following filtration in the universal fibration $EG\to BG$:
$$
\begin{array}{ccccccccc}
\mathcal Z_{K_0}(G) & \hookrightarrow & \mathcal Z_{K_1}(G) &
\hookrightarrow & \cdots & \hookrightarrow & \mathcal Z_{K_i}(G) &
\hookrightarrow &\cdots\\
\downarrow && \downarrow &&&&\downarrow\\
\mathcal Z_{K_0}(G)/G & \hookrightarrow & \mathcal Z_{K_1}(G)/G &
\hookrightarrow & \cdots & \hookrightarrow & \mathcal Z_{K_i}(G)/G &
\hookrightarrow & \cdots.
\end{array}
$$
The well-known {\it Milnor filtration\/} in the universal fibration of the
group $G$ corresponds to the case $K_i=\D^i$.
\end{example}

\section{Cohomology of moment-angle complexes and combinatorics of
simplicial manifolds}
\subsection{Eilenberg--Moore spectral sequence}
\label{eile}
In their 1966 paper~\cite{EM} Eilenberg and Moore constructed a
spectral sequence of great importance for algebraic topology. This spectral
sequence can be considered as an extension of Adams' approach to calculating
the cohomology of loop spaces~\cite{Ad0}. In 1960-70s the Eilenberg--Moore
spectral sequence allowed to obtain many important results on the cohomology
of loop spaces and homogeneous spaces for Lie group actions. In our paper we
describe new applications of this spectral sequence to combinatorial
problems. This section contains the necessary information about the spectral
sequence; we follow Smith's paper~\cite{Sm} in this description.

To be precise, there are two Eilenberg--Moore spectral sequences, the
algebraic and the topological.

\begin{theorem}[{Eilenberg--Moore \cite[Theorem~1.2]{Sm}}]
\label{algemss}
Let $A$ be a differential graded $\k$-algebra, and $M$,~$N$ differential
graded $A$-modules. Then there exists a spectral sequence
$\{E_r,d_r\}$ that converges to $\Tor_A(M,N)$ and has the $E_2$-term
$$
  E_2^{-i,j}=\Tor^{-i,j}_{H[A]}\bigl(H[M],H[N]\bigr),\quad i,j\ge0,
$$
where $H[\cdot]$ denotes the cohomology algebra (module).
\end{theorem}

The spectral sequence of Theorem~\ref{algemss} lives in the {\it second\/}
quadrant and the differential $d_r$ adds $(r,1-r)$ to bidegree. This spectral
sequence is called the (algebraic) {\it Eilenberg--Moore spectral
sequence\/}. It gives a decreasing filtration $\{F^{-p}\Tor_A(M,N)\}$ on
$\Tor_A(M,N)$ with the property that
$$
  E_\infty^{-p,n+p}=F^{-p}\biggl( \sum_{-i+j=n}\Tor^{-i,j}_A(M,N) \biggr)
  \biggl/ F^{-p+1}\biggl( \sum_{-i+j=n}\Tor^{-i,j}_A(M,N) \biggr).
$$

Topological applications of Theorem~\ref{algemss} arise in the case when
$A,M,N$ are singular (or cellular) cochain algebras of certain topological
spaces. The classical situation is described by the commutative diagram
\begin{equation}
\begin{CD}
  E @>>> E_0\\
  @VVV @VVV\\
  B @>>> B_0,
\end{CD}
\label{comsq}
\end{equation}
where $E_0\to B_0$ is a Serre fibre bundle with fibre $F$ over the simply
connected base~$B_0$, and $E\to B$ is the pullback along a continuous map
$B\to B_0$. For any space~$X$, let $C^*(X)$ denote either the singular
cochain algebra of $X$ or (in the case when $X$ is a cellular complex) the
cellular cochain algebra of~$X$. Obviously, $C^*(E_0)$ and $C^*(B)$ are
$C^*(B_0)$-modules. Under these assumptions the following statement holds.

\begin{lemma}[{\cite[Proposition~3.4]{Sm}}]
\label{gentoralg}
$\Tor_{C^*(B_0)}(C^*(E_0),C^*(B))$ is an algebra in a natural way, and the
is a canonical isomorphism of algebras
$$
  \Tor_{C^*(B_0)}\bigl(C^*(E_0),C^*(B)\bigr)\to H^*(E).
$$
\end{lemma}

Applying Theorem~\ref{algemss} in the case $A=C^*(B_0)$, $M=C^*(E_0)$,
$N=C^*(B)$ and taking into account Lemma~\ref{gentoralg} we obtain
\begin{theorem}[{Eilenberg--Moore}]
\label{topemss}
There exists a spectral sequence of commutative algebras $\{E_r,d_r\}$ with
\begin{enumerate}
\item[(a)] $E_r\Rightarrow H^*(E)$;
\item[(b)] $E_2^{-i,j}=\Tor^{-i,j}_{H^*(B_0)}(H^*(E_0),H^*(B))$.
\end{enumerate}
\end{theorem}

The spectral sequence of Theorem~\ref{topemss} is called the (topological)
{\it Eilenberg--Moore spectral sequence\/}. In the very important
particular case when $B$ is a point (see~(\ref{comsq})) we get
\begin{corollary}
\label{onefib}
  Let $E\to B$ be a fibration over the simply connected space $B$ with
  fibre~$F$. Then there exists a spectral sequence of commutative algebras
  $\{E_r,d_r\}$ with
  \begin{enumerate}
    \item[(a)] $E_r\Rightarrow H^*(F)$,
    \item[(b)] $E_2=\Tor_{H^*(B)}(H^*(E),\k)$.
  \end{enumerate}
\end{corollary}
We refer to the spectral sequence of Corollary~\ref{topemss} as the
{\it Eilenberg--Moore spectral sequence of fibration $E\to B$\/}.

\begin{example}\label{qtem}
Let $M^{2n}$ be a quasitoric manifold over~$P^n$.  Consider the
Eilenberg--Moore spectral sequence of the bundle $ET^n\times_{T^n}M^{2n}\to
BT^n$ with fibre~$M^{2n}$. By Proposition~\ref{qtbc},
$H^*(ET^n\times_{T^n}M^{2n})=H^*(\bp)\cong\k(P^n)$. The monomorphism
$$
  \k[t_1,\ldots,t_n]=H^*(BT^n)\to H^*(ET^n\times_{T^n}M^{2n})=\k(P^n)
$$
takes $t_i$ to $\theta_i$, $i=1,\ldots,n$, see~(\ref{theta}). The
$E_2$ term of the Eilenberg--Moore spectral sequence is
$$
  E_2^{*,*}=\Tor^{*,*}_{H^*(BT^n)}\bigl(H^*(ET^n\times_{T^n}M^{2n}),\k\bigr)=
  \Tor^{*,*}_{\k[t_1,\ldots,t_n]}\bigl(\k(P^n),\k\bigr).
$$
Since $\k(P^n)$ is a free $\k[t_1,\ldots,t_n]$-module, we have
\begin{multline*}
  \Tor^{*,*}_{\k[t_1,\ldots,t_n]}\bigl(\k(P^n),\k\bigr)=
  \Tor^{0,*}_{\k[t_1,\ldots,t_n]}\bigl(\k(P^n),\k\bigr)\\
  =\k(P^n)\otimes_{\k[t_1,\ldots,t_n]}\k=\k(P^n)/(\t_1,\ldots,\t_n).
\end{multline*}
Therefore, $E_2^{0,*}=\k(P^n)/\mathcal J_\ell$ and $E_2^{-p,*}=0$ for $p>0$.
It follows that the Eilenberg--Moore spectral sequence collapses at the $E_2$
term and $H^*(M^{2n})=\k(P^n)/J_\ell$, in accordance with
Theorem~\ref{qtcoh}.
\end{example}

\subsection{Cohomology of moment-angle complex $\zk$: the case of general $K$}
\label{coh1}
Here we apply the Eilenberg--Moore spectral sequence to calculating the
cohomology algebra of the moment-angle complex~$\zk$. This describes
$H^*(\zk)$ as a {\it bigraded algebra\/}. The corresponding bigraded Betti
numbers are important combinatorial invariants of~$K$.

\begin{theorem}
\label{cohom1}
  The following isomorphism of graded algebras holds:
  $$
    H^{*}(\zk)\cong\Tor_{\k[v_1,\ldots,v_m]}\bigl(\k(K),\k\bigr).
  $$
  In particular,
  $$
    H^p(\zk)\cong
    \sum_{-i+2j=p}\Tor^{-i,2j}_{\k[v_1,\ldots,v_m]}\bigl(\k(K),\k\bigr).
  $$
\end{theorem}
\begin{proof}
Let us consider the Eilenberg--Moore spectral sequence of the commutative
square
\begin{equation}
\label{wuk}
  \begin{CD}
    E @>>> ET^m\\
    @VVV @VVV\\
    SR(K) @>i>> BT^m,
  \end{CD}
\end{equation}
where the left vertical arrow is the pullback along~$i$.
Corollary~\ref{homfib} shows that $E$ is homotopy equivalent to~$\zk$.

By Proposition~\ref{homsrs}, the map $i:SR(K)\hookrightarrow BT^m$ induces
the quotient epimorphism
$i^*:C^*(BT^m)=\k[v_1,\ldots,v_m]\to\k(K)=C^*(SR(K))$, where $C^*(\cdot)$
denotes the cellular cochain algebra. Since $ET^m$ is contractible, there is
a chain equivalence $C^{*}(ET^m)\simeq\k$. Therefore, there is an isomorphism
\begin{equation}
\label{chain}
 \Tor_{C^{*}(BT^m)}\bigl(C^{*}(SR(K)),C^{*}(ET^m)\bigr)\cong
 \Tor_{\k[v_1,\ldots,v_m]}\bigl(\k(K),\k\bigr).
\end{equation}

The Eilenberg--Moore spectral sequence of commutative square~(\ref{wuk}) has
$$
  E_2=\Tor_{H^{*}(BT^m)}\bigl(H^{*}(SR(K)),H^{*}(ET^m)\bigr)
$$
and converges to $\Tor_{C^{*}(BT^m)}(C^{*}(SR(K)),C^{*}(ET^m))$
(Theorem~\ref{algemss}). Since
$$
  \Tor_{H^{*}(BT^m)}\bigl(H^{*}(SR(K)),H^{*}(ET^m)\bigr)=
  \Tor_{\k[v_1,\ldots,v_m]}\bigl(\k(K),\k\bigr),
$$
it follows from~(\ref{chain}) that the spectral sequence collapses at the
$E_2$ term, that is, $E_2=E_{\infty}$. Lemma~\ref{gentoralg} shows that the
module $\Tor_{C^{*}(BT^m)}(C^{*}(SR(K)),C^{*}(ET^m))$ is an algebra
isomorphic to $H^{*}(\zk)$, which concludes the proof.
\end{proof}

Theorem~\ref{cohom1} describes the cohomology of $\zk$ as a {\it bigraded\/}
algebra and shows that the corresponding bigraded Betti numbers
$b^{-i,2j}(\zk)$ coincide with that of~$\k(K)$, see~(\ref{bbnfr}). Using the
Koszul resolution for $\k$ and Lemma~\ref{koscom}, we get
\begin{theorem}
\label{cohom2}
  The following isomorphism of bigraded algebras holds:
  $$
    H^{*,*}(\zk)\cong
    H\bigl[\Lambda[u_1,\ldots,u_m]\otimes\k(K),d\bigr],
  $$
  where the bigraded structure and the differential in the right hand side are
  defined by~{\rm(\ref{diff})}.
\end{theorem}

In the sequel we denote square-free monomials
$$
  u_{i_1}\dots u_{i_p}v_{j_1}\dots v_{j_q}\in
  \bigl[\Lambda[u_1,\ldots,u_m]\otimes\k(K),d\bigr]
$$
by $u_Iv_J$, where $I=\{i_1,\ldots,i_p\}$, $J=\{j_1,\ldots,j_q\}$. Note that
$\bideg u_Iv_J=(-p,2(p+q))$.

\begin{remark}
Since the differential $d$ in~(\ref{diff}) does not change the second
degree, the differential bigraded algebra
$[\Lambda[u_1,\ldots,u_m]\otimes\k(K),d]$ splits into the sum of
differential algebras consisting of elements of fixed second degree.
\end{remark}

\begin{corollary}
\label{ls}
  The Leray--Serre spectral sequence of the principal $T^m$-bundle
  $ET^m\times\zk\to\bk$ collapses at the $E_3$ term.
\end{corollary}
\begin{proof}
The spectral sequence under consideration converges to
$H^{*}(ET^m\times\zk)=H^{*}(\zk)$ and has
$$
  E_2=H^{*}(\bk)\otimes H^{*}(T^m)=\Lambda[u_1,\ldots,u_m]\otimes\k(K).
$$
The differential in the $E_2$ term acts as in~(\ref{diff}). Hence,
$$
  E_3=H[E_2,d]=H\bigl[\Lambda[u_1,\ldots,u_m]\otimes\k(K)\bigr]=
  H^{*}(\zk)
$$
(the last identity follows from Theorem~\ref{cohom2}).
\end{proof}

\begin{construction}
\label{astar}
Let $A^{-q}(K)\subset\L[u_1,\ldots,u_m]\otimes\k(K)$ be the subspace
generated by monomials $u_I$ and $u_Iv_J$ such that $J$ is a simplex of~$K$,
$\#I=q$ and $I\cap J=\emptyset$. Define
$$
  A^{*}(K)=\bigoplus_{q=0}^{m}A^{-q}(K).
$$
Since $d(u_i)=v_i$, we have $d\bigl(A^{-q}(K)\bigr)\subset A^{-q+1}(K)$.
Therefore, $A^{*}(K)$ is a cochain subcomplex in
$[\L[u_1,\ldots,u_m]\otimes\k(K),d]$. Moreover, $A^*(K)$ inherits the
bigraded module structure from $\L[u_1,\ldots,u_m]\otimes\k(K)$, with
differential $d$ adding $(1,0)$ to bidegree. Hence, we have the additive
inclusion (i.e. the monomorphism of bigraded modules)
$i_a:A^{*}(K)\hookrightarrow\L[u_1,\ldots,u_m]\otimes\k(K)$. Finally,
$A^*(K)$ is an algebra in the obvious way, but is {\it not\/} a subalgebra of
$\L[u_1,\ldots,u_m]\otimes\k(K)$. (Indeed, for instance, $v_1^2=0$ in
$A^*(K)$ but not in $\L[u_1,\ldots,u_m]\otimes\k(K)$.) Nevertheless, we have
the multiplicative {\it projection\/} (an epimorphism of bigraded algebras)
$j_m:\L[u_1,\ldots,u_m]\otimes\k(K)\to A^*(K)$. The additive inclusion $i_a$
and the multiplicative projection $j_m$ obviously satisfy $j_m\cdot i_a=\id$.
\end{construction}

\begin{lemma}
\label{iscoh}
  Cochain complexes $[\L[u_1,\ldots,u_m]\otimes\k(K),d]$
  and $[A^*(K),d]$ have the same cohomology. This implies the following
  isomorphism of bigraded $\k$-modules:
  $$
    H[A^{*}(K),d]\cong\Tor_{\k[v_1,\ldots,v_m]}\bigl(\k(K),\k\bigr).
  $$
\end{lemma}
\begin{proof}
A routine check shows that the cochain homotopy operator $s$ for the Koszul
resolution (see the proof of Proposition~VII.2.1 in~\cite{Mac})
establishes a cochain homotopy equivalence between the maps $\id$ and
$i_a\cdot j_m$ of the algebra $[\L[u_1,\ldots,u_m]\otimes\k(K),d]$ to itself.
That is,
$$
  ds+sd=\id-i_a\cdot j_m.
$$
We just illustrate the above identity on few simple examples.
\begin{align*}
  1)\quad& s(u_1v_2)=u_1u_2,\quad ds(u_1v_2)=u_2v_1-u_1v_2,\quad
  sd(u_1v_2)=u_1v_2-u_2v_1,\\
  &\text{hence,}\quad (ds+sd)(u_1v_2)=0=(\id-i_a\cdot j_m)(u_1v_2);\\
  2)\quad& s(u_1v_1)=u_1^2=0,\quad ds(u_1v_1)=0,\quad
  d(u_1v_1)=v_1^2,\quad sd(u_1v_1)=u_1v_1,\\
  &\text{hence,}\quad (ds+sd)(u_1v_1)=u_1v_1=(\id-i_a\cdot j_m)(u_1v_1);\\
  3)\quad& s(v_1^2)=u_1v_1,\quad ds(v_1^2)=v_1^2,\quad
  d(v_1^2)=0,\\
  &\text{hence,}\quad (ds+sd)(v_1^2)=v_1^2=(\id-i_a\cdot j_m)(v_1^2).
\end{align*}
\end{proof}

Now we recall our cellular decomposition of~$\zk$, see Lemma~\ref{zkcell}.
The cells are $D_IT_J$ with $J\subset[m]$, $I\in K$, and $I\cap J=\emptyset$.
Let $C_{*}(\zk)$ and $C^{*}(\zk)$ denote the corresponding cellular chain and
cochain complexes respectively. Both complexes $C^*(\zk)$ and $A^{*}(K)$ have
the same cohomology $H^{*}(\zk)$. The complex $C^{*}(\zk)$ has the basis
consisting of cochains $(D_IT_J)^{*}$. As an algebra, $C^*(\zk)$ is generated
by the cochains $T_i^{*}$, $D_i^{*}$ of dimension 1 and 2 respectively dual
to the cells $T_i$ and $D_i$, $i,j=1,\ldots,m$. At the same time, $A^*(K)$ is
multiplicatively generated by $u_i$, $v_i$, $i,j=1,\ldots,m$.

\begin{theorem} \label{cellcom}
  The correspondence $v_Iu_J\mapsto(D_IT_J)^{*}$ establishes a canonical
  isomorphism of differential graded algebras $A^*(K)$ and~$C^*(\zk)$.
\end{theorem}
\begin{proof}
It follows directly from the definitions of $A^*(K)$ and $C^{*}(\zk)$ that
the map is an isomorphism of graded algebras. Thus, it remains to prove that
it commutes with differentials. Let $d$, $d_c$ and $\partial_c$ denote the
differentials in $A^*(K)$, $C^{*}(\zk)$ and $C_*(\zk)$ respectively. Since
$d(v_i)=0$, $d(u_i)=v_i$, we need to show that $d_c(D_i^{*})=0$,
$d_c(T_i^{*})=D_i^{*}$.  We have $\partial_c(D_i)=T_i$, $\partial_c(T_i)=0$.
Any 2-cell of $\zk$ is either $D_j$ or $T_{jk}$, $k\ne j$. Then
$$
  (d_cT_i^{*},D_j)=(T_i^{*},\partial_cD_j)=(T_i^{*},T_j)=\delta_{ij},\quad
  (d_cT_i^{*},T_{jk})=(T_i^{*},\partial_cT_{jk})=0,
$$
where $\delta_{ij}=1$ if $i=j$ and $\delta_{ij}=0$ otherwise. Hence,
$d_c(T_i^{*})=D_i^{*}$. Further, any 3-cell of $\zk$
is either $D_jT_k$ or $T_{j_1j_2j_3}$. Then
\begin{align*}
  &(d_cD_i^{*},D_jT_k)=(D_i^{*},\partial_c(D_jT_k))=
  (D_i^{*},T_{jk})=0,\\
  &(d_cD_i^{*},T_{j_1j_2j_3})=(D_i^{*},\partial_cT_{j_1j_2j_3})=0.
\end{align*}
Hence, $d_c(D_i^{*})=0$.
\end{proof}
Theorem~\ref{cellcom} provides a topological interpretation for the
differential algebra $[A^*(K),d]$. In the sequel we do not distinguish the
cochain complexes $A^*(K)$ and $C^*(\zk)$, and identify $u_i$ with~$T_i^{*}$,
$v_i$ with~$D_i^{*}$.

Now we recall that the algebra $[A^*(K),d]$ is bigraded.  The isomorphism of
Theorem~\ref{cellcom} provides a bigraded structure for the cellular chain
complex $C_{*}(\zk),\partial_c]$, with
\begin{equation}\label{bgcellz}
  \bideg(D_i)=(0,2),\quad \bideg(T_i)=(-1,2),\quad \bideg(1_i)=(0,0).
\end{equation}
The differential $\partial_c$ adds $(-1,0)$ to bidegree, and the cellular
homology of $\zk$ also acquires a bigraded structure.

Let us assume now that the ground field $\k$ is of zero characteristic (e.g.
$\k=\Q$, the field of rational numbers). Define the bigraded Betti numbers
\begin{equation}
\label{bbn}
  b_{-q,2p}(\zk)=\dim H_{-q,2p}[C_*(\zk),\partial_c],
  \quad q,p=0,\ldots,m.
\end{equation}
Theorem~\ref{cellcom} and Lemma~\ref{iscoh} show that
\begin{equation}
\label{bbntor}
  b_{-q,2p}(\zk)=\dim\Tor^{-q,2p}_{\k[v_1,\ldots,v_m]}\bigl(\k(K),\k\bigr)=
  \b^{-q,2p}\bigl(\k(K)\bigr)
\end{equation}
(see~(\ref{bbnfr})). Alternatively, $b_{-q,2p}(\zk)$ equals the dimension of
$(-q,2p)$-th bigraded component of the cohomology algebra
$H[\L[u_1,\ldots,u_m]\otimes\k(K),d]$). For the ordinary Betti numbers
$b_k(\zk)$ holds
\begin{equation}\label{ordb}
  b_k(\zk)=\sum_{-q+2p=k}b_{-q,2p}(\zk),\quad k=0,\ldots,m+n.
\end{equation}

Below we describe some basic properties of bigraded Betti
numbers~(\ref{bbn}).

\begin{lemma}\label{bbgen}
  Let $K^{n-1}$ be a simplicial complex with $m=f_0$ vertices and $f_1$
  edges, and let $\zk$ be the corresponding moment-angle complex,
  $\dim\zk=m+n$. Then
  \begin{itemize}
  \item[(a)] $b_{0,0}(\zk)=b_0(\zk)=1$,\quad $b_{0,2p}(\zk)=0$ for $p>0${\rm;}
  \item[(b)] $b_{-q,2p}=0$ for $p>m$ or $q>p${\rm;}
  \item[(c)] $b_1(\zk)=b_2(\zk)=0${\rm;}
  \item[(d)] $b_3(\zk)=b_{-1,4}(\zk)=\binom{f_0}2-f_1${\rm;}
  \item[(e)] $b_{-q,2p}(\zk)=0$ for $q\ge p>0$ or $p-q>n${\rm;}
  \item[(f)] $b_{m+n}(\zk)=b_{-(m-n),2m}(\zk)$.
  \end{itemize}
\end{lemma}
\begin{proof}
We make calculations with the cochain complex
$A^{*}(K)\subset\L[u_1,\ldots,u_m]\otimes\k(K)$. The module $A^*(K)$ has the
basis consisting of monomials $u_Jv_I$ with $I\in K$ and $I\cap J=\emptyset$.
Since $\bideg v_i=(0,2)$, $\bideg u_j=(-1,2)$, the bigraded component
$A^{-q,2p}(K)$ is generated by monomials $u_Jv_I$ with $\#I=p-q$ and $\#J=q$.
In particular, $A^{-q,2p}(K)=0$ if $p>m$ or $q>p$, whence the assertion~(b)
follows. To prove~(a) we observe that $A^{0,0}(K)$ is generated by 1, while
any $v_I\in A^{0,2p}(K)$, $p>0$, is a coboundary, whence
$H^{0,2p}(\zk)=0$, $p>0$.

Now we prove the assertion (e). Any $u_Jv_I\in A^{-q,2p}(K)$ has $I\in K$,
while any simplex of $K$ is at most $(n-1)$-dimensional. It follows that
$A^{-q,2p}(K)=0$ for $p-q>n$. By~(b), $b_{-q,2p}(\zk)=0$ for $q>p$, so it
remains to prove that $b_{-q,2q}(\zk)=0$ for~$q>0$. The module $A^{-q,2q}(K)$
is generated by monomials $u_J$, $\#J=q$. Since $d(u_i)=v_i$, it follows
easily that there are no non-zero cocycles in $A^{-q,2q}(K)$. Hence,
$H^{-q,2q}(\zk)=0$.

The assertion~(c) follows from~(e) and~(\ref{ordb}).

It also follows from (e) that $H^{3}(\zk)=H^{-1,4}(\zk)$.  The basis for
$A^{-1,4}(K)$ consists of monomials $u_jv_i$, $i\ne j$. We have
$d(u_jv_i)=v_iv_j$ and $d(u_iu_j)=u_jv_i-u_iv_j$. It follows that $u_jv_i$ is
a cocycle if and only if $\{i,j\}$ is not a 1-simplex in~$K$; in this case
two cocycles $u_jv_i$ and $u_iv_j$ represent the same cohomology class. The
assertion~(d) now follows easily.

The remaining assertion (f) follows from the fact that the monomial
$u_Iv_J\in A^*(K)$ of maximal total degree $m+n$ necessarily has
$\#I+\#J=m$, $\#J=n$, $\#I=m-n$.
\end{proof}

Lemma~\ref{bbgen} shows that non-zero bigraded Betti numbers $b_{r,2p}(\zk)$,
$r\ne0$ appear only in the ``strip" bounded by the lines $p=m$, $r=-1$,
$p+r=1$ and $p+r=n$ in the second quadrant (see Figure~8~(a)).
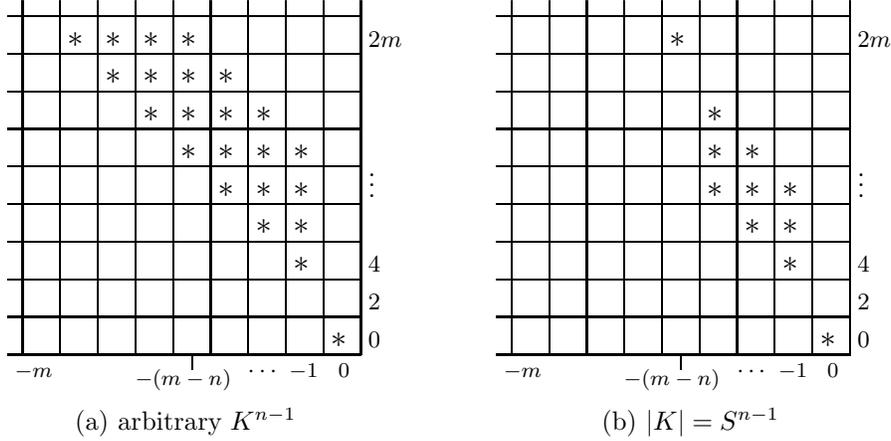
\begin{figure}
  \begin{picture}(115,60)
  \multiput(45,10)(0,5){10}{\line(-1,0){47}}
  \multiput(45,10)(-5,0){10}{\line(0,1){47}}
  \put(46,11){\small 0}
  \put(46,16){\small 2}
  \put(46,21){\small 4}
  \put(46,31){\vdots}
  \put(46,51){\small $2m$}
  \put(42,7){\footnotesize 0}
  \put(35.5,7){\footnotesize $-1$}
  \put(30,7){\small $\cdots$}
  \put(22.5,10){\line(0,-1){2}}
  \put(15,6){\footnotesize $-(m-n)$}
  \put(-1,7){\footnotesize $-m$}
  \put(41,11){\Large $*$}
  \multiput(36,21)(-5,5){7}{\Large $*$}
  \multiput(36,26)(-5,5){6}{\Large $*$}
  \multiput(36,31)(-5,5){5}{\Large $*$}
  \multiput(36,36)(-5,5){4}{\Large $*$}
  \put(7,0){(a)\ arbitrary $K^{n-1}$}
  \multiput(110,10)(0,5){10}{\line(-1,0){47}}
  \multiput(110,10)(-5,0){10}{\line(0,1){47}}
  \put(111,11){\small 0}
  \put(111,16){\small 2}
  \put(111,21){\small 4}
  \put(111,31){\vdots}
  \put(111,51){\small $2m$}
  \put(107,7){\footnotesize 0}
  \put(100.5,7){\footnotesize $-1$}
  \put(95,7){\small $\cdots$}
  \put(87.5,10){\line(0,-1){2}}
  \put(80,6){\footnotesize $-(m-n)$}
  \put(64,7){\footnotesize $-m$}
  \put(106,11){\Large $*$}
  \multiput(101,21)(-5,5){3}{\Large $*$}
  \multiput(101,26)(-5,5){3}{\Large $*$}
  \multiput(101,31)(-5,5){3}{\Large $*$}
  \put(86,51){\Large $*$}
  \put(77,0){(b)\ $|K|=S^{n-1}$}
  \end{picture}
  \caption{Possible locations of non-zero bigraded Betti numbers
  $b_{-q,2p}(\zk)$ (marked by $*$).}
\end{figure}

The homogeneous component $C_{-q,2p}(\zk)$ has the basis consisting of
cellular chains $D_IT_J$ with $I\in K$, $\#I=p-q$, $\#J=q$. It follows that
\begin{equation}\label{dimc}
  \dim C_{-q,2p}(\zk)=f_{p-q-1}\bin{m-p+q}q
\end{equation}
(we assume $\binom ij=0$ if $i<j$ or $j<0$), where
$(f_0,f_1,\ldots,f_{n-1})$ is the $f$-vector of $K^{n-1}$ and $f_{-1}=1$.
The differential $\partial_c$ does not change the second degree, i.e.
$$
  \partial_c:C_{-q,2p}(\zk)\to C_{-q-1,2p}(\zk).
$$
Hence, the chain complex $C_{*,*}(\zk)$ splits as follows:
$$
  [C_{*,*}(\zk),\partial_c]=\bigoplus_{p=0}^m
  [C_{*,2p}(\zk),\partial_c].
$$

\begin{remark}
The similar decomposition holds also for the cellular cochain complex
$[C^{*,*}(\zk),d_c]\cong[A^{*,*}(K),d]$.
\end{remark}

Let us consider the Euler characteristic of the complex
$[C_{*,2p}(\zk),\partial_c]$:
\begin{equation}\label{chip}
  \chi_p(\zk):=\sum_{q=0}^m(-1)^q\dim C_{-q,2p}(\zk)
  =\sum_{q=0}^m(-1)^qb_{-q,2p}(\zk).
\end{equation}
Define the generating polynomial $\chi(\zk;t)$ as
$$
  \chi(\zk;t)=\sum_{p=0}^m\chi_p(\zk)t^{2p}.
$$
The following theorem calculates this polynomial in terms of the $h$-vector
of~$K$.

\begin{theorem}
\label{gpz}
  For any $(n-1)$-dimensional simplicial complex $K$ with $m$ vertices holds
  \begin{equation}
  \label{hchi}
    \chi(\zk;t)=(1-t^2)^{m-n}(h_0+h_1t^2+\cdots+h_nt^{2n}),
  \end{equation}
  where $(h_0,h_1,\ldots,h_n)$ is the $h$-vector of~$K$.
\end{theorem}
\begin{proof}
It follows from (\ref{chip}) and (\ref{dimc}) that
\begin{equation}
\label{chipzk}
  \chi_p(\zk)=\sum_{j=0}^m(-1)^{p-j}f_{j-1}\binom{m-j}{p-j},
\end{equation}
Then
\begin{multline}
\label{chidir}
  \chi(\zk;t)=\sum_{p=0}^m\chi_p(K)t^{2p}=
  \sum_{p=0}^m\sum_{j=0}^mt^{2j}t^{2(p-j)}(-1)^{p-j}f_{j-1}
  \binom{m-j}{p-j}\\
  =\sum_{j=0}^mf_{j-1}t^{2j}(1-t^2)^{m-j}=
  (1-t^2)^m\sum_{j=0}^nf_{j-1}(t^{-2}-1)^{-j}.
\end{multline}
Denote $h(t)=h_0+h_1t+\cdots+h_nt^n$. Then it follows from~(\ref{hvector})
that
$$
  t^nh(t^{-1})=(t-1)^n\sum_{i=0}^nf_{i-1}(t-1)^{-i}.
$$
Substituting $t^{-2}$ for $t$ above, we finally obtain from~(\ref{chidir})
$$
  \frac{\chi(\zk;t)}{(1-t^2)^m}=\frac{t^{-2n}h(t^2)}{(t^{-2}-1)^n}=
  \frac{h(t^2)}{(1-t^2)^n},
$$
which is equivalent to (\ref{hchi}).
\end{proof}

Theorem~\ref{gpz} allows to express the numbers of faces of a simplicial
complex in terms of the bigraded Betti numbers of the corresponding
moment-angle complex~$\zk$.

\begin{corollary}\label{chizk}
  For any simplicial complex $K$ the Euler characteristic of the corresponding
  moment-angle complex $\zk$ is zero.
\end{corollary}
\begin{proof}
We have
$$
  \chi(\zk)=\sum_{p,q=0}^m(-1)^{-q+2p}b_{-q,2p}(\zk)=
  \sum_{p=0}^m\chi_p(\zk)=\chi(\zk;1)
$$
Now the statement follows from (\ref{hchi}).
\end{proof}

\begin{remark}
Another proof of the above corollary follows from the observation
that the diagonal subgroup $S^1\subset T^m$ always acts freely on
$\zk$ (see also section~\ref{part}). Hence, there exists a principal
$S^1$-bundle $\zk\to\zk/S^1$, which implies $\chi(\zk)=0$.
\end{remark}

The torus $T^m=\rho^{-1}(1,\ldots,1)$ is a cellular subcomplex of $\zk$ (see
Lemma~\ref{the incl}). The cellular cochain subcomplex
$C^{*}(T^m)\subset C^*(\zk)\cong A^*(K)$ has the basis consisting of cochains
$(T_I)^*$ and is mapped to the exterior algebra
$\Lambda[u_1,\ldots,u_m]\subset A^*(K)$ under the isomorphism of
Theorem~\ref{cellcom}. It follows that there is an isomorphism of modules
\begin{equation}\label{pair}
  C^{*}(\zk,T^m)\cong A^*(K)/\Lambda[u_1,\ldots,u_m].
\end{equation}
Likewise, we introduce relative bigraded Betti numbers
\begin{equation}\label{rbbn}
  b_{-q,2p}(\zk,T^m)=\dim H^{-q,2p}\bigl[C^{*}(\zk,T^m),d\bigr],
  \quad q,p=0,\ldots,m,
\end{equation}
define the $p$-th relative Euler characteristic $\chi_p(\zk,T^m)$ as
\begin{equation}\label{rchip}
  \chi_p(\zk,T^m)=\sum_{q=0}^m(-1)^q\dim C^{-q,2p}(\zk,T^m)
  =\sum_{q=0}^m(-1)^qb_{-q,2p}(\zk,T^m),
\end{equation}
and define the generating polynomial $\chi(\zk,T^m;t)$ as
$$
  \chi(\zk,T^m;t)=\sum_{p=0}^m\chi_p(\zk,T^m)t^{2p}.
$$

\begin{theorem}
\label{rgp}
  For any $(n-1)$-dimensional simplicial complex $K$ with $m$ vertices holds
  \begin{equation}
  \label{rhchi}
    \chi(\zk,T^m;t)=(1-t^2)^{m-n}(h_0+h_1t^2+\cdots+h_nt^{2n})-(1-t^2)^m.
  \end{equation}
\end{theorem}
\begin{proof}
Since $C^{*}(T^m)=\Lambda[u_1,\ldots,u_m]$ and $\bideg u_i=(-1,2)$,
we have
$$
  \dim C^{-q}(T^m)=\dim C^{-q,2q}(T^m)=\textstyle\binom mq.
$$
Combining (\ref{pair}), (\ref{chip}) and (\ref{rchip}) we get
$$
  \chi_p(\zk,T^m)=\chi_p(\zk)-(-1)^p\dim C^{-p,2p}(T^m).
$$
Hence,
\begin{align*}
  \chi(\zk,T^m;t)&=\chi(\zk;t)-\sum_{p=0}^m(-1)^p{\textstyle\binom
  mpt^{2p}}\\
  &=(1-t^2)^{m-n}(h_0+h_1t^2+\cdots+h_nt^{2n})-(1-t^2)^m,
\end{align*}
where we used (\ref{hchi}).
\end{proof}

We will use Theorem~\ref{rgp} in section~\ref{coh3}.

\begin{theorem}\label{reduc}
Suppose that $K^{n-1}$ is Cohen--Macaulay, and let $\mathcal J$ be the ideal
in $\k(K)$ generated by a degree-two regular sequence of length~$n$. Then the
following isomorphism of algebras holds:
$$
  H^*(\zk)\cong\Tor_{\k[v_1,\ldots,v_m]/\mathcal J}
  \bigl(\k(K)/\mathcal J,\k\bigr).
$$
\end{theorem}
\begin{proof}
This follows from Theorem~\ref{cohom1} and Lemma~\ref{tortor}.
\end{proof}
Note that the $\k$-algebra $\k(K)/\mathcal J$ is finite-dimensional. This
fact sometimes makes Theorem~\ref{reduc} more convenient for
calculations (in the Cohen--Macaulay case) than general
Theorem~\ref{cohom1}.

\subsection{Cohomology of moment-angle complex $\zk$: the case of
spherical~$K$}
\label{coh2}
If $K$ is a simplicial sphere, then the complex $\zk$ is a manifold
(Lemma~\ref{maman}). This imposes additional conditions on the cohomology
of~$\zk$; the corresponding results are described in this section together
with some interesting interpretations of combinatorial problems reviewed
in chapter~1.

\begin{theorem}\label{fc}
  Let $K$ be an $(n-1)$-dimensional simplicial sphere, and $\zk$ the
  corresponding moment-angle manifold, $\dim\zk=m+n$. Then the fundamental
  cohomology class of $\zk$ is represented by any monomial $\pm v_Iu_J\in
  A^*(K)$ of bidegree $(-(m-n),2m)$ such that $I$ is an $(n-1)$-simplex of
  $K$ and $I\cap J=\emptyset$. The sign depends on the orientation of~$\zk$.
\end{theorem}
\begin{proof}
It follows from Lemma~\ref{bbgen}~(f) that $H^{m+n}(\zk)=H^{-(m-n),2m}(\zk)$.
By definition, the module $A^{-(m-n),2m}(K)$ is spanned by monomials $v_Iu_J$
such that $I\in K^{n-1}$, $\#I=n$, $J=[m]\setminus I$. Any such monomial is a
cocycle. Suppose that $I,I'$ are two $(n-1)$-simplices of $K^{n-1}$ that
share a common $(n-2)$-face. We claim that the corresponding cocycles
$v_Iu_J$, $v_{I'}u_{J'}$, where $J=[m]\setminus I$, $J'=[m]\setminus I'$,
represent the same cohomology class (up to a sign). Indeed, let
$v_Iu_J=v_{i_1}\cdots v_{i_n}u_{j_1}\cdots u_{j_{m-n}}$, $v_{I'}u_{J'}=
v_{i_1}\cdots v_{i_{n-1}}v_{j_1}u_{i_n}u_{j_2}\cdots u_{j_{m-n}}$. Since any
$(n-2)$-face of $K$ is contained in exactly two $(n-1)$-faces, the identity
\begin{multline*}
  d(v_{i_1}\cdots v_{i_{n-1}}u_{i_n}u_{j_1}u_{j_2}\cdots u_{j_{m-n}})\\
  =v_{i_1}\cdots v_{i_n}u_{j_1}\cdots u_{j_{m-n}}-
  v_{i_1}\cdots v_{i_{n-1}}v_{j_1}u_{i_n}u_{j_2}\cdots u_{j_{m-n}}
\end{multline*}
holds in $A^*(K)\subset\L[u_1,\ldots,u_m]\otimes\k(K)$. Hence,
$[v_Iu_J]=[v_{I'}u_{J'}]$ (as cohomology classes). Since $K^{n-1}$ is a
simplicial sphere, any two $(n-1)$-simplices can be connected by a chain of
simplices such that any two successive simplices share a common $(n-2)$-face.
Thus, all monomials $v_Iu_J$ in $A^{-(m-n),2m}(K)$ represent the same
cohomology class (up to a sign). This class is the generator of
$H^{m+n}(\zk)$, i.e. the fundamental cohomology class of~$\zk$.
\end{proof}

\begin{remark}
In the proof of the above theorem we have used two combinatorial properties
of~$K^{n-1}$. The first one is that any $(n-2)$-face is contained in exactly
two $(n-1)$-faces, and the second is that any two $(n-1)$-simplices can be
connected by a chain of simplices such that any two successive simplices
share a common $(n-2)$-face. Both properties hold for any simplicial
manifold. Hence, for any simplicial manifold $K^{n-1}$ we have
$b_{m+n}(\zk)=b_{-(m-n),2m}(\zk)=1$, and the generator of $H^{m+n}(\zk)$ can
be chosen as described in Theorem~\ref{fc}.
\end{remark}

\begin{corollary}\label{bpd}
  The Poincar\'e duality for the moment angle manifold $\zk$ defined by a
  simplicial sphere $K^{n-1}$ regards the bigraded structure in the
  (co)homology, i.e.
  $$
    H^{-q,2p}(\zk)\cong H_{-(m-n)+q,2(m-p)}(\zk).
  $$
  In particular,
  \begin{equation}\label{bpdbn}
    b_{-q,2p}(\zk)=b_{-(m-n)+q,2(m-p)}(\zk).\qquad\square
  \end{equation}
\end{corollary}

\begin{corollary}\label{bbss}
  Let $K^{n-1}$ be an $(n-1)$-dimensional simplicial sphere, and $\zk$
  the corresponding moment-angle complex, $\dim\zk=m+n$. Then
  \begin{itemize}
  \item[(a)] $b_{-q,2p}(\zk)=0$ for $q\ge m-n$, with only exception
   $b_{-(m-n),2m}=1$;
  \item[(b)] $b_{-q,2p}(\zk)=0$ for $p-q\ge n$, with only exception
   $b_{-(m-n),2m}=1$.
  \end{itemize}
\end{corollary}

It follows that if $K^{n-1}$ is a simplicial sphere, then non-zero
bigraded Betti numbers $b_{r,2p}(\zk)$, $r\ne0$, $r\ne m-n$, appear only in
the ``strip" bounded by the lines $r=-(m-n-1)$, $r=-1$, $p+r=1$ and $p+r=n-1$
in the second quadrant (see Figure~8~(b)). Compare this with Figure~8~(a)
corresponding to the case of general~$K$.

\begin{example}
Let $K=\partial\D^{m-1}$. Then $\k(K)=\k[v_1,\ldots,v_m]/(v_1\cdots v_m)$
(Example~\ref{simbsim}). It is easy to see that the cohomology groups
$H[\k(K)\otimes\Lambda[u_1,\ldots,u_m],d]$ (see Theorem~\ref{cohom2}) are
generated by the classes $1$ and $[v_1v_2\cdots v_{m-1}u_m]$. We have $\deg
(v_1v_2\cdots v_{m-1}u_m)=2m-1$, and Theorem~\ref{fc} shows that
$v_1v_2\cdots v_{m-1}u_m$ represents the fundamental cohomology class of
$\zk\cong S^{2m-1}$ (see Example~\ref{zpsphere}).
\end{example}

\begin{example}\label{mgon}
Let $K$ be the boundary complex of an $m$-gon $P^2$ with $m\ge4$. We have
$\k(K)=\k[v_1,\ldots,v_m]/\mathcal I_P$, where $\mathcal I_P$ is generated by
the monomials $v_iv_j$, $i-j\ne0,1\mod m$. The complex $\zk=\zp$ is a smooth
manifold of dimension~$m+2$. The Betti numbers and the cohomology rings of
these manifolds were calculated in~\cite{BP2}. Namely,
$$
  \dim H^k(\zp)=\left\{
  \begin{array}{l}
    1\quad\text{for }k=0,m+2;\\[1mm]
    0\quad\text{for }k=1,2,m,m+1;\\[1mm]
    (m-2)\binom{m-2}{k-2}-\binom{m-2}{k-1}-\binom{m-2}{k-3}
    \quad\text{for }3\le k\le m-1.
  \end{array}
  \right.
$$
For example, in the case $m=5$ the group $H^3(\zp)$ has 5 generators
represented by the cocycles
$v_iu_{i+2}\in\k(K)\otimes\Lambda[u_1,\ldots,u_5]$, $i=1,\ldots,5$, while the
group $H^4(\zp)$ has 5 generators represented by the cocycles
$v_ju_{j+2}u_{j+3}$, $j=1,\ldots,5$. As it follows from Theorem~\ref{fc}, the
product of cocycles $v_iu_{i+2}$ and $v_ju_{j+2}u_{j+3}$ represents a
non-zero cohomology class in $H^7(\zp)$ if and only if all indices
$i,i+2,j,j+2,j+3$ are different. Thus, for each of the 5 cohomology classes
$[v_iu_{i+2}]$ there is the unique (Poincar\'e dual) cohomology class
$[v_ju_{j+2}u_{j+3}]$ such that the product
$[v_iu_{i+2}]\cdot[v_ju_{j+2}u_{j+3}]$ is non-zero.
\end{example}

It follows from (\ref{chip}) and (\ref{bpdbn}) that for any simplicial sphere
$K$ holds
$$
  \chi_p(\zk)=(-1)^{m-n}\chi_{m-p}(\zk).
$$
From this and (\ref{hchi}) we get
\begin{multline*}
  \frac{h_0+h_1t^2+\cdots+h_nt^{2n}}{(1-t^2)^n}=
  (-1)^{m-n}\frac{\chi_m+\chi_{m-1}t^2+\cdots+\chi_0t^{2m}}{(1-t^2)^m}\\=
  (-1)^n\frac{\chi_0+\chi_1t^{-2}+\cdots+\chi_mt^{-2m}}{(1-t^{-2})^m}
  =(-1)^n\frac{h_0+h_1t^{-2}+\cdots+h_nt^{-2n}}{(1-t^{-2})^n}\\=
  \frac{h_0t^{2n}+h_1t^{2(n-1)}+\cdots+h_n}{(1-t^2)^n}.
\end{multline*}
Hence, $h_i=h_{n-i}$. Thus, the Dehn--Sommerville equations are
a corollary of the bigraded Poincar\'e duality~(\ref{bpdbn}).

The identity~(\ref{hchi}) also allows to interpret different inequalities for
the $f$-vectors of simplicial spheres (respectively, simplicial manifolds) in
terms of topological invariants (bigraded Betti numbers) of the corresponding
moment-angle manifolds (respectively, complexes)~$\zk$.

\begin{example}
It follows from Lemma~\ref{bbgen} that for any $K$ we have
\begin{align*}
  &\chi_0(\zk)=1,&&\chi_1(\zk)=0,\\
  &\chi_2(\zk)=-b_{-1,4}(\zk)=-b_3(\zk),&&
  \chi_3(\zk)=b_{-2,6}(\zk)-b_{-1,6}(\zk)
\end{align*}
(note that $b_4(\zk)=b_{-2,6}(\zk)$, while
$b_5(\zk)=b_{-1,6}(\zk)+b_{-3,8}(\zk)$). Now, identity~(\ref{hchi}) shows
that
\begin{align*}
  &h_0=1,\\
  &h_1=m-n,\\
  &h_2={\textstyle\binom{m-n+1}2}-b_3(\zk),\\
  &h_3={\textstyle\binom{m-n+2}3}-(m-n)b_{-1,4}(\zk)+
   b_{-2,6}(\zk)-b_{-1,6}(\zk).
\end{align*}
It follows that the inequality $h_1\le h_2$ ($n\ge4$) from the Generalized
Lower Bound hypothesis~(\ref{glbt}) for simplicial spheres is equivalent to
the following:
\begin{equation}
\label{h12}
  b_3(\zk)\le\textstyle\binom{m-n}2.
\end{equation}
The next inequality $h_2\le h_3$ ($n\ge6$) from~(\ref{glbt}) is equivalent to
the following inequality for the bigraded Betti numbers of~$\zk$:
\begin{equation}
\label{h23}
  {\textstyle\binom{m-n+1}3}-(m-n-1)b_{-1,4}(\zk)+
  b_{-2,6}(\zk)-b_{-1,6}(\zk)\ge0.
\end{equation}
\end{example}

We see that the combinatorial Generalized Lower Bound inequalities are
interpreted as ``topological" inequalities for the (bigraded) Betti numbers
of a certain manifold. So, one can try to use topological methods (such as
the equivariant topology or Morse theory) to prove inequalities
like~(\ref{h12}) or~(\ref{h23}).  Such topological approach to problems like
$g$-conjecture or Generalized Lower Bound has the advantage of being
independent on whether the simplicial sphere $K$ is polytopal or not. Indeed,
all known proofs of the necessity of $g$-theorem for simplicial polytopes
(including the original one by Stanley given in section~\ref{tori},
McMullen's proof~\cite{McM2}, and the recent proof by Timorin~\cite{Ti})
follow the same scheme. Namely, the numbers $h_i$, $i=1,\ldots,n$, are
interpreted as the dimensions of graded components $A^i$ of a certain graded
algebra $A$ satisfying the Hard Lefschetz Theorem. The latter means that
there is an element $\omega\in A^1$ such that the multiplication by $\omega$
defines a {\it monomorphism\/} $A^i\to A^{i+1}$ for $i<\bigl[\frac n2\bigr]$.
This implies $h_i\le h_{i+1}$ for $i<\bigl[\frac n2\bigr]$ (see
section~\ref{tori}). However, such an element $\omega$ is lacking for
non-polytopal~$K$, which means that a new technique has to be developed for
proving the $g$-conjecture for simplicial spheres.

As it was mentioned in section~\ref{sim2}, simplicial spheres are Gorenstein*
complexes. Using theorems~\ref{gorencom},~\ref{tordual} and our
Theorem~\ref{cohom1} we obtain the following solution of the analogue of
Problem~\ref{zkmanprob}.

\begin{theorem}
The complex $\zk$ is a Poincar\'e duality complex (over~$\k$) if and only if
for any simplex $I\in K$ (including $I=\emptyset$) the subcomplex $\link I$
has the homology of a sphere of dimension $\dim\,(\link I)$.
\end{theorem}

\subsection{Partial quotients of manifold $\zp$}
\label{part}
Here we return to the case of polytopal $K$ (i.e. $K=K_P$) and study
quotients of $\zp$ by freely acting subgroups $H\subset T^m$.

For any combinatorial simple polytope $P^n$ denote by $s(P^n)$ the maximal
dimension of subgroups $H\subset T^m$ that act freely on~$\zp$. The number
$s(P^n)$ is obviously a combinatorial invariant of~$P^n$.

\begin{problem}[V.\,M. Buchstaber]\label{sp}
  Express $s(P^n)$ via known combinatorial invariants of~$P^n$.
\end{problem}

\begin{proposition}
If $P^n$ has $m$ facets, then $s(P^n)\le m-n$.
\end{proposition}
\begin{proof}
Every subgroup of $T^m$ of dimension $>m-n$ intersects non-trivially with any
$n$-dimensional isotropy subgroup, and therefore can not act freely on~$\zp$.
\end{proof}

\begin{proposition}\label{diags}
The diagonal circle subgroup $S_d:=\{(e^{2\pi i\f},\ldots,e^{2\pi i\f})\in
T^m\}$, $\f\in\R$, acts freely on any~$\zp$. Thus, $s(P^n)\ge1$.
\end{proposition}
\begin{proof}
Since any isotropy subgroup for $\zp$ is coordinate (see
Definition~\ref{zp}), it intersects with $S_d$ only at the unit.
\end{proof}

Another lower bound for the number $s(P^n)$ was proposed in~\cite{Iz2}. Let
$\F=\{F_1,\ldots,F_m\}$ be the set of facets of~$P^n$. The surjective map
$\varrho:\F\to[k]$ (where $[k]=\{1,\ldots,k\}$) is called the {\it regular
$k$-paint colouring of facets of\/}~$P^n$ if $\varrho(F_i)\ne\varrho(F_j)$
whenever $F_i\cap F_j\ne\emptyset$. The {\it chromatic number\/}
$\gamma(P^n)$ of a polytope $P^n$ is the minimal $k$ for which there exists a
regular $k$-paint colouring of facets of~$P^n$.

\begin{example}
Suppose $P^n$ is a 2-neighbourly simple polytope with
$m$ facets. Then $\gamma(P^n)=m$.
\end{example}

\begin{proposition}[{\cite{Iz2}}]
The following inequality holds:
$$
  s(P^n)\ge m-\gamma(P^n).
$$
\end{proposition}
\begin{proof}
The map $\varrho:\F\to[k]$ defines the epimorphism of tori
$\tilde\varrho:T^m\to T^k$. It is easy to see that if $\varrho$ is a regular
colouring, then $\Ker\tilde\varrho\cong T^{m-k}$ acts freely on~$\zp$.
\end{proof}

Let $H\subset T^m$ be a subgroup of dimension $r\le m-n$. Choosing a basis
in~$H$, we can write it in the form
\begin{equation}\label{h}
  H=\bigl\{
  (e^{2\pi i(s_{11}\f_1+\dots+s_{1r}\f_r)},\ldots, e^{2\pi
     i(s_{m1}\f_1+\dots+s_{mr}\f_r)})\in T^m \bigr\},
\end{equation}
where $\f_i\in\R$, $i=1,\ldots,r$. The $m\times r$ integer matrix $S=(s_{ij})$
defines a monomorphism $\Z^r\to\Z^m$ whose image is a direct summand
in~$\Z^m$. For any subset $\{i_1,\ldots,i_n\}\subset[m]$ denote by $S_{\hat
i_1,\ldots,\hat i_n}$ the $(m-n)\times r$ submatrix of $S$ obtained by
deleting the rows $i_1,\ldots,i_n$. Recall that any vertex $v\in P^n$ is the
intersection of $n$ facets (see~(\ref{vert})). The following criterion of
freeness for the action of $H$ on $\zp$ holds.

\begin{lemma}
\label{free}
  The subgroup~{\rm(\ref{h})} acts freely on $\zp$ if and only if for any
  vertex $v=F_{i_1}\cap\ldots\cap F_{i_n}$ of $P^n$ the $(m-n)\times
  r$-submatrix $S_{\hat i_1,\ldots,\hat i_n}$ defines the monomorphism
  $\Z^r\hookrightarrow\Z^{m-n}$ to a direct summand.
\end{lemma}
\begin{proof}
It follows from Definition~\ref{zp} that the orbits of the $T^m$-action on
$\zp$ corresponding to the vertices of $P^n$ have maximal (rank~$n$)
isotropy subgroups. The isotropy subgroup corresponding to the vertex
$v=F_{i_1}\cap\ldots\cap F_{i_n}$ is the coordinate subtorus
$T^n_{i_1,\ldots,i_n}\subset T^m$. The subgroup~(\ref{h}) acts freely on $\zp$
if and only if it intersects each isotropy subgroup only at the unit. This is
equivalent to the condition that the map $H\times T^n_{i_1,\ldots,i_n}\to
T^m$ is injective for any $v=F_{i_1}\cap\ldots\cap F_{i_n}$. The latter map is
injective whenever the image of the corresponding map
$\Z^{r+n}\hookrightarrow\Z^m$ is a direct summand of~$\Z^m$. The matrix of
the map $\Z^{r+n}\hookrightarrow\Z^m$ is obtained by adding $n$ columns
$(0,\ldots,0,1,0,\ldots,0)^t$ (1 stands at the place $i_j$, $j=1,\ldots,n$)
to~$S$. This matrix defines a direct summand of $\Z^m$ exactly when the same
is true for each $S_{\hat i_1,\ldots,\hat i_n}$.
\end{proof}

In particular, for subgroups of rank $m-n$ we get

\begin{corollary}
\label{maxfree}
  The subgroup~{\rm(\ref{h})} of rank $r=m-n$ acts freely on $\zp$ if and
  only if for any vertex $v=F_{i_1}\cap\ldots\cap F_{i_n}$ of $P^n$ the minor
  $S_{\hat i_1,\ldots,\hat i_n}$ satisfies ${\det S_{\hat i_1\ldots
  \hat i_n}=\pm1}$.
\end{corollary}

\begin{proposition}\label{chars}
A simple polytope $P^n$ admits a characteristic map
if and only if $s(P^n)=m-n$.
\end{proposition}
\begin{proof}
Proposition~\ref{zpbun} shows that if $P^n$ admits a characteristic
map~$\ell$, then the $(m-n)$-dimensional subgroup $H(\ell)$ acts freely
on~$\zp$, whence $s(P^n)=m-n$. Now suppose $s(P^n)=m-n$, i.e. there exists a
subgroup~(\ref{h}) of rank $r=m-n$ that acts freely on~$\zp$. The
corresponding $m\times(m-n)$-matrix $S$ defines a monomorphism
$\Z^{m-n}\to\Z^m$ whose image is a direct summand. It follows that there is
an $n\times m$-matrix $\L$ such that the sequence
$$
\begin{CD}
  0 @>>> \Z^{m-n} @>S>> \Z^m @>\L>> \Z^n @>>> 0
\end{CD}
$$
is exact. Since $S$ satisfies the condition of Corollary~\ref{maxfree}, the
matrix $\L$ satisfies the condition~(\ref{L}), thus defining a characteristic
map for~$P^n$.
\end{proof}

For any subgroup~(\ref{h}) of dimension $r=m-n$ define the following
linear forms in $\k[v_1,\ldots,v_m]$:
\begin{equation}\label{w}
w_i=s_{1i}v_1+\dots+s_{mi}v_m,\quad i=1,\ldots,m-n.
\end{equation}

Suppose $M^{2n}$ is a quasitoric manifold over $P^n$ with characteristic
map~$\ell$. Write the subgroup $H(\ell)$ in the form~(\ref{h}); this
defines elements~(\ref{w}). Under these assumptions the following statement
holds.
\begin{lemma}
\label{zpqt}
  The following isomorphism of algebras holds:
  $$
    H^{*}(\zp)\cong\Tor_{\k[w_1,\ldots,w_{m-n}]}\bigl(H^{*}(M^{2n}),\k\bigr),
  $$
  where the $\k[w_1,\ldots,w_{m-n}]$-module structure on
  $H^{*}(M^{2n})=\k[v_1,\ldots,v_m]/\mathcal I_P+\mathcal J_\ell$ is defined
  by~{\rm(\ref{w})}.
\end{lemma}
\begin{proof}
Theorem~\ref{reduc} shows that
$$
  H^*(\zk)\cong\Tor_{\k[v_1,\ldots,v_m]/\mathcal J_\ell}
  \bigl(\k(K)/\mathcal J_\ell,\k\bigr).
$$
The quotient $\k[v_1,\ldots,v_m]/\mathcal J_\ell$ is identified with
$\k[w_1,\ldots,w_{m-n}]$.
\end{proof}

\begin{theorem}
\label{degene3}
  The Leray--Serre spectral sequence of the $T^{m-n}$-bundle
  $\zp\to M^{2n}$ collapses at the $E_3$ term. Furthermore, the
  following isomorphism of algebras holds:
  \begin{gather*}
    H^{*}(\zp)\cong H\bigl[\Lambda[u_1,\ldots,u_{m-n}]\otimes
    (\k(P)/\mathcal J_\ell),d\bigr],\\
    \bideg v_i=(0,2),\quad\bideg u_i=(-1,2);\\
    d(u_i)=w_i,\quad d(v_i)=0.
  \end{gather*}
\end{theorem}
\begin{proof}
Since $H^*(T^{m-n})=\L[u_1,\ldots,u_{m-n}]$ and
$H^*(M^{2n})=\k(P)/\mathcal J_\ell$, we have
$$
E_3\cong H\bigl[(\k(P)/\mathcal J_\ell)\otimes\L[u_1,\ldots,u_{m-n}],d\bigr].
$$
By Lemma~\ref{koscom},
$$
  H\bigl[(\k(P)/\mathcal J_\ell)\otimes\L[u_1,\ldots,u_{m-n}],d\bigr]\cong
  \Tor_{\k[w_1,\ldots,w_{m-n}]}\bigl(H^{*}(M^{2n}),\k\bigr).
$$
Combining the above two identities with Lemma~\ref{zpqt} we get
$E_3=H^*(\zp)$, which concludes the proof.
\end{proof}

Now we are going to describe the cohomology of the quotient $\zp/H$ for
arbitrary freely acting subgroup~$H$. First, we write $H$ in the
form~(\ref{h}) and choose an $(m-r)\times m$-matrix $T=(t_{ij})$ of rank
$(m-r)$ satisfying $T\cdot S=0$. This is done in the same way as in the proof
of Proposition~\ref{chars} (in particular, $T$ is the characteristic matrix
for the quasitoric manifold $\zp/H$ in the case $r=m-n$).

\begin{theorem}
\label{quot}
  The following isomorphism of algebras holds:
  $$
    H^{*}(\zp/H)\cong\Tor_{\k[t_1,\ldots,t_{m-r}]}
    \bigl(\k(P),\k\bigr),
  $$
  where the $\k[t_1,\ldots,t_{m-r}]$-module structure on
  $\k(P)=\k[v_1,\ldots,v_m]/\mathcal I_P$ is given by the map
  $$
  \begin{array}{rcl}
    k[t_1,\ldots,t_{m-r}]&\to&k[v_1,\ldots,v_m]\\[1mm]
    t_i&\to&t_{i1}v_1+\dots+t_{im}v_m.
  \end{array}
  $$
\end{theorem}

\begin{remark}
  Theorem~\ref{quot} reduces to Theorem~\ref{cohom1} in the case $r=0$ and
  to Example~\ref{qtem} in the case $r=m-n$.
\end{remark}

\begin{proof}[Proof of Theorem~{\rm\ref{quot}}]
The inclusion of the subgroup $T^r\cong H\hookrightarrow T^m$ defines the
map of classifying spaces $h:BT^r\to BT^m$.
Let us consider the commutative square
$$
  \begin{CD}
    E @>>> B_TP\\
    @VVV @VVpV\\
    BT^r @>h>> BT^m,
  \end{CD}
$$
where the left vertical arrow is the pullback along~$h$. It can be easily
seen that $E$ is homotopy equivalent to the quotient~$\zp/H$.
The Eilenberg--Moore spectral sequence of the above square converges to the
cohomology of $\zp/H$ and has
$$
  E_2=\Tor_{\k[v_1,\ldots,v_m]}\bigl(\k(P),\k[w_1,\ldots,w_r]\bigr),
$$
where the $k[v_1,\ldots,v_m]$-module structure on $k[w_1,\ldots,w_r]$ is
defined by the matrix~$S$, i.e. by the map $v_i\to
s_{i1}w_1+\ldots+s_{ir}w_r$. It can be shown in the same way as in the proof
of Theorem~\ref{cohom1} (using cellular decompositions) that the spectral
sequence collapses at the $E_2$ term and the following isomorphism of
algebras holds:
\begin{equation} \label{Ytor1}
  H^{*}(\zp/H)=\Tor_{\k[v_1,\ldots,v_m]}
  \bigl(\k(P),\k[w_1,\ldots,w_r]\bigr).
\end{equation}
Now put $\Lambda=\k[v_1,\ldots,v_m]$, $\Gamma=\k[t_1,\ldots,t_{m-r}]$,
$A=\k[w_1,\ldots,w_r]$, and $C=\k(P)$ in Theorem~\ref{change}. Since
$\Lambda$ here is a free $\Gamma$-module and
$\Omega=\Lambda/\Gamma\cong\k[w_1,\ldots,w_r]$, the spectral sequence
$\{\widetilde{E}_s,\widetilde{d}_s\}$ arises. Its $E_2$ term is
$$
  \widetilde{E}_2=\Tor_{\k[w_1,\ldots,w_r]}
  \Bigl(\k[w_1,\ldots,w_r],
  \Tor_{\k[t_1,\ldots,t_{m-r}]}\bigl(\k(P),\k\bigr)\Bigr),
$$
and it converges to
$\Tor_{\k[v_1,\ldots,v_m]}(\k(P),\k[w_1,\ldots,w_r])$.
Since $\k[w_1,\ldots,w_r]$ is a free $\k[w_1,\ldots,w_r]$-module, we have
$$
  \widetilde{E}_2^{p,q}=0\;\mbox{ for }p\ne0,\quad
  \widetilde{E}_2^{0,*}=\Tor_{\k[t_1,\ldots,t_{m-r}]}\bigl(\k(P),\k\bigr).
$$
Thus, the spectral sequence collapses at the $E_2$ term, and the following
isomorphism of algebras holds:
$$
  \Tor_{\k[v_1,\ldots,v_m]}\bigl(\k(P),\k[w_1,\ldots,w_r]\bigr)\cong
  \Tor_{\k[t_1,\ldots,t_{m-r}]}\bigl(\k(P),\k\bigr),
$$
which together with~(\ref{Ytor1}) proves the theorem.
\end{proof}

\begin{corollary}
$H^{*}(\zp/H)\cong H\bigl[\Lambda[u_1,\ldots,u_{m-r}]\otimes\k,d\bigr]$,
where $du_i=(t_{i1}v_1+\ldots+t_{im}v_m)$, $dv_i=0$,
$\bideg v_i=(0,2)$, $\bideg u_i=(-1,2)$.
\end{corollary}

\begin{example}
Let $H=S_d$ (see Proposition~\ref{diags}).  In this situation the matrix $S$
is the column of $m$ units.  By Theorem~\ref{quot},
\begin{equation}
\label{y1}
  H^{*}(\zp/S_d)\cong\Tor_{\k[t_1,\ldots,t_{m-1}]}
  \bigl(\k(P),\k\bigr),
\end{equation}
where the $\k[t_1,\ldots,t_{m-1}]$-module structure
on $\k(P)=\k[v_1,\ldots,v_m]/I$ is defined by
$$
    t_i\longrightarrow v_i-v_m,\quad i=1,\ldots,m-1.
$$
Suppose that the $S^1$-bundle $\zp\to\zp/S_d$ is classified by the map
$c:\zp/S_d\to BT^1\cong\C P^\infty$. Since $H^{*}(\C P^{\infty})=\k[w]$, the
element $c^{*}(w)\in H^2(\zp/S_d)$ is defined.

\begin{lemma}
\label{neib}
  $P^n$ is $q$-neighbourly if and only if $(c^{*}(w))^q\ne0$.
\end{lemma}
\begin{proof}
The map $c^{*}$ takes the cohomology ring $H^*(BT^1)\cong\k[w]$ to the
subring $\k(P)\otimes_{\k[t_1,\ldots,t_{m-1}]}\k=
\Tor^0_{\k[t_1,\ldots,t_{m-1}]}(\k(P),\k)$ of $H^*(\zp/H)$. This subring is
isomorphic to the quotient $\k(P)/(v_1=\dots=v_m)$. Now the assertion
follows from the fact that a polytope $P^n$ is $q$-neighbourly if and only if
the ideal $\mathcal I_P$ does not contain monomials of degree~$<q+1$.
\end{proof}
\end{example}

\subsection{Bigraded Poincar\'e duality and analogues of Dehn--Sommerville
equations for simplicial manifolds}\label{coh3}
Here we assume that $K^{n-1}$ is a simplicial manifold. In this case the
moment-angle complex $\zk$ is not a manifold, however, its singularities can
be easily treated. Indeed, the cubical complex $\cc(K)$
(Construction~\ref{cck}) is homeomorphic to $|\cone(K)|$, and the vertex of
the cone is the point $p=(1,\ldots,1)\in\cc(K)\subset I^m$.  Let
$U_\varepsilon(p)\subset\cc(K)$ be a small neighbourhood of $p$ in~$\cc(K)$.
Then the closure of $U_\varepsilon(p)$ is also homeomorphic to $|\cone(K)|$.
It follows from the definition of $\zk$ (see~(\ref{zkwk})) that
$U_\varepsilon(T^m):=\rho^{-1}(U_\varepsilon(p))\subset\zk$ is a small
invariant neighbourhood of the torus $T^m=\rho^{-1}(p)$ in~$\zk$. For
small $\varepsilon$ the closure of $U_\varepsilon(T^m)$ is
homeomorphic to $|\cone(K)|\times T^m$. Removing $U_\varepsilon(T^m)$ from
$\zk$ we obtain a manifold with boundary, which we denote~$W_K$. Thus, we
have
$$
  W_K=\zk\setminus U_\varepsilon(T^m),\quad
  \partial W_K=|K|\times T^m.
$$
Note that since the neighbourhood $U_\varepsilon(T^m)$ is $T^m$-stable,
the torus $T^m$ acts on~$W_K$.

\begin{theorem}
\label{homotwk}
  The manifold with boundary $W_K$ is equivariantly homotopy equivalent to
  the moment-angle complex $\wk$ (see~{\rm(\ref{zkwk})}). There is a
  canonical relative homeomorphism of pairs $(W_K,\partial W_K)\to(\zk,T^m)$.
\end{theorem}
\begin{proof}
To prove the first assertion we construct homotopy equivalence
$\cc(K)\setminus U_\varepsilon(p)\to\cub(K)$ as it is shown on Figure~9. This
map is covered by an equivariant homotopy equivalence $W_K=\zk\setminus
U_\varepsilon(T^m)\to\wk$. The second assertion follows easily from the
definition of~$W_K$.
\begin{figure}
  \begin{picture}(120,45)
  \put(70,5){\circle*{2}}
  \put(80,15){\circle*{2}}
  \put(45,30){\circle*{2}}
  \put(55,40){\circle*{2}}
  \put(45,5){\circle*{2}}
  \put(80,40){\circle*{2}}
  \multiput(45,29.3)(0,0.1){16}{\line(1,1){10}}
  \multiput(70,4.3)(0,0.1){16}{\line(1,1){10}}
  \put(73,33){\vector(1,1){4}}
  \put(73,33){\line(1,1){7}}
  \put(45,30){\line(1,0){20}}
  \put(70,5){\line(0,1){20}}
  \put(65,30){\vector(-1,0){15}}
  \put(70,25){\vector(0,-1){15}}
  \put(65.5,28){\line(-3,-2){20.5}}
  \put(68,25.6){\line(-2,-3){13.5}}
  \put(65.5,28){\vector(-3,-2){13}}
  \put(68,25.6){\vector(-2,-3){9}}
  \put(67,31.5){\line(-4,1){16}}
  \put(70,32.5){\line(-2,3){5}}
  \put(67,31.5){\vector(-4,1){12.5}}
  \put(70,32.5){\vector(-2,3){3}}
  \put(71.5,27){\line(1,-4){4}}
  \put(72.5,30){\line(1,-1){7}}
  \put(71.5,27){\vector(1,-4){3}}
  \put(72.5,30){\vector(1,-1){5}}
  \put(70,30){\oval(10,10)[lb]}
  \qbezier(65,30)(68,33)(73,33)
  \qbezier(70,25)(73,28)(73,33)
  \put(77,32.5){\line(4,1){17}}
  \put(96,36){$\cc(K)\setminus U_\varepsilon(p)$}
  \put(45,20){\line(-2,3){7}}
  \put(50,35){\line(-4,-1){10}}
  \put(24,32.5){$\cub(K)$}
  \linethickness{1mm}
  \put(45,5){\line(1,0){25}}
  \put(55,40){\line(1,0){25}}
  \put(45,5){\line(0,1){25}}
  \put(80,15){\line(0,1){25}}
  \end{picture}
  \caption{Homotopy equivalence
  $\cc(K)\setminus U_\varepsilon(p)\to\cub(K)$.}
  \end{figure}
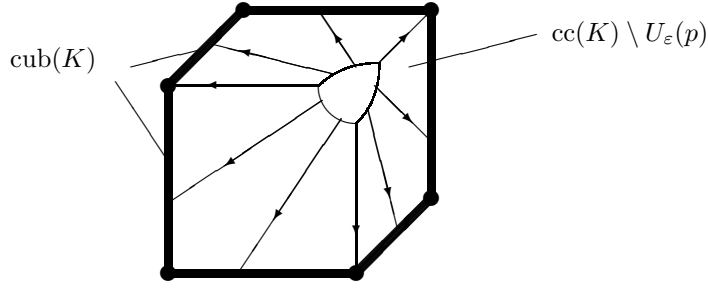
\end{proof}

According to Lemma \ref{macell}, the moment-angle complex $\wk\subset(D^2)^m$
has a cellular structure with 5 different cell types $D_i$, $I_i$, $0_i$,
$T_i$, $1_i$, $i=1,\ldots,m$ (see Figure~7). The homology of $\wk$ (and
therefore of~$W_K$) can be calculated by means of the corresponding cellular
chain complex, which we denote $[\mathcal C_{*}(\wk),\partial_c]$. Though
$\wk$ has more types of cells than $\zk$ (recall that $\zk$ has only 3 cell
types $D_i$, $T_i$,~$1_i$), the cellular chain complex $[\mathcal
C_{*}(\wk),\partial_c]$ can be canonically made {\it bigraded\/} as well.
Namely, the following statement holds (compare with~(\ref{bgcellz})).

\begin{lemma}
\label{bgw}
  Put
  \begin{gather}
  \label{bgcellw}
    \bideg{D_i}=(0,2),\quad\bideg{T_i}=(-1,2),\quad\bideg{I_i}=(1,0),\\
    \bideg{0_i}=\bideg{1_i}=(0,0),\quad i=1,\ldots,m.\notag
  \end{gather}
  This makes the cellular chain complex $[\mathcal C_{*}(\wk),\partial_c]$ a
  bigraded differential module with differential $\partial_c$ adding $(-1,0)$
  to bidegree. The original grading of $\mathcal C_{*}(\wk)$ by dimensions of
  cells corresponds to the total degree (i.e. the dimension of a cell equals
  the sum of its two degrees).
\end{lemma}
\begin{proof}
We need only to check that the differential $\partial_c$ adds $(-1,0)$ to
bidegree. This follows from~(\ref{bgcellw}) and the following formulae:
$$
  \partial_cD_i=T_i,\quad\partial_cI_i=1_i-0_i,\quad
  \partial_cT_i=\partial_c1_i=\partial_c0_i=0.
$$
\end{proof}
Note that unlike bigraded structure in $\mathcal C_{*}(\zk)$ elements of
$\mathcal C_{*,*}(\wk)$ may have {\it positive\/} first degree (due to the
positive first degree of~$I_i$). Nevertheless, the differential
$\partial_c$ does not change the second degree (as in the case of~$\zk$),
which allows to split the bigraded complex $\mathcal C_{*,*}(\wk)$ into the
sum of complexes $\mathcal C_{*,2p}(\wk)$, $p=0,\ldots,m$.

In the same way as we did for $\zk$ and for the pair $(\zk,T^m)$
we define
\begin{gather}
  \label{bbnw}
  b_{q,2p}(\wk)=\dim H_{q,2p}\bigl[\mathcal C_{*,*}(\wk),\partial_c\bigr],
  \quad -m\le q\le m,\;0\le p\le m;\\
  \label{chipw}
  \chi_p(\wk)=\sum_{q=-m}^m(-1)^q\dim\mathcal C_{q,2p}(\wk)
  =\sum_{q=-m}^m(-1)^qb_{q,2p}(\wk);\\
  \chi(\wk;t)=\sum_{p=0}^m\chi_p(\wk)t^{2p}
  \notag
\end{gather}
(note that $q$ above may be both positive and negative).

The following theorem provides the exact formula for the generating polynomial
$\chi(\wk;t)$ and is analogous to theorems~\ref{gpz} and~\ref{rgp}.

\begin{theorem}
\label{chitwk}
  For any simplicial complex $K^{n-1}$ with $m$ vertices holds
  \begin{align*}
    \chi(\wk;t)&=
    (1-t^2)^{m-n}(h_0+h_1t^2+\cdots+h_nt^{2n})+
    \bigl(\chi(K)-1\bigr)(1-t^2)^m\\
    &=(1-t^2)^{m-n}(h_0+h_1t^2+\cdots+h_nt^{2n})+(-1)^{n-1}h_n(1-t^2)^m,
  \end{align*}
  where $\chi(K)=f_0-f_1+\ldots+(-1)^{n-1}f_{n-1}=1+(-1)^{n-1}h_n$ is the
  Euler characteristic of~$K$.
\end{theorem}
\begin{proof}
The definition of $\wk$ (see~(\ref{zkwk})) shows that
$D_II_J0_LT_P1_Q$ (see section~\ref{cell}) is a cell of
$\wk$ if and only if the following two conditions are satisfied:
\begin{enumerate}
\item[(a)] The set $I\cup J\cup L$ is a simplex of~$K^{n-1}$.
\item[(b)] $\#L\ge1$.
\end{enumerate}
Let $c_{ijlpq}(\wk)$ denote the number of cells $D_II_J0_LT_P1_Q\subset\wk$
with $i=\#I$, $j=\#J$, $l=\#L$, $p=\#P$, $q=\#Q$, $i+j+l+p+q=m$. It follows
that
\begin{equation}
\label{ci}
  c_{ijlpq}(\wk)=f_{i+j+l-1}\textstyle\binom{i+j+l}i\binom{j+l}l
  \binom{m-i-j-l}p,
\end{equation}
where $(f_0,\ldots,f_{n-1})$ is the $f$-vector of $K$ (we also assume
$f_{-1}=1$ and $f_k=0$ for $k<-1$ or~$k>n-1$). By (\ref{bgcellw}),
$$
  \bideg(D_II_J0_LT_P1_Q)=(j-p,2(i+p)).
$$
Now we calculate $\chi_r(\wk)$ using (\ref{chipw}) and (\ref{ci}):
$$
  \chi_r(\wk)=\mathop{\sum_{i,j,l,p}}\limits_{i+p=r,l\ge1}(-1)^{j-p}
  f_{i+j+l-1}\textstyle\binom{i+j+l}i\binom{j+l}l\binom{m-i-j-l}p.
$$
Substituting $s=i+j+l$ above we obtain
\begin{align*}
  \chi_r(\wk)&=\mathop{\sum_{l,s,p}}\limits_{l\ge1}(-1)^{s-r-l}
  f_{s-1}\textstyle{\binom{s}{r-p}\binom{s-r+p}l\binom{m-s}p}\\
  &=\sum_{s,p}\Bigl((-1)^{s-r}f_{s-1}{\textstyle\binom{s}{r-p}\binom{m-s}p}
  \sum_{l\ge1}(-1)^l{\textstyle\binom{s-r+p}l}\Bigl)
\end{align*}
Since
$$
  \sum_{l\ge1}(-1)^l{\textstyle\binom{s-r+p}l}=
  \left\{
  \begin{aligned}
    -1,&\quad s>r-p,\\
    0,&\quad s\le r-p
  \end{aligned},
  \right.
$$
we get
\begin{align*}
  \chi_r(\wk)&=-\mathop{\sum_{s,p}}\limits_{s>r-p}(-1)^{s-r}
  f_{s-1}{\textstyle\binom{s}{r-p}\binom{m-s}p}\\
  &=-\sum_{s,p}(-1)^{r-s}
  f_{s-1}{\textstyle\binom{s}{r-p}\binom{m-s}p}+
  \sum_s(-1)^{r-s}f_{s-1}{\textstyle\binom{m-s}{r-s}}.
\end{align*}
The second sum in the above formula is exactly $\chi_r(\zk)$
(see~(\ref{chipzk})). To calculate the first sum we observe that
$\sum_p\bin{s}{r-p}\bin{m-s}p=\bin mr$ (this follows from calculating the
coefficient of $\alpha^r$ in both sides of the identity
$(1+\alpha)^s(1+\alpha)^{m-s}=(1+\alpha)^m$). Hence,
$$
  \chi_r(\wk)=-\sum_s(-1)^{r-s}f_{s-1}{\textstyle\binom mr}+\chi_r(\zk)=
  (-1)^r{\textstyle\binom mr}\bigl( \chi(K)-1 \bigr)+\chi_r(\zk),
$$
since $-\sum_s(-1)^sf_{s-1}=\chi(K)-1$ (remember that $f_{-1}=1$). Finally,
using~(\ref{hchi}) we calculate
\begin{align*}
  \chi(\wk;t)=\sum_{r=0}^m\chi_r(\wk)t^{2r}=
  \sum_{r=0}^m(-1)^r{\textstyle\binom mr}\bigl( \chi(K)-1 \bigr)t^{2r}+
  \sum_{r=0}^m\chi_r(\zk)t^{2r}\\
  =\bigl( \chi(K)-1 \bigr)(1-t^2)^m+
  (1-t^2)^{m-n}(h_0+h_1t^2+\cdots+h_nt^{2n}).
\end{align*}
\end{proof}

Suppose now that $K$ is an orientable simplicial manifold. It is easy
to see that in this case $W_K$ is also orientable. Hence, there
are relative Poincar\'e duality isomorphisms:
\begin{equation}\label{rpd}
  H_k(W_K)\cong H^{m+n-k}(W_K,\partial W_K), \quad k=0,\ldots,m.
\end{equation}

\begin{corollary}[Generalised Dehn--Sommerville equations]
\label{DSsm}
The following relations hold for the $h$-vector $(h_0,h_1,\ldots,h_n)$ of
any orientable simplicial manifold~$K^{n-1}$:
$$
  h_{n-i}-h_i=(-1)^i\bigl(\chi(K^{n-1})-\chi(S^{n-1})\bigr)
  {\textstyle\binom ni},\quad i=0,1,\ldots,n,
$$
where $\chi(S^{n-1})=1+(-1)^{n-1}$ is the Euler characteristic of an
$(n-1)$-sphere.
\end{corollary}
\begin{proof}
By Theorem~\ref{homotwk}, $H_k(W_K)=H_k(\wk)$ and $H^{m+n-k}(W_K,\partial_c
W_K)= H^{m+n-k}(\zk,T^m)$. Moreover, it can be seen in the same way as in
Corollary~\ref{bpd} that relative Poincar\'e duality isomorphisms~(\ref{rpd})
regard the bigraded structures in the (co)homology of $\wk$ and $(\zk,T^m)$.
It follows that
\begin{gather}
  \notag
  b_{-q,2p}(\wk)=b_{-(m-n)+q,2(m-p)}(\zk,T^m),\\
  \notag
  \chi_p(\wk)=(-1)^{m-n}\chi_{m-p}(\zk,T^m),\\
  \label{relchid}
  \chi(\wk;t)=(-1)^{m-n}t^{2m}\chi(\zk,T^m;\textstyle\frac1t).
\end{gather}
Using (\ref{rhchi}) we calculate
\begin{multline*}
  (-1)^{m-n}t^{2m}\chi(\zk,T^m;{\textstyle\frac1t})\\
      =(-1)^{m-n}t^{2m}(1-t^{-2})^{m-n}(h_0+h_1t^{-2}+\cdots+h_nt^{-2n})\\
      -(-1)^{m-n}t^{2m}(1-t^{-2})^m\\
    =(1-t^2)^{m-n}(h_0t^{2n}+h_1t^{2n-2}+\cdots+h_n)+
    (-1)^{n-1}(1-t^2)^m.
\end{multline*}
Substituting the formula for $\chi(\wk;t)$ from Theorem~\ref{chitwk} and the
above expression into~(\ref{relchid}) we obtain
\begin{multline*}
  (1-t^2)^{m-n}(h_0+h_1t^2+\cdots+h_nt^{2n})+\bigl(\chi(K)-1\bigr)
  (1-t^2)^m\\=(1-t^2)^{m-n}(h_0t^{2n}+h_1t^{2n-2}+\cdots+h_n)+
  (-1)^{n-1}(1-t^2)^m.
\end{multline*}
Calculating the coefficient of $t^{2i}$ in both sides after dividing the
above identity by $(1-t^2)^{m-n}$, we obtain
$h_{n-i}-h_i=(-1)^i\bigl(\chi(K^{n-1})-\chi(S^{n-1})\bigr)
{\textstyle\binom ni}$, as needed.
\end{proof}
If $|K|=S^{n-1}$ or $n-1$ is odd, Corollary~\ref{DSsm} gives the classical
equations $h_{n-i}=h_i$.

\begin{corollary}
If $K^{n-1}$ is a simplicial manifold with $h$-vector $(h_0,\ldots,h_n)$ then
$$
  h_{n-i}-h_i=(-1)^i(h_n-1){\textstyle\binom ni},\quad i=0,1,\ldots,n.
$$
\end{corollary}
\begin{proof}
Since
$\chi(K^{n-1})=1+(-1)^{n-1}h_n$, $\chi(S^{n-1})=1+(-1)^{n-1}$, we have
$$
  \chi(K^{n-1})-\chi(S^{n-1})=(-1)^{n-1}(h_n-1)=(h_n-1)
$$
(the coefficient $(-1)^{n-1}$ can be dropped since for odd $n-1$
the left hand side is zero).
\end{proof}

\begin{corollary}
  For any $(n-1)$-dimensional orientable simplicial manifold the numbers
  $h_{n-i}-h_i$, $i=0,1,\ldots,n$, are homotopy invariants. In particular,
  they do not depend on a triangulation.
\end{corollary}
In the particular case of $PL$-manifolds the topological invariance of numbers
$h_{n-i}-h_i$ was firstly observed by Pachner in~\cite[(7.11)]{Pac2}.

\begin{figure}
  \begin{center}
  \begin{picture}(30,30)
  \multiput(0,0)(0,10){4}{\line(1,0){30}}
  \multiput(0,0)(10,0){4}{\line(0,1){30}}
  \put(0,20){\line(1,1){10}}
  \put(0,10){\line(1,1){20}}
  \put(0,0){\line(1,1){30}}
  \put(10,0){\line(1,1){20}}
  \put(20,0){\line(1,1){10}}
  \end{picture}
  \end{center}
  \caption{Triangulation of $T^2$ with $\mb f=(9,27,18)$, $\mb h=(1,6,12,-1)$.}
\end{figure}
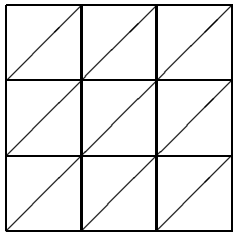
\begin{example}
Consider triangulations of the 2-torus~$T^2$. We have $n=3$, $\chi(T^2)=0$.
From $\chi(K^{n-1})=1+(-1)^{n-1}h_n$ we deduce $h_3=-1$.
Corollary~\ref{DSsm} gives
$$
  h_3-h_0=-2,\quad h_2-h_1=6.
$$
For instance, the triangulation on Figure~10 has $f_0=9$ vertices, $f_1=27$
edges and $f_2=18$ triangles. The corresponding $h$-vector is $(1,6,12,-1)$.
\end{example}

\section{Subspace arrangements and cohomology rings of their complements}
\subsection{Summary of results on the cohomology of general arrangement
complements}
\label{gene}
An {\it arrangement\/} is a finite set $\A=\{L_1,\ldots,L_r\}$ of planes
(affine subspaces) in some affine space (either real or complex). For any
arrangement $\A=\{L_1,\ldots,L_r\}$ in $\C^m$ define its {\it support\/}
$|\A|$ as
$$
  |\A|:=\bigcup_{i=1}^rL_i\subset\C^m,
$$
and its {\it complement\/} $U(\A)$ as
$$
  U(\A):=\C^m\setminus|\A|,
$$
and similarly for arrangements in~$\R^m$.

Arrangements and their complements play a pivotal r\^ole in many
constructions of combinatorics, algebraic and symplectic geometry etc.; they
also arise as configuration spaces for different classical mechanical
systems. In the study of arrangements it is very important to get a
sufficiently detailed description of the topology of complements $U(\A)$
(this includes number of connected components, homotopy type, homology
groups, cohomology ring etc.). A host of elegant results in this direction
appeared during the last three decades, however, the whole picture is far
from being complete. The theory ascends to work of Arnold~\cite{Ar}, which
described the classifying space for the braid group $B_n$ as the complement
of arrangement of all diagonal hyperplanes $\{z_i=z_j\}$ in~$\C^n$. The
cohomology ring of this complement was also calculated there. This result was
generalised by Brieskorn~\cite{Br} and motivated the further development of
the theory of {\it complex hyperplane arrangements\/} (i.e. arrangements of
codimension-one complex affine subspaces). One of the main results here is
the following.

\begin{theorem}[\cite{Ar}, \cite{Br}, \cite{OS}]
 Let $\A=\{L_1,\ldots,L_r\}$ be an arrangement of complex hyperplanes
 in~$\C^m$, and the hyperplane $L_i$ is the zero set of linear function
 $l_i$, $i=1,\ldots,r$. Then the integer cohomology algebra of the complement
 $\C^m\setminus|\A|$ is isomorphic to the algebra generated by closed
 differential 1-forms $\frac 1{2\pi i}\frac{dl_i}{l_i}$.
\end{theorem}
\noindent Relations between the forms $\frac 1{2\pi i}\frac{dl_i}{l_i}$ are
also explicitly described. In the case of diagonal hyperplanes $\{z_i=z_j\}$
we have the forms $\omega_{ij}=\frac1{2\pi i}\frac{d(z_i-z_j)}{z_i-z_j}$,
which are subject to the following relations:
$$
  \omega_{ij}\wedge\omega_{jk}+\omega_{jk}\wedge\omega_{ki}+
  \omega_{ki}\wedge\omega_{ij}=0.
$$

The theory of complex hyperplane arrangements is probably the most well
understood part of the whole study. Several surveys and monographs are
available; we mention just~\cite{OT}, where further references can be found.
Relationships of {\it real\/} hyperplane arrangements with polytopes and {\it
oriented matroids\/} are discussed in~\cite[Lecture~7]{Zi}.

In the general situation, the Goresky--MacPherson theorem~\cite[Part~III]{GM}
expresses the cohomology groups $H^{i}(U(\mathcal A))$ (without ring
structure) as a sum of homology groups of subcomplexes of a certain
simplicial complex. We formulate this result below. For a detailed survey of
general arrangements we refer to~\cite{Bj}. Some important results in this
direction can be found in monograph~\cite{Va}.

Let $\A=\{L_1,\ldots,L_r\}$ be an arrangement of planes in~$\R^n$. The
intersections
$$
  v=L_{i_1}\cap\dots\cap L_{i_k}
$$
form a poset $(\mathcal P,<)$ with respect to the inclusion (i.e. $v<w$ if
and only if $v$ and $w$ are different and $v$ is contained in~$w$). The poset
$\mathcal P$ is assumed to have a unique maximal element $T$ corresponding to
the ambient space of the arrangement. The {\it rank function\/} $d$ on
$\mathcal P$ is defined by $d(v)=\dim v$.  The order complex $K(\P)$ (see
Example~\ref{oc}) is called the {\it order complex of arrangement\/}~$\A$.
Define
$$
  \P_{(v,w)}=\{x\in\P\::\:v<x<w\},\quad\P_{>v}=\{x\in\P\::\:x>v\}.
$$

\begin{theorem}[{Goresky and MacPherson~\cite[Part~III]{GM}}]\label{GMf}
The following formula holds for the homology of the complement~$U(\A)${\rm:}
$$
  H_i\bigl(U(\A);\Z\bigr)=\bigoplus_{v\in\P}
  H^{n-d(v)-i-1}\bigl(K(\P_{>v}),K(\P_{v,T});\Z\bigr),
$$
with the agreement that $H^{-1}(\emptyset,\emptyset)=\Z$.
\end{theorem}
\noindent The proof of this theorem uses the {\it stratified Morse theory\/}
developed in~\cite{GM}.

\begin{remark}
The homology groups of a complex arrangement in $\C^n$ can be
calculated by regarding it as a real arrangement in~$\R^{2n}$.
\end{remark}

The cohomology {\it rings\/} of the complements of arrangements are much more
subtle. In general, the integer cohomology ring of $U(\A)$ {\it is not\/}
determined by the poset~$\P$. An approach to calculating the cohomology
algebra of the complement $U(\A)$ was proposed by De Concini and
Procesi~\cite{dCP}. In particular, they proved that the {\it rational\/}
cohomology ring of $U(\A)$ is determined by the combinatorics of
intersections. This result was extended by Yuzvinsky in~\cite{Yu}.

\subsection{Coordinate subspace arrangements and cohomology of~$\zk$.}
\label{coor}

An arrangement $\A=\{L_1,\ldots,L_r\}$ is called {\it coordinate\/} if every
plane~$L_i$, $i=1,\ldots,r$, is a coordinate subspace. In this section we
apply the results of chapter~4 to cohomology algebras of the complements of
complex coordinate subspace arrangements. The case of real coordinate
arrangements is discussed at the end of the section.

Any coordinate subspace of $\C^m$ has the form
\begin{equation}
\label{li}
  L_I=\{(z_1,\ldots,z_m)\in\C^m\::\:z_{i_1}=\cdots=z_{i_k}=0\},
\end{equation}
where $I=\{i_1,\ldots,i_k\}$ is a subset of~$[m]$. Obviously, $\dim
L_I=m-\#I$.

\begin{construction}\label{casim}
For each simplicial complex $K$ on the set $[m]$ define the
complex coordinate subspace arrangement $\mathcal{CA}(K)$ as the
set of subspaces $L_I$ such that $I$ is not a simplex of~$K$:
$$
  \mathcal{CA}(K)=\{L_I\::\:I\notin K\}.
$$
Denote the complement of $\mathcal{CA}(K)$ by~$U(K)$, that is
\begin{equation}
\label{compl}
  U(K)=\C^m\setminus\bigcup_{I\notin K}L_I.
\end{equation}
If $K'\subset K$ is a subcomplex, then $U(K')\subset U(K)$.
\end{construction}

\begin{proposition}
The assignment $K\mapsto U(K)$ defines a one-to-one order-preserving
correspondence between simplicial complexes on the set~$[m]$ and
complements of coordinate subspace arrangements in~$\C^m$.
\end{proposition}
\begin{proof}
Suppose $\mathcal{CA}$ is a coordinate subspace arrangement in~$\C^m$. Define
\begin{equation}\label{ka}
  K(\mathcal{CA}):=\{I\subset[m]\::\:L_I\not\subset|\mathcal{CA}|\}.
\end{equation}
Obviously, $K(\mathcal{CA})$ is a simplicial complex. By the construction,
$K(\mathcal{CA})$ depends only on $|\mathcal{CA}|$ (i.e. on
$U(\mathcal{CA})$) and $U(K(\mathcal{CA}))=U(\mathcal{CA})$, whence the
proposition follows.
\end{proof}

If a coordinate subspace arrangement $\mathcal A$ contains a hyperplane, say
$\{z_i=0\}$, then its complement $U(\mathcal A)$ is factorised as $U(\mathcal
A_0)\times\C^{*}$, where $\mathcal A_0$ is a coordinate subspace
arrangement in the hyperplane $\{z_i=0\}$ and $\C^{*}=\C\setminus\{0\}$.
Thus, for any coordinate subspace arrangement $\mathcal A$ the complement
$U(\mathcal A)$ decomposes as
$$
  U(\mathcal A)=U(\mathcal A')\times(\C^{*})^k,
$$
were $\mathcal A'$ is a coordinate arrangement in $\C^{m-k}$ that does not
contain hyperplanes. On the other hand,~(\ref{ka}) shows that $\mathcal{CA}$
contains the hyperplane $\{z_i=0\}$ if and only if $\{i\}$ is not a vertex of
$K(\mathcal{CA})$. It follows that $U(K)$ is the complement of a coordinate
arrangement without hyperplanes if and only if the vertex set of $K$ is the
whole~$[m]$. Keeping in mind these remarks, we restrict ourselves to
coordinate subspace arrangements without hyperplanes and simplicial complexes
on the vertex set~$[m]$.

\begin{remark}
  In terms of Construction \ref{nsc} we have $U(K)=K_\bullet(\C,\C^*)$.
\end{remark}

\begin{example}
\label{uk}
1. If $K=\D^{m-1}$ (the $(m-1)$-simplex) then $U(K)=\C^m$.

2. If $K=\partial\D^{m-1}$ (the boundary of simplex) then
$U(K)=\C^m\setminus\{0\}$.

3. If $K$ is a disjoint union of $m$ vertices, then $U(K)$ is obtained
by removing all codimension-two coordinate subspaces
$z_i=z_j=0$, $i,j=1,\ldots,m$, from~$\C^m$.
\end{example}

The action of the algebraic torus $(\C^*)^m$ on $\C^m$ descends to~$U(K)$. In
particular, the standard action of the torus $T^m$ is defined on~$U(K)$. The
quotient $U(K)/T^m$ can be identified with $U(K)\cap\R^m_+$, where $\R^m_+$
is regarded as a subset of~$\C^m$.

\begin{lemma}
\label{zu}
  $\cc(K)\subset U(K)\cap\R^m_+$ and $\zk\subset U(K)$
  (see Construction~{\rm\ref{cck}} and~{\rm(\ref{zkwk})}).
\end{lemma}
\begin{proof}
Take $y=(y_1,\ldots,y_m)\in\cc(K)$. Let $I=\{i_1,\ldots,i_k\}$ be the maximal
subset of $[m]$ such that $y\in L_I\cap\R^n_+$ (i.e.
$y_{i_1}=\dots=y_{i_k}=0$). Then it follows from the definition of $\cc(K)$
(see~(\ref{fcck})) that $I$ is a simplex of~$K$. Hence,
$L_I\notin\mathcal{CA}(K)$ and $y\in U(K)$. Thus, the first statement is
proved. Since $\cc(K)$ is the quotient of~$\zk$, the second assertion follows
from the first one.
\end{proof}

\begin{theorem}
\label{he1}
  There is an equivariant deformation retraction $U(K)\to\zk$.
\end{theorem}
\begin{proof}
First, we construct a deformation retraction $r:U(K)\cap\R^m_+\to\cc(K)$.
This is done inductively. We start from the boundary complex of an
$(m-1)$-simplex and remove simplices of positive dimensions until we
obtain~$K$. On each step we construct a deformation retraction, and the
composite map would be a required retraction~$r$.

If $K=\partial\D^{m-1}$ is the boundary complex of an $(m-1)$-simplex,
then $U(K)\cap\R^m_+=\R^m_+\setminus\{0\}$. In this case the retraction
$r$ is shown on Figure~11.
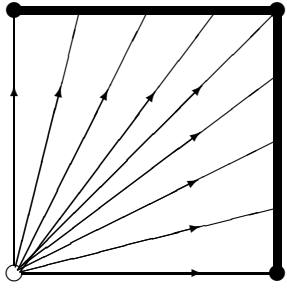
\begin{figure}
  \begin{picture}(120,45)
  \put(45,5){\circle{2}}
  \put(80,5){\circle*{2}}
  \put(45,40){\circle*{2}}
  \put(80,40){\circle*{2}}
  \put(46,5){\vector(1,0){24}}
  \put(46,5){\line(1,0){34}}
  \put(45.8,5.2){\vector(4,1){24}}
  \put(45.8,5.2){\line(4,1){34}}
  \put(45.5,5.5){\vector(2,1){24}}
  \put(45.5,5.5){\line(2,1){34}}
  \put(45.8,5.8){\vector(4,3){24}}
  \put(45.8,5.8){\line(4,3){34}}
  \put(46,6){\vector(1,1){24}}
  \put(46,6){\line(1,1){34}}
  \put(45,6){\vector(0,1){24}}
  \put(45,6){\line(0,1){34}}
  \put(45.2,5.8){\vector(1,4){6}}
  \put(45.2,5.8){\line(1,4){8.5}}
  \put(45.5,5.5){\vector(1,2){12}}
  \put(45.5,5.5){\line(1,2){17}}
  \put(45.8,5.2){\vector(3,4){18}}
  \put(45.8,5.2){\line(3,4){26}}
  \linethickness{1mm}
  \put(80,5){\line(0,1){35}}
  \put(45,40){\line(1,0){35}}
  \end{picture}
  \caption{The retraction $r:U(K)\cap\R^m_+\to\cc(K)$ for
  $K=\partial\D^{m-1}$.}
\end{figure}
Now suppose that $K$ is obtained by removing one $(k-1)$-dimensional simplex
$J=\{j_1,\ldots,j_k\}$ from simplicial complex~$K'$, that is $K\cup J=K'$. By
the inductive hypothesis, the there is a deformation retraction
$r':U(K')\cap\R^m_+\to\cc(K')$. Let $a\in\R^m_+$ be the point with
coordinates $y_{j_1}=\ldots=y_{j_k}=0$ and $y_i=1$ for $i\notin J$. Since $J$
is not a simplex of~$K$, we have $a\notin U(K)\cap\R^m_+$. At the same time,
$a\in C_J$ (see~(\ref{ijface})). Hence, we may apply the retraction from
Figure~11 on the face $C_J\subset I^m$, with centre at~$a$. Denote this
retraction by~$r_J$. Then $r=r_J\circ r'$ is the required deformation
retraction.

The deformation retraction $r:U(K)\cap\R^m_+\to\cc(K)$ is covered by
an equivariant deformation retraction $U(K)\to\zk$.
\end{proof}

In the case $K=K_P$ (i.e. $K$ is a polytopal simplicial sphere corresponding
to a simple polytope~$P^n$) the deformation retraction $U(K_P)\to\zp$ from
Theorem~\ref{he1} can be realised as the orbit map for an action of a
contractible group. We denote $U(P^n):=U(K_P)$. Set
$$
  \R^m_>=\{(y_1,\ldots,y_m)\in\R^m\::\:y_i>0,\; i=1,\ldots,m\}\subset\R^m_+.
$$
Then $\R^m_>$ is a group with respect to the multiplication, and it acts on
$\R^m$, $\C^m$ and $U(P^n)$ by coordinatewise multiplication. There is the
isomorphism $\exp:\R^m\to\R^m_>$ between the additive and the multiplicative
groups taking $(y_1,\ldots,y_m)\in\R^m$ to
$(e^{y_1},\ldots,e^{y_m})\in\R^m_>$.

Let us consider the $m\times(m-n)$-matrix $W$ introduced in
Construction~\ref{dist} for any simple polytope~(\ref{ptope}).

\begin{proposition}
\label{wprop}
  For any vertex $v=F_{i_1}\cap\cdots\cap F_{i_n}$ of $P^n$ the
  maximal minor $W_{\hat i_1\ldots\hat i_n}$ which is obtained by deleting
  $n$ rows $i_1,\ldots,i_n$ from $W$ is non-degenerate: $\det W_{\hat
  i_1\ldots \hat i_n}\ne0$.
\end{proposition}
\begin{proof}
If $\det W_{\hat i_1\ldots \hat i_n}=0$ then vectors
$\mb l_{i_1},\ldots,\mb l_{i_n}$ (see~(\ref{ptope})) are linearly dependent,
which is impossible.
\end{proof}

The matrix $W$ defines the subgroup
\begin{equation}\label{rw}
  R_W=\bigr\{(e^{w_{11}\tau_1+\cdots+w_{1,m-n}\tau_{m-n}},\ldots,
  e^{w_{m1}\tau_1+\cdots+w_{m,m-n}\tau_{m-n}})\bigl\}\subset\R^m_>,
\end{equation}
where the parameters $\tau_1,\ldots,\tau_{m-n}$ vary over~$\R^{m-n}$.
Obviously, $R_W\cong\R^{m-n}_>$.

\begin{theorem}[{\cite[Theorem~2.3]{BP4}}]
\label{zpos}
The subgroup $R_W\subset\R^m_>$ acts freely on $U(P^n)\subset\C^m$. The
composite map $\zp\hookrightarrow U(P^n)\to U(P^n)/R_W$ of the embedding
$i_e$ (Lemma~{\rm\ref{ie}}) and the orbit map is an equivariant diffeomorphism
(with respect to the $T^m$-actions).
\end{theorem}
\begin{proof}
A point from $\C^m$ has the non-trivial isotropy subgroup with respect to the
action of $\R^m_>$ on $\C^m$ if and only if at least one of its coordinates
vanishes.  It follows from~(\ref{compl}) that if a point $x\in U(P^n)$ has
some zero coordinates, then the corresponding facets of $P^n$ have at least
one common vertex $v\in P^n$. Let $v=F_{i_1}\cap\cdots\cap
F_{i_n}$. The isotropy subgroup of the point $x$ with respect to the
action of the subgroup $R_W$ is non-trivial only if some linear combination
of columns of $W$ lies in the coordinate subspace spanned by $\mb
e_{i_1},\ldots,\mb e_{i_n}$. But this implies that $\det W_{\hat i_1\ldots
\hat i_n}=0$, which contradicts Proposition~\ref{wprop}.  Thus, $R_W$ acts
freely on $U(P^n)$.

To prove the second part of the theorem it is sufficient to show that each
orbit of the action of $R_W$ on $U(P^n)\subset\C^m$ intersects the image
$i_e(\zp)$ at a single point. Since the embedding $i_e$ is equivariant with
respect to the $T^m$-actions, the latter statement is equivalent to that each
orbit of the action of $R_W$ on $U(P^n)\cap\R^m_+$ intersects the image
$i_P(P^n)$ (see Theorem~\ref{thcubpol}) at a single point. Let $y\in
i_P(P^n)\subset\R^m$. Then $y=(y_1,\ldots,y_m)$ lies in some $n$-face
$i_P(C^n_v)$ of the cube $I^m\subset\R^m$, see~(\ref{cubpolmap}). We need to
show that the $(m-n)$-dimensional subspace spanned by the vectors
$(w_{11}y_1,\ldots,w_{m1}y_m)^t,\ldots,(w_{1,m-n}y_1,\ldots,w_{m,m-n}y_m)^t$
is in general position with the $n$-face $i_P(C^n_v)$ of~$I^m$. This
follows directly from~(\ref{cubpolmap}) and Proposition~\ref{wprop}.
\end{proof}

Suppose now that $P^n$ is a lattice simple polytope, and let $M_P$ be the
corresponding toric variety (Construction~\ref{nf}). Along with the real
subgroup $R_W\subset\R^m_>$~(\ref{rw}) we define
$$
  C_W=\bigl\{(e^{w_{11}\phi_1+\cdots+w_{1,m-n}\phi_{m-n}},\ldots,
  e^{w_{m1}\phi_1+\cdots+w_{m,m-n}\phi_{m-n}})\bigr\}
  \subset (\C^{*})^m,
$$
where the parameters $\phi_1,\ldots,\phi_{m-n}$ vary over~$\C^{m-n}$.
Obviously, $C_W\cong(\C^*)^{m-n}$. It is shown in \cite{Au},
\cite{Ba},~\cite{Co2} that $C_W$ acts freely on $U(P^n)$ and the toric
variety $M_P$ can be identified with the orbit space (or {\it geometric
quotient\/}) $U(P^n)/C_W$. Thus, we have the following commutative diagram:
\begin{equation}\label{uptv}
\begin{CD}
  U(P^n) @>R_W\cong\R^{m-n}_{>}>> \zp\\
  @VC_W\cong(\C^{*})^{m-n}VV @VVT^{m-n}V\\
  M_P @= M_P.
\end{CD}
\end{equation}

\begin{remark}
It can be shown~\cite[Theorem~2.1]{Co2} that {\it any\/} toric variety
$M_\Sigma$ corresponding to a fan $\Sigma\subset\R^n$ with $m$
one-dimensional cones can be identified with the universal categorical
quotient $U(\mathcal{CA}_\Sigma)/G$, where $U(\mathcal{CA}_\Sigma)$ is the
complement of a certain coordinate arrangement (determined by the
fan~$\Sigma$) and $G\cong(\C^*)^{m-n}$. The categorical quotient becomes the
geometric quotient if and only if the fan $\Sigma$ is simplicial. In this
case $U(\mathcal{CA}_\Sigma)=U(K_\Sigma)$.
\end{remark}

On the other hand, if the projective toric variety $M_P$ is non-singular,
then $M_P$ is a symplectic manifold of dimension~$2n$, and the action of
$T^n$ on it is Hamiltonian~\cite{Au}. In this case the diagram~(\ref{uptv})
displays $M_P$ as the result of the process of {\it symplectic reduction\/}.
Namely, let $H_W\cong T^{m-n}$ be the maximal compact subgroup in~$C_W$, and
$\mu:\C^m\to\R^{m-n}$ the {\it moment map\/} for the Hamiltonian action of
$H_W$ on~$\C^m$. Then for any regular value $a\in\R^{m-n}$ of the map $\mu$
there is the following diffeomorphism:
$$
  \mu^{-1}(a)/H_W\longrightarrow U(P^n)/C_W=M_P
$$
(details can be found in~\cite{Au}). In this situation $\mu^{-1}(a)$ is
exactly our manifold~$\zp$. This gives us another interpretation of the
manifold $\zp$ as the level surface for the moment map (in the case when
$P^n$ can be realised as the quotient of a non-singular projective toric
variety).

\begin{example}
Let $P^n=\D^n$ (the $n$-simplex). Then $m=n+1$,
$U(P^n)=\C^{n+1}\setminus\{0\}$. $R_W\cong\R_>$, $C_W\cong\C^{*}$ and
$H_W\cong S^1$ are the diagonal subgroups in $\R^{n+1}_>$, $(\C^{*})^{n+1}$
and $T^{m+1}$ respectively (see Example~\ref{simdist}). Hence, $\zp\cong
S^{2n+1}=(\C^{n+1}\setminus\{0\})/\R_>$ and
$M_P=(\C^{n+1}\setminus\{0\})/\C^{*}=\C P^n$. The moment map $\mu:\C^m\to\R$
takes $(z_1,\ldots,z_m)\in\C^m$ to $\frac12(|z_1|^2+\ldots+|z_m|^2)$, and for
$a\ne0$ we have $\mu^{-1}(a)\cong S^{2n+1}\cong\zk$.
\end{example}
The previous discussion illustrates the importance of calculating the
cohomology of subspace arrangement complements.

\begin{theorem}[Buchstaber and Panov]
\label{cohar}
  The following isomorphism of graded algebras holds:
  \begin{align*}
    H^{*}\bigl(U(K)\bigr)
    &\cong\Tor_{\k[v_1,\ldots,v_m]}\bigl(\k(K),\k\bigr)\\
    &=H\bigl[\L[u_1,\ldots,u_m]\otimes\k(K),d\bigr].
  \end{align*}
\end{theorem}
\begin{proof}
This follows from theorems~\ref{he1}, \ref{cohom1} and \ref{cohom2}.
\end{proof}

Theorem~\ref{cohar} provides an extremely effective way to calculate the
cohomology algebra of the complement of any complex coordinate subspace
arrangement. The De Concini and Procesi~\cite{dCP} and Yuzvinsky~\cite{Yu}
rational models of the cohomology algebra of an arrangement complement also
can be interpreted as an application of the Koszul resolution. However, these
author did not discuss the relationships with the Stanley--Reisner ring in
the case of coordinate subspace arrangements.

\begin{problem}
Calculate the cohomology algebra {\sl with $\Z$ coefficients\/} of a
coordinate subspace arrangement complement and describe its relationships
with the corresponding $\Tor$-algebra $\Tor_{\Z[v_1,\ldots,v_m]}(\Z(K),\Z)$.
\end{problem}

\begin{example}
Let $K$ be a disjoint union of $m$ vertices.  Then $U(K)$ is obtained by
removing all codimension-two coordinate subspaces $z_i=z_j=0$,
$i,j=1,\ldots,m$ from $\C^m$ (see Example~\ref{uk}).  The face ring is
$\k(K)=\k[v_1,\ldots,v_m]/\mathcal I_K$, where $\mathcal I_K$ is generated by
monomials $v_iv_j$, $i\ne j$.  An easy calculation using
Corollary~\ref{cohar} shows that the subspace of cocycles in
$\k(K)\otimes\Lambda[u_1,\ldots,u_m]$ has the basis consisting of monomials
$v_{i_1}u_{i_2}u_{i_3}\cdots u_{i_k}$ with $k\ge2$, $i_p\ne i_q$ for
$p\ne q$. Since $\deg(v_{i_1}u_{i_2}u_{i_3}\cdots u_{i_k})=k+1$, the space of
$(k+1)$-dimensional cocycles has dimension $m\binom{m-1}{k-1}$.  The space of
$(k+1)$-dimensional coboundaries is $\binom mk$-dimensional (it is spanned by
the coboundaries of the form $d(u_{i_1}\cdots u_{i_k})$).  Hence,
\begin{align*}
  &\dim H^{0}\bigl(U(K)\bigr)=1,\quad
  H^{1}\bigl(U(K)\bigr)=H^{2}\bigl(U(K)\bigr)=0,\\
  &\dim H^{k+1}\bigl(U(K)\bigr)=
  m\bin{m-1}{k-1}-\bin mk=(k-1)\bin mk,\quad2\le k\le m,
\end{align*}
and the multiplication in the cohomology is trivial.

In particular, for $m=3$ we have 6 three-dimensional cohomology classes
$[v_iu_j]$, $i\ne j$ subject to 3 relations $[v_iu_j]=[v_ju_i]$, and 3
four-dimensional cohomology classes $[v_1u_2u_3]$, $[v_2u_1u_3]$,
$[v_3u_1u_2]$ subject to one relation
$$
  [v_1u_2u_3]-[v_2u_1u_3]+[v_3u_1u_2]=0.
$$
Hence, $\dim H^{3}(U(K))=3$, $\dim H^4(U(K))=2$,
and the multiplication is trivial.
\end{example}

\begin{example}
  Let $K$ be the boundary of an $m$-gon, $m>3$. Then
  $$
    U(K)=\C^m\setminus\bigcup_{i-j\ne0,1\mod m}\{z_i=z_j=0\}.
  $$
  By Theorem~\ref{cohar}, the cohomology ring of $H^*(U(K);\k)$ is isomorphic
  to the ring described in Example~\ref{mgon}.
\end{example}

As it is shown in~\cite{GPW}, in the case of arrangements of {\it real\/}
coordinate subspaces only {\it additive\/} analogue of our
Theorem~\ref{cohar} holds. Namely, let us consider the polynomial ring
$\k[x_1,\ldots,x_m]$ with $\deg x_i=1$, $i=1,\ldots,m$. Then the graded
structure in the face ring $\k(K)$ changes accordingly. The Betti numbers of
the real coordinate subspace arrangement $U_\R(K)$ are calculated by means of
the following result.

\begin{theorem}[{\cite[Theorem~3.1]{GPW}}]
The following isomorphism hold:
$$
  H^p\bigl(U_\R(K)\bigr)
  \cong\sum_{-i+j=p}\Tor^{-i,j}_{\k[x_1,\ldots,x_m]}\bigl(\k(K),\k\bigr)
  =H^{-i,j}\bigl[\L[u_1,\ldots,u_m]\otimes\k(K),d\bigr],
$$
where $\bideg u_i=(-1,1)$, $\bideg v_i=(0,1)$, $du_i=x_i$, $dx_i=0$.
\end{theorem}

As it was observed in~\cite{GPW}, there is {\it no\/} multiplicative
isomorphism analogous to Theorem~\ref{cohar} in the case of real
arrangements, that is, the algebras $H^*(U_\R(K))$ and
$\Tor_{\k[x_1,\ldots,x_m]}(\k(K),\k)$ are not isomorphic in general. The
paper~\cite{GPW} also contains the formulation of the first multiplicative
isomorphism of our Theorem~\ref{cohar} for complex coordinate subspace
arrangements (see~\cite[Theorem~3.6]{GPW}), with a reference to yet
unpublished paper by Babson and Chan.

Until now, we have used the description of coordinate subspaces by means of
equations (see~(\ref{li})). On the other hand, a coordinate subspace can be
defined as the linear span of some subset of the standard basis
$\{\mb e_1,\ldots,\mb e_m\}$. This leads to the dual approach to the
description of coordinate subspace arrangements, which corresponds to the
passage from simplicial complex $K$ to the associated complex $\widehat{K}$
(Example~\ref{dual}). This approach was used in~\cite{dL}. It was shown
there that the summands in the Goresky--MacPherson formula in the coordinate
subspace arrangement case are homology groups of links of simplices
of~$\widehat{K}$. This allowed to interpret the product of cohomology classes
of the complement of a coordinate subspace arrangement (either real or
complex) in terms of the combinatorics of links of simplices in~$\widehat{K}$
(see~\cite[Theorem~1.1]{dL}).

We mention that our theorems~\ref{cohom2} and~\ref{he1} show that the
Goresky--\-Mac\-Pher\-son result (Theorem~\ref{GMf}) in the case of coordinate
subspace arrangements is equivalent to the Hochster theorem
(Theorem~\ref{hoch}).

\subsection{Diagonal subspace arrangements and cohomology of loop
space~$\O\zk$.}
\label{diag}
In this section we establish relationships between the results of~\cite{PRW}
on the cohomology of real diagonal arrangement complements and the cohomology
of the loop spaces $\O\bk$ and~$\O\zk$.

For each subset $I=\{i_1,\ldots,i_k\}\subset[m]$
define the {\it diagonal subspace\/} $D_I$ in $\R^m$ as
$$
  D_I=\{(y_1,\ldots,y_m)\in\R^m\::\:y_{i_1}=\cdots=y_{i_k}\}.
$$
Diagonal subspaces in $\C^m$ are defined similarly. An arrangement of planes
$\A=\{L_1,\ldots,L_r\}$ (either real or complex) is called {\it diagonal\/}
if all planes $L_i$, $i=1,\ldots,r$, are diagonal subspaces. The classical
example of a diagonal subspace arrangement is given by the arrangement of all
diagonal hyperplanes $\{z_i=z_j\}$ in~$\C^m$; its complement is the
classifying space for the braid group~$B_m$, see~\cite{Ar}.

\begin{construction}\label{dasim}
Given a simplicial complex $K$ on the vertex set~$[m]$, introduce the
diagonal subspace arrangement $\mathcal{DA}(K)$ as the set of subspaces
$D_I$ such that $I$ is not a simplex of~$K$:
$$
  \mathcal{DA}(K)=\{D_I\::\:I\notin K\}.
$$
Denote the complement of the arrangement $\mathcal{DA}(K)$ by~$M(K)$.
\end{construction}

The proof of the following statement is similar to the proof of the
corresponding statement for coordinate subspace arrangements from
section~\ref{coor}.

\begin{proposition}
The assignment $K\mapsto M(K)$ defines a one-to-one order-preserving
correspondence between simplicial complexes on the vertex set~$[m]$ and the
complements of diagonal subspace arrangements in~$\R^m$.
\end{proposition}

Here we still assume that $\k$ is a field. The multigraded (or $\N^m$-graded)
structure in the ring $\k[v_1,\ldots,v_m]$ (Construction~\ref{mgrad}) defines
an $\N^m$-grading in the Stanley--Reisner ring~$\k(K)$. The monomial
$v_1^{i_1}\cdots v_m^{i_m}$ acquires the multidegree $(2i_1,\ldots,2i_m)$.
Let us consider the modules $\Tor_{\k(K)}(\k,\k)$. They can be calculated,
for example, by means of the minimal free resolution (Example~\ref{minimal})
of the field~$\k$ regarded as a $\k(K)$-module. The minimal resolution also
carries a natural $\N^m$-grading, and we denote the subgroup of elements of
multidegree $(2i_1,\ldots,2i_m)$ in $\Tor_{\k(K)}(\k,\k)$ by
$\Tor_{\k(K)}(\k,\k)_{(2i_1,\ldots,2i_m)}$.

\begin{theorem}[{\cite[Theorem 1.3]{PRW}}]\label{dacoh}
The following isomorphism holds for the cohomology groups of the complement
$M(K)$ of a real diagonal subspace arrangement:
$$
  H^i\bigl(M(K);\k\bigr)\cong\Tor^{-(m-i)}_{\k(K)}(\k,\k)_{(2,\ldots,2)}.
$$
\end{theorem}

\begin{remark}
Instead of simplicial complexes $K$ on the vertex set~$[m]$ the authors
of~\cite{PRW} considered square-free monomial ideals $\mathcal
I\subset\k[v_1,\ldots,v_m]$. Proposition~\ref{sfmi} shows that these two
approaches are equivalent.
\end{remark}

\begin{theorem}\label{looptor}
The following additive isomorphism holds:
$$
  H^*(\O\bk;\k)\cong\Tor_{\k(K)}(\k,\k).
$$
\end{theorem}
\begin{proof}
Let us consider the Eilenberg--Moore spectral sequence of the Serre fibration
$P\to SR(K)$ with fibre $\O SR(K)$, where $SR(K)$ is the Stanley--Reisner
space (Definition~\ref{srspace}) and $P$ is the path space over~$SR(K)$. By
Corollary~\ref{onefib},
\begin{equation}\label{de2t}
  E_2=\Tor_{H^*(SR(K))}\bigl(H^*(P),\k\bigr)\cong\Tor_{\k(K)}(\k,\k),
\end{equation}
and the spectral sequence converges to $\Tor_{C^*(SR(K))}(C^*(P),\k)\cong
H^*(\O SR(K))$. Since $P$ is contractible, there is a cochain equivalence
$C^*(P)\simeq\k$. We have $C^*(SR(K))\cong\k(K)$. Therefore,
$$
  \Tor_{C^*(SR(K))}\bigl(C^*(P),\k\bigr)\cong\Tor_{\k(K)}(\k,\k),
$$
which together with (\ref{de2t}) shows that the spectral sequence collapses
ate the $E_2$ term. Hence, $H^*(\O SR(K))\cong\Tor_{\k(K)}(\k,\k)$. Finally,
Theorem~\ref{homeq1} shows that $H^*(\O SR(K))\cong H^*(\O\bk)$, which
concludes the proof.
\end{proof}

\begin{proposition}\label{loopi}
The following isomorphism of algebras holds
$$
  H^*(\O\bk)\cong H^*(\O\zk)\otimes\L[u_1,\ldots,u_m].
$$
\end{proposition}
\begin{proof}
Consider the bundle $\bk\to BT^m$ with fibre~$\zk$. It is not hard to prove
that the corresponding loop bundle $\O\bk\to T^m$ with fibre $\O\zk$ is
trivial (note that $\O BT^m\cong T^m$). To finish the proof it remains to
mention that $H^*(T^m)\cong\L[u_1,\ldots,u_m]$.
\end{proof}

Theorems \ref{he1} and \ref{cohar} give an application of the theory of
moment-angle complexes to calculating the cohomology ring of a coordinate
subspace arrangement complement. Similarly,
theorems~\ref{dacoh},~\ref{looptor} and Proposition~\ref{loopi} establish the
connection between the cohomology of a diagonal subspace arrangement
complement and the cohomology of the loop space over the moment-angle
complex~$\zk$. However, in this case the situation is much more subtle than
in the case of coordinate subspace arrangements. For instance, we do not have
an analogue of the {\it multiplicative\/} isomorphism from
Theorem~\ref{cohar}. That is why we consider this concluding section only as
the first step in applying the theory of moment-angle complexes to studying
the complements of diagonal subspace arrangements.

\end{document}